\documentclass[12pt,a4paper,fleqn]{article}
\pdfoutput=1

\usepackage[english]{babel}
\RequirePackage{ucs}
\RequirePackage[utf8x]{inputenc}
\RequirePackage[T1]{fontenc}

\RequirePackage{amsthm,amsmath}
\RequirePackage{amssymb,amstext}
\RequirePackage{amsfonts,stmaryrd}
\RequirePackage{amsbsy}
\RequirePackage{mathrsfs}
\RequirePackage{bm}
\RequirePackage{bbm}
\RequirePackage{aliascnt,ifthen}
\RequirePackage{graphicx}
\graphicspath{{./img/}}
\RequirePackage{color,calc}
\RequirePackage{times}
\usepackage{multirow}
\usepackage{supertabular}
\usepackage{booktabs}
\usepackage{caption}
\usepackage{subcaption}
\usepackage{algorithm}
\usepackage[noend]{algpseudocode}

\makeatletter
\def\BState{\State\hskip-\ALG@thistlm}
\makeatother

\usepackage{url}
\usepackage{verbatim}

\usepackage[pdfpagemode=UseOutlines ,plainpages=false,
hyperindex=true,colorlinks=true]{hyperref}
	\makeatletter
	\Hy@breaklinkstrue
	\makeatother

\definecolor{darkred}{rgb}{0.6,0,0.1}
\definecolor{darkgreen}{rgb}{0,0.5,0}
\definecolor{darkblue}{rgb}{0,0,0.5}

\hypersetup{colorlinks ,linkcolor=gray ,filecolor=darkgreen, urlcolor=darkblue ,citecolor=black}

\RequirePackage{paralist}
\usepackage{enumitem}

\usepackage{natbib}
\usepackage{hypernat}
\renewcommand{\cite}{\citet}
\usepackage{DACD-abrev-package}

\definecolor{dgreen}{rgb}{0,0.5,0}
\definecolor{dblue}{rgb}{0,0,0.5}
\definecolor{dred}{rgb}{0.6,0.0,0.1}
\definecolor{dgold}{rgb}{0.5,0.3,0.0}
\definecolor{dvio}{rgb}{0.6,0.3,0.5}
\definecolor{gray}{rgb}{0.5,0.5,0.5}
\definecolor{dbraun}{rgb}{.5,0.2,0}

\newcommand{\dr}{\color{dred}}
\newcommand{\db}{\color{dblue}}
\newcommand{\dg}{\color{dgreen}}

\newcommand{\dgrau}{\color{gray}}

\newcommand{\colre}{dred}
\newcommand{\colas}{dblue}
\newcommand{\colde}{dblue}
\newcommand{\colrem}{dgold}
\newcommand{\colil}{dgreen}

\newtheoremstyle{styre}
  {1.1\topsep}
  {\topsep}
  {\normalfont\itshape}
  {}
  {\color{\colre}}
  {.}
  {.5em}
  {\thmname{\textbf{#1}\xspace}{\xspace\dgrau\thmnumber{#2}}\thmnote{\xspace\textit{\small(#3)}}}

\newtheoremstyle{styas}
  {1.1\topsep}
  {\topsep}
  {\normalfont\itshape}
  {}
  {\color{\colas}}
  {.}
  {.5em}
  {}

\newtheoremstyle{styrem}
  {1.1\topsep}
  {\topsep}
  {\normalfont\itshape}
  {}
  {\color{\colrem}}
  {.}
  {.5em}
  {}

\newtheoremstyle{styil}
  {1.1\topsep}
  {\topsep}
  {\normalfont\rmfamily}
  {}
  {\color{\colil}}
  {.}
  {.5em}
  {\thmname{\textbf{#1}\xspace}{\xspace\db\thmnumber{#2}}\thmnote{\xspace\textit{\small(#3)}}}

\newtheoremstyle{styte}
        {}
        {}
        {\normalfont}
        {}
        {}
        {}
        {0ex}
        {}

\newtheoremstyle{stypro}
	{0.5\topsep}
	{1.1\topsep}
	{\upshape}
	{}
	{\color{\colre}}
	{}
	{.5em}
	{\thmnote{\textit{\textcolor{\colre}#3}}}

\theoremstyle{styre}\newtheorem{pr}{Proposition}[section]
\newaliascnt{co}{pr}
\theoremstyle{styre}\newtheorem{co}[co]{Corollary}
\aliascntresetthe{co}
\newaliascnt{thm}{pr}
\theoremstyle{styre}\newtheorem{thm}[thm]{Theorem}
\aliascntresetthe{thm}
\newaliascnt{lem}{pr}
\theoremstyle{styre}\newtheorem{lem}[lem]{Lemma}
\aliascntresetthe{lem}
\newaliascnt{rem}{pr}
\theoremstyle{styrem}\newtheorem{rem}[rem]{Remark}
\aliascntresetthe{rem}
\newaliascnt{il}{pr}
\theoremstyle{styil}\newtheorem{il}[il]{Illustration}
\aliascntresetthe{il}
\theoremstyle{styas}\newtheorem{ass}{Assumption}

\theoremstyle{styas}

\theoremstyle{stypro}\newtheorem*{pro}{}

\newcommand{\remEnd}{{\scriptsize\textcolor{\colrem}{\qed}}}
\newcommand{\proEnd}{{\scriptsize\textcolor{\colre}{\qed}}}
\newcommand{\ilEnd}{{\scriptsize\textcolor{\colil}{\qed}}}

\usepackage{cleveref}
\crefname{pr}{\color{\colre}Proposition}{\color{\colre}Propositions}
\crefname{co}{\color{\colre}Corollary}{\color{\colre}Corollaries}
\crefname{thm}{\color{\colre}Theorem}{\color{\colre}Theorems}
\crefname{lem}{\color{\colre}Lemma}{\color{\colre}Lemmata}

\crefname{ass}{assumption}{assumptions}
\crefformat{ass}{#2{\color{\colas}Assumption {#1}}#3}

\crefname{de}{definition}{definitions}
\crefformat{de}{#2{\color{\colas}Definition {#1}}#3}

\crefname{rem}{{\color{\colrem}Remark}}{{\color{\colrem}Remarks}}

\crefname{il}{\color{\colil}Illustration}{\color{\colil}Illustrations}

\usepackage{environ}
\NewEnviron{te}{
{\BODY}
}

\makeatletter
\newcommand{\mylabel}[2]{#2\def\@currentlabel{#2}\label{#1}}
\makeatother

\numberwithin{equation}{section}

\newcounter{nc}
\makeatletter
\def\@fnsymbol#1{\ensuremath{\ifcase#1\or 1 \or 2 \or 3 \or  *\or  \star \or 4\or  , \or
g\or h\or i\else\@ctrerr\fi}}
\makeatother

\author{\begin{minipage}{.28\textwidth}\center {\sc Jan Johannes}\thanks{Institut f\"ur Angewandte
    Mathematik, M$\Lambda$THEM$\Lambda$TIKON, Im Neuenheimer Feld 205,
  D-69120 Heidelberg, Germany, e-mail: \url{johannes@math.uni-heidelberg.de}}\\[.5ex]\small Ruprecht-Karls-Universität Heidelberg\\\null\end{minipage}\and
\begin{minipage}{.28\textwidth}\center
{\sc Xavier Loizeau}\thanks{National Centre of Excellence in Mass Spectrometry Imaging (NiCE-MSI), National Physical Laboratory (NPL), Hampton Road, Teddington TW11 0LW, Royaume-Uni,    e-mail: \url{xavier.loizeau@npl.co.uk}}\\[.5ex]\small
 National Physical Laboratory (NPL)\\\null\end{minipage}\\[-6ex]}

\date{Preliminary version: \today}
\title{Data-driven aggregation in circular deconvolution}

\begin{document}
\maketitle

\begin{abstract}
  In a circular deconvolution model we consider the fully data driven density estimation of a circular random variable where the density of the additive independent measurement error is unknown.
  We have at hand two independent iid samples, one of the contaminated version of the variable of interest, and the other of the additive noise.
  We show optimality, in an oracle and minimax sense, of a fully data-driven weighted sum of orthogonal series density estimators.
  Two shapes of random weights are considered, one motivated by a Bayesian approach and the other by a well known model selection method.
  We derive non-asymptotic upper bounds for the quadratic risk and the maximal quadratic risk over Sobolev-like ellipsoids of the fully data-driven estimator.
  We compute rates which can be obtained in different configurations for the smoothness of the density of interest and the error density.
  The rates (strictly) match the optimal oracle or minimax rates for a large variety of cases, and feature otherwise at most a deterioration by a logarithmic factor.
  We illustrate the performance of the fully data-driven weighted sum of orthogonal series estimators by a simulation study.
\end{abstract}
{\footnotesize
\begin{tabbing}
\>Aggregation, Oracle inequality, Adaptation\\[.2ex]
\end{tabbing}}

\section{Introduction}\label{s:i}

\begin{te}
  In a circular convolution model one objective is to estimate non-parametrically the density of a random variable taking values on the unit circle from observations blurred by an additive noise.
  Here we show optimality, in an oracle and minimax sense, of a fully data-driven weighted sum of orthogonal series estimators (OSE's).
  Two shapes of random weights are considered, one motivated by a Bayesian approach and the other by a well known model selection method.
  Circular data are met in a variety of applications, such as data representing a direction on a compass in handwriting recognition (\cite{bahlmann2006directional}) and in meteorology (\cite{carnicero2013non}), or anything from  opinions on a political compass to time reading on a clock face (\cite{gill2010circular}) in political sciences.
  The non-parametric density estimation in a circular deconvolution model has been considered for example in \cite{Efromovich97, ComteTaupin, JohannesSchwarz2013a}, while \cite{SchluttenhoferJohannes2020a, SchluttenhoferJohannes2020b}, for example, study minimax testing.
  For an overview of convolutional phenomenons met in other models the reader may refer to \cite{Meister2009}.
\end{te}

\begin{te}
  Throughout this work we will tacitly identify the circle with the unit interval $[0,1)$, for notational convenience.
  Let $\rY := \rX + \rE - \lfloor \rX + \rE \rfloor = \rX + \rE \mod 1$ be the observable contaminated random variable and $\ydf$ its density.
  If we denote by $\xdf$ and $\edf$ the respective circular densities of the random variable of interest $\rX$ and of the additive and independent noise $\rE$, then, we have
  \begin{equation*}
    \ydf(y) = (\xdf \oast \edf)(y):= \int_{[0,1)} \xdf((y-s) - \gauss{y-s})\,\edf(s)\,ds,\quad y\in[0,1),
  \end{equation*}
such that $\oast$ stands for the circular convolution.
  Therefore, the estimation of $\xdf$ is called a circular deconvolution problem.
\end{te}

\begin{te}
  We highlight hereafter that, thanks to the convolution theorem, an estimator of the circular density $\xdf$ is usually based on the Fourier transforms of $\edf$, and $\ydf$ which may be estimated from the data.
  For any complex number $z$, denote $\overline{z}$ its complex conjugate, and $\vert z \vert$ its modulus.
  Let $\Lp[2]:=\Lp[2]([0,1))$ be the Hilbert space of square integrable complex-valued functions defined on $[0,1)$ endowed with the usual inner product $\VskalarLp{h_{1},h_{2}}=\int_{[0,1)} h_{1}(x) \overline{ h_{2}(x)}dx$, and associated norm $\VnormLp{\mbullet}$.
  Each $h \in \Lp[2]$ admits a representation as discrete Fourier series $\He = \sum_{j \in \Zz} \fHe{j} \bas_j$ with respect to the exponential basis $\{\bas_j\}_{j\in\Zz}$, where $\fHe{j}:=\VskalarLp{h,\bas_j}$ is the $j$-th Fourier coefficient of $\He$, and $\bas_j(x): = \exp(-\iota2\pi j x)$ for $x\in[0,1)$, and a square root $\iota$ of $-1$.
\end{te}

\begin{te}
  In this work we suppose that $\xdf$, $\edf$, and hence $\ydf$, belong to the subset $\cD$ of real-valued Lebesgue densities in $\Lp[2]$.
  We denote the expectation associated with $\ydf$ and $\edf$ by $\Ex[\ydf]$, and $\Ex[\edf]$ respectively.
  We note that $\fydf[0] = 1$, and $\Ex[\ydf][\bas_{j}(-\rY)] = \fydf[j] = \ofou[-j]{\ydf}$ for any $j \in \Zz$ as it is the case for any density.
  The key to our analysis is the convolution theorem which states that, in a circular model, $\ydf = \edf \oast \xdf$ holds if and only if $\fou[j]{\ydf} = \fou[j]{\edf} \cdot \fou[j]{\xdf}$ for all $j\in\Zz$.
  Therefore and as long as $\fedf[j] \ne 0$ for all $j\in\Zz$, which is assumed from now on, we have
  \begin{equation} \label{i:rep}
    \xdf = \bas_0 + \sum\nolimits_{|j|\in\Nz} \fedfI[j]\,\fydf[j]\, \bas_j\; \quad \mbox{ with } \fydf[j]=\Ex[\ydf][\bas_j(-\rY)]\mbox{ and } \fedf[j]=\Ex[\edf][\bas_j(-\rE)].
  \end{equation}
  Note that an analogous representation holds in the case of deconvolution on the real line with compactly supported $\rX$-density, i.e. when the error term $\rE$, and hence $\rY$, take their values in $\Rz$.
  In this situation, the deconvolution density still admits a discrete representation as in \eqref{i:rep}, but involving the characteristic functions of $\edf$ and $\ydf$ rather than their discrete Fourier coefficients.
  For a more detailed study of the Fourier analysis of probability distributions, the reader is referred, for example, to \citet{bremaud2014fourier}, Chapter 2.
\end{te}

\begin{te}
  In this paper we do not know neither the density $\ydf=\xdf\oast\edf$ of the contaminated observations nor the error density $\edf$, but we have at our disposal two independent samples of independent and identically distributed (\iid) random variables of size $\ssY\in\Nz$ and $\ssE\in\Nz$, respectively:
  \begin{equation}\label{i:obs}
    \rY_i\sim\ydf,\quad i\in\nset{\ssY} := \nset{1, \ssY} := [1, \ssY] \cap \Zz, \quad\text{and}\quad
    \rE_i\sim\edf,\quad i\in\nset{\ssE}.
  \end{equation}
  In this situation, for each dimension parameter $\Di \in \Nz$ an OSE of $\xdf$ is given by
  \begin{multline}\label{i:hxdfPr}
    \hxdfPr:= \bas_0 + \sum\nolimits_{|j|\in\nset{\Di}} \hfedfmpI[j]\hfydf[j] \bas_j, \quad \text{with }\hfydf[j] := \ssY^{-1}\sum\nolimits_{i \in \nset{\ssY}} \bas_j(-\rY_i),\\
    \hfedfmpI[j] :=\hfedfI[j]\Ind{\{|\hfedf[j]|^2\geq1/\ssE\}} \text{ and } \hfedf[j] := \ssE^{-1}\sum\nolimits_{i\in \nset{\ssE}} \bas_j(-\rE_i).
  \end{multline}
  The threshold using the indicator function $\Ind{\{|\hfedf[j]|^2\geq1/\ssE\}}$, accounts for the uncertainty caused by estimating $\fedf[j]$ by $\hfedf[j]$.
  It corresponds to $\hfedf[j]$'s noise level as an estimator of $\fedf[j]$ which is a natural choice (cf. \cite{Neumann1997}, p. 310f.).
\end{te}
\begin{te}
  Thanks to the properties of the sequences $(\hfydf[j])_{j \in \Zz}$, and $(\hfedfmpI[j])_{j \in \Zz}$, for any $k$ in $\Nz$, the estimator $\hxdfPr$ is a real valued function integrating to $1$.
  It is not necessarily positive valued, however, one might project the estimator on $\cD$, leading to an even smaller quadratic error.
  Nevertheless $\hxdfPr$ depends on a dimension parameter $\Di$ whose choice essentially determines the estimation accuracy.
\end{te}

\begin{te}
  In \cite{JohannesSchwarz2013a}, a minimax criterion is used to formulate optimality.
  It is shown that, by choosing the dimension parameter properly, the maximal risk of an OSE as in \eqref{i:hxdfPr} reaches the lower bound over Sobolev-like ellipsoids.
  However, the optimal choice of the dimension depends on the unknown ellipsoids.
  A fully data-driven selection based on a penalised contrast method is proposed and it is shown to yield minimax optimal rates for a large family of such ellipsoids.
  This selection procedure is inspired by the work of \citet{BarronBirgeMassart1999}, which was applied in the case of known error density by \citet{ComteTaupin}.
  For an extensive overview of model selection by penalised contrast, the reader may refer to \citet{Massart2007}.
  More precisely, \cite{JohannesSchwarz2013a} introduce an upper bound $\widehat{M}$ for the dimension parameter, and penalties $(\peneSv)_{\Di \in \nset{\widehat{M}}}$, depending on the samples $(\rY_{i})_{i \in \nset{\ssY}}$, and $(\rE_{i})_{i \in \nset{\ssE}}$, but neither on $\xdf$ nor $\edf$.
  Then, the fully data-driven estimator is defined as
  \begin{equation}\label{i:hDi}
    \hxdfPr[\tDi]:=\bas_{0} + \sum\nolimits_{\vert j \vert \in \nset{\tDi}} \hfedfmpI[j]\hfydf[j]
    \bas_j\quad\text{with }
    \tDi:=\argmin_{\Di\in\nset{\hM}}\{-\VnormLp{\hxdfPr}^2+\peneSv\}.
  \end{equation}
  The empirical upper bound $\widehat{M}$ proposed in \cite{JohannesSchwarz2013a} is technically rather involved and more importantly simulations suggest that it leads to values which are often much too restrictive.
\end{te}

\begin{te}
  Here, rather than a data-driven selection of a dimension parameter, we propose to sum the OSE's with positive data-driven weights adding up to one.
  Namely, given for each $\Di \in \nset{\ssY}$, the OSE's as in \eqref{i:hxdfPr}, and a random weight $\We[\Di] \in [0, 1]$, we consider the convex sum
  \begin{equation}\label{i:hxdfAg}
    \hxdf[{\We[]}]=\sum\nolimits_{\Di \in \nset{\ssY}} \We\hxdfPr, \quad \text{with } \sum\nolimits_{\Di\in \nset{\ssY}}\We[\Di]=1.
  \end{equation}
  Introducing the model selection weights,
  \begin{equation}\label{au:de:msWe}
    \msWe := \Ind{\{\Di = \hDi\}}, \quad \Di \in \nset{\ssY}, \quad \text{ with } \hDi:=\argmin_{\Di\in\nset{\ssY}}\{-\VnormLp{\hxdfPr}^2+\peneSv\}
  \end{equation}
  allows us to consider the model selected estimator $\hxdfPr[\hDi]=\hxdfAg[{\msWe[]}]=\sum_{\Di\in\nset{\ssY}}\msWe\hxdfPr$ as a data-driven weighted sum, avoiding a restrictive empirical upper bound $\widehat{M}$ as in \eqref{i:hDi}.
\end{te}

\begin{te}
  We study a second shape of random weights, motivated by a Bayesian approach in the context of an inverse Gaussian sequence space model and its iterative extension respectively described in \cite{johannes2020adaptive} and \cite{loizeau2020hierarchical}.
  For some constant $\rWc\in\Nz$ we define Bayesian weights
  \begin{equation}\label{au:de:erWe}
    \erWe:=\frac{\exp(-\rWn\{-\VnormLp{\hxdfPr}^2+\peneSv\})}{\sum_{l=1}^{\ssY}\exp(-\rWn\{-\VnormLp{\hxdfPr[l]}^2+\peneSv[l]\})}, \quad \Di\in\nset{\ssY}.
  \end{equation}
  Note that in \eqref{au:de:msWe} and \eqref{au:de:erWe} the quantity $\VnormLp{\hxdfPr}^2=\sum_{j=-\Di}^{\Di}|\hfedfmpI[j]|^2|\hfydf[j]|^2$ can be calculated from the data without any prior knowledge about the error density $\edf$.
  Thereby, as the sequence of penalties $(\peneSv[\Di])_{\Di \in \nset{\ssY}}$ given in bellow \eqref{au:de:peneSv} does not involve any prior knowledge neither of $\xdf$ nor $\edf$, the weights in \eqref{au:de:msWe} and \eqref{au:de:erWe} are fully data-driven.
\end{te}

\begin{te}
  Let us emphasise the role of the parameter $\rWc$ used in \eqref{au:de:erWe}.
  If $\hDi$ as in \eqref{au:de:msWe} minimises uniquely the penalised contrast function, then it is easily seen that for each $\Di \in \nset{\ssY}$ the Bayesian weight $\erWe[\Di]$ converges to the model selection weight $\msWe$ as $\rWc \to \infty$.
  We shall see that the fully data-driven weighted sum $\hxdf[{\We[]}]$ with Bayesian weights $\We[] = \erWe[]$ or model selection weights $\We[] = \msWe[]$ yields minimax optimal convergence rates over Sobolev-like ellipsoids.
  Thus, the theory presented here does not give a way to chose the parameter $\rWc$.
  However, simulations suggest that the Bayesian weights lead to more stable results as it is often recorded in the field of estimator aggregation.
\end{te}

\begin{te}
  The shape of the weighted sum $\hxdf[{\We[]}]$ is similar to the form studied in the estimator aggregation literature.
  Aggregation in the context of regression problems is considered, for instance, in \citet{dalalyan2008aggregation, tsybakov2014aggregation, rigollet2012kullback, dalalyan2012sparse, bellec2015sharp}), while \citet{rigollet2007linear} study density estimation.
  Traditionally, the aggregation of a family of arbitrary estimators is performed through an optimisation program for the random weights, and the goal is to compare the convergence rate of the aggregation to the one of the best estimator in the family.
  Here, while we restrict ourselves to OSE's, their number is as large as the sample size.
  The random weights are given explicitly without an optimisation program and do not rely on a sample splitting.
  In addition, we allow for a degenerated cases where one OSE receives all the weight of the sum.
\end{te}

\begin{te}
  This paper is organised as follows.
  In \cref{ak} assuming that the error density $\edf$ is known, we introduce a family of OSE's.
  We briefly recall the oracle and minimax theory before introducing model selection and Bayesian weights respectively similar to \eqref{au:de:msWe}, and \eqref{au:de:erWe}, which still depend on characteristics of the error density. The weighted sum of the OSE's is thus only partially data-driven.
  We derive non-asymptotic upper bounds for the quadratic risk and the maximal quadratic risk over Sobolev-like ellipsoids of the partially data-driven estimator.
  In \cref{au}, dismissing the knowledge of the density $\edf$ an additional sample of the noise is observed.
  Choosing the weights in \eqref{au:de:msWe}, and \eqref{au:de:erWe} fully data-driven we derive non-asymptotic upper risk bounds for the now fully data-driven weighted sums of OSE's.
  In \cref{ak,au} we compute rates which can be obtained in different configurations for the smoothness of the density of interest $\xdf$ and the error density $\edf$.
  The rates (strictly) match the optimal oracle or minimax rates for a large variety of cases, and feature otherwise at most a deterioration by a logarithmic factor.
  We illustrate in \cref{si} the reasonable performance of the fully data-driven weighted sum of OSE's by a simulation study.
  All technical proofs are deferred to the Appendix.
\end{te}

\section{Partially data-driven aggregation: known error density}\label{ak}

\begin{te}\paragraph{Notations.} Throughout this section the error density $\edf\in\cD$ is known. Therefore,
 given an \iid $\ssY$-sample $(\rY_i)_{i\in\nset[1,]{\ssY}}$ from
 $\ydf=\xdf\oast\edf$  we denote by
 $\FuEx[\ssY]{\xdf,\edf}$ the expectation with respect to their
joint distribution $\FuVg[\ssY]{\xdf,\edf}$. The estimation of the unknown circular density $\xdf$ is based  on a
dimension reduction which we briefly elaborate first. Given the
exponential basis $\Zset{\bas_j}$
and a dimension parameter $\Di\in\Nz_0:=\Nz\cup\{0\}$ we have
the subspace $\mHiH$ spanned by the $2\Di+1$ basis functions
$\Kset[j\in\nset{-\Di,\Di}]{\bas_j}$ at our disposal. For
abbreviation, we denote by $\Proj[\Di]^{}$ and $\ProjC[\Di]$ the
orthogonal projections on $\mLp$ and its orthogonal complement
$\mLp^{\perp}$ in $\Lp[2]$, respectively. For each  $\He\in\Lp[2]$ we consider its orthogonal
projection $\HePr:=\Proj[\Di]\He$
 and its associated approximation error $\VnormLp{\HePr-\He}=\VnormLp{\ProjC[\Di]\He}$.
Note that for any
density $p\in\cD\cap\Lp[2]$ holds
$\ProjC[0]p=p-\bas_0$ and we define $\bFdfS[p]:=\Ksuite[\Di\in\Nz_0]{\bFdf{p}}\in\Rz^{\Nz_0}$ with
\begin{equation}\label{bFdf}
1\geq\bFdf{p}:=\VnormLp{\ProjC[\Di]p}/\VnormLp{\ProjC[0]p} \quad \text{ (with the convention $0/0=0$)}
\end{equation}
where
$\lim_{\Di\to\infty}\bFdf{p}=0$  due to the dominated convergence
theorem.\
\end{te}

\begin{te}
\paragraph{Risk bound.}
Keeping in mind that the error density satisfies $|\fedf[k]|>0$ for
all $k\in\Zz$, we define $\iSvS=\Ksuite[k\in\Nz_0]{\iSv[k]}\in\Rz^{\Nz_0}$, and, for any $x_{\mbullet}\in\Rz^{\Nz_0}$ let us introduce $\SiSy[\mbullet]{x}=\Ksuite[\Di\in\Nz_0]{\SiSy{x}}\in\Rz^{\Nz_0}$ with
\begin{equation}\label{iSv}
\SiSy[0]{x}:=0, \quad  \SiSy{x}:=\Di^{-1}\sum\nolimits_{j\in\nset[1,]{\Di}} x_j; \quad \text{ and } \quad 1\leq\iSv[k]:=|\fedf[k]|^{-2}=|\fedf[-k]|^{-2}.
\end{equation}
We define the OSE's in the present case similarly to \eqref{i:hxdfPr} by
\begin{equation}\label{i:txdfPr}
  \txdfPr = \bas_{0} + \sum\nolimits_{\vert j \vert \in \nset{\Di}} \fedfI[j]\hfydf[j]\bas_j,
\end{equation}
By elementary calculations for each $\Di\in\Nz_0$ the risk of $\txdfPr$  in \eqref{i:txdfPr} satisfies
\begin{equation}\label{oo:e1}
  \noRi{\txdfPr}{\xdf}{\iSv}
  + n^{-1} \VnormLp{\ProjC[0]\xdf}^2=2 \ssY^{-1}\Di\SiSv +\tfrac{\ssY+1}{\ssY}\VnormLp{\ProjC[0]\xdf}^2\sbFxdf.
\end{equation}
\end{te}

\begin{te}The quadratic risk  in the last display depends on the
dimension parameter $\Di$ and hence by selecting an optimal value it
will be minimised, which we formulate next.
For a sequence $\Ksuite[\Di\in\Nz_0]{a_{\Di}}$ of real numbers with minimal value in a set
$A\subset{\Nz_0}$ we set
$\argmin\set{a_{\Di},\Di\in A}:=\min\{\Di \in A : a_{\Di} \leq a_{j},\;\forall j \in A\}$.  For  any non-negative sequence
 $x_{\mbullet}:=\Ksuite[\Di\in\Nz_0]{x_{\Di}}$, $y_{\mbullet}:=\Ksuite[\Di\in\Nz_0]{y_{\Di}}$ and each $\Di\in\Nz_0$ define
\begin{multline}\label{oo:de:doRao}
  \doRao[x_{\mbullet},y_{\mbullet}]{\Di}:=[x_{\Di}\vee\ssY^{-1} \,\Di\, y_{\Di}]
  :=\max\vectB{x_{\Di}, \ssY^{-1}\,\Di\, y_{\Di} },\\
  \hfill
  \noDio[x_{\mbullet},y_{\mbullet}]:=\argmin\Nset[{\Di\in\Nz_0}]{\doRao[x_{\mbullet},y_{\mbullet}]{\Di}}
  \quad\text{ and }\hfill\\
  \noRao[x_{\mbullet},y_{\mbullet}]:=\min\Nset[{\Di\in\Nz_0}]{\doRao[x_{\mbullet},y_{\mbullet}]{\Di}}=\doRao[x_{\mbullet},y_{\mbullet}]{\noDio[x_{\mbullet},y_{\mbullet}]}.
\end{multline}
\end{te}

\begin{rem}\label{oo:rem:ora}Here and subsequently, our upper bounds of the risk derived from
  \eqref{oo:e1} make use of the definitions \eqref{oo:de:doRao}, for example, replacing the sequences
  $x_{\mbullet}$ and $y_{\mbullet}$ by $\sbFxdfS$ and
  $\SiSvS$, respectively.  However, in what follows the sequence
  $x_{\mbullet}$ and $y_{\mbullet}$ is always monotonically non-increasing
  and non-decreasing, respectively, with $x_0\leq1\leq y_0$ and
   $\lim_{\Di\to\infty}x_{\Di}=0=\lim_{\Di\to\infty}y_{\Di}^{-1}$.  In this
   situations by construction hold  $\noDio[x_{\mbullet},y_{\mbullet}]\in\nset[1,]{\ssY}$ and $\noRao[x_{\mbullet},y_{\mbullet}]\geq \ssY^{-1}$ for all
  $\ssY\in\Nz$, and $\lim_{\ssY\to\infty}\noRao[x_{\mbullet},y_{\mbullet}]=0$. For the latter observe
  that for each $\delta>0$ there is $\Di_{\delta}\in\Nz$ and
  $\ssY_\delta\in\Nz$ such that $x_{\Di_\delta}\leq \delta$ and
  $\Di_{\delta} y_{\Di_\delta} \ssY^{-1}\leq\delta$, $\noRao[x_{\mbullet},y_{\mbullet}]\leq\doRao[x_{\mbullet},y_{\mbullet}]{\Di_\delta}\leq \delta$,  for all $\ssY\geq \ssY_{\delta}$.
  We shall use those
  elementary findings in the sequel without further reference. \remEnd
\end{rem}

\begin{te}Throughout the paper we  shall distinguish for the deconvolution
density $\xdf$ and hence it's associated sequence
$\bFS:=\bFxdfS\in\pRz^{\Nz_0}$ of approximation errors as in \eqref{bFdf}
the two cases: \begin{inparaenum}[i]
\item[\mylabel{oo:xdf:p}{\dgrau\bfseries{(p)}}] there is $K\in\Nz_0$
  with  $\bF[K]=0$ and $\bF[K-1]>0$ (with the convention
  $\bF[-1]:=1$), and
\item[\mylabel{oo:xdf:np}{\dgrau\bfseries{(np)}}] for all $K\in\Nz$
  holds $\bF[K]>0$. \end{inparaenum} Let us stress, that for any
monotonically non-decreasing sequence $y_{\mbullet}$ with
  $y_0\geq1$, the order of the rate
  $\Nsuite[n]{\noRao[\sbFS,y_{\mbullet}]}$ defined in \eqref{oo:de:doRao}
  with $x_{\mbullet}$ replaced by $\bFS$  in case \ref{oo:xdf:p} and \ref{oo:xdf:np}
is parametric and non-parametric, respectively. More precisely, in
case \ref{oo:xdf:p} it holds $\noDio[\bFS,y_{\mbullet}]=K$ and
$\noRao[\bFS,y_{\mbullet}]=\ssY^{-1} K y_{K}$ for all
$\ssY>{Ky_{K}}/{\sbF[K-1]}$, while  in case
\ref{oo:xdf:np} holds
$\lim_{\ssY\to\infty}\noDio[\bFS,y_{\mbullet}]=\infty$ and $\lim_{\ssY\to\infty}\ssY\noRao[\sbFS,y_{\mbullet}]=\infty$.
\end{te}

\begin{te}\paragraph{Oracle optimality.} Coming back to the identity \eqref{oo:e1}  and  exploiting the definition
\eqref{oo:de:doRao} with $x_{\mbullet}$ and $y_{\mbullet}$, respectively,
 replaced by $\sbFS:=\sbFxdfS$ and  $\oiSvS$ as in \eqref{iSv} it follows immediately
  \begin{equation}\label{oo:oub}
\inf\NsetB[\Di\in\Nz]{\noRi{\txdfPr}{\xdf}{\iSv}}\hfill
    \leq \noRi{\txdfPr[\noDio]}{\xdf}{\iSv}
    \leq 2[1\vee\VnormLp{\ProjC[0]\xdf}^2]\noRao.
  \end{equation}
On the other hand with
  $[a\wedge b]:=\min(a,b)$ for $a,b\in\Rz$ from
  \eqref{oo:e1}  we also conclude
  \begin{multline}\label{oo:olb}
    \inf\NsetB[\Di\in\Nz]{\noRi{\txdfPr}{\xdf}{\iSv}}
    \geq
    \big([\VnormLp{\ProjC[0]\xdf}^2\wedge2]-\tfrac{\VnormLp{\ProjC[0]\xdf}^2}{\ssY\noRao}\big)\;\noRao.\hfill
  \end{multline}
 For each $\ssY\in\Nz$ combining  \eqref{oo:oub} and \eqref{oo:olb}   $\noRao$,
 $\noDio$  and   $\txdfPr[\noDio]$, respectively, is an oracle rate,
 oracle dimension and oracle optimal estimator (up to a constant), if
the leading factor on the right hand side in \eqref{oo:olb} is uniformly in
$\ssY$ bounded away from zero.
Note that $\noRao$ is in case \ref{oo:xdf:np}  always an orale rate,
while  in
case \ref{oo:xdf:p} whenever $K\oiSv[K]>
\,[1\vee\tfrac{1}{2}\VnormLp{\ProjC[0]\xdf}^2]$.
\end{te}

\begin{te}\paragraph{Aggregation.}
We call aggregation weights any $\We[]:=(\We)_{\Di\in\nset[1,]{\ssY}} \in[0,1]^\ssY$ defining on the set $\nset[1,]{\ssY}$ a
discrete probability measure $\FuVg{\We[]}(\{\Di\}):=\We$, $\Di\in\nset[1,]{\ssY}$.
We consider here and subsequently a weighted sum $\txdf[{\We[]}]:=\sum\nolimits_{\Di\in\nset[1,]{\ssY}}\We\txdfPr$ of the orthogonal series estimators defined in \eqref{i:txdfPr}.
Clearly, the coefficients $\Zsuite{\ftxdfPr[{\We[]}]{j}}$ of $\txdf[{\We[]}]$ satisfy $\ftxdfPr[{\We[]}]{j}=0$ for
  $|j|>\ssY$, and  for any  $|j|\in\nset[1,]{\ssY}$ holds
$\ftxdfPr[{\We[]}]{j}=\sum_{\Di\in\nset[1,]{\ssY}}\We\ftxdfPr{j}
=\FuVg{\rWe[]}(\nset{|j|,\ssY})\fedfI[j]\hfydf[j]$.
We  note that by construction
$\ftxdfPr[{\We[]}]{0}=1$, $\ftxdfPr[{\We[]}]{-j}=\overline{\ftxdfPr[{\We[]}]{j}}$ and $1\geq
|\ftxdfPr[{\We[]}]{j}|$. Hence, $\txdfPr[{\We[]}]$ is real and integrates to one, however,
it is not necessary non-negative.   Our aim is to prove an upper bound for its risk
  $\noRi{\txdfPr[{\We[]}]}{\xdf}{\iSv}$ as well its maximal risk over Sobolev-like ellipsoids. For
  arbitrary aggregation weights and penalty sequence the next lemma establishes an
  upper bound for the loss of the aggregated estimator. Selecting suitably
  the weights and penalties  this bound provides in the sequel our key
  argument.
\end{te}

\begin{lem}\label{co:agg}
  Consider a weighted sum $\txdfPr[{\We[]}]$ with arbitrary aggregation weights $\We[]$ and non-negative penalty terms $(\pen)_{\Di\in\nset[1,]{\ssY}}$. For any $\mDi,\pDi\in\nset[1,]{\ssY}$  holds
  \begin{multline}\label{co:agg:e1}
    \VnormLp{\txdfPr[{\We[]}]-\xdf}^2\leq \tfrac{2}{7}\pen[\pDi] +2\VnormLp{\ProjC[0]\xdf}^2\bias[\mDi]^2(\xdf)\\\hfill
    +2\VnormLp{\ProjC[0]\xdf}^2\FuVg{\We[]}(\nsetro[1,]{\mDi})+\tfrac{2}{7}\sum_{\Di\in\nsetlo{\pDi,\ssY}}\pen\We\Ind{\{\VnormLp{\txdfPr-\xdfPr}^2<\pen/7\}}\\
+2\sum_{\Di\in\nset{\pDi,\ssY}}\vectp{\VnormLp{\txdfPr-\xdfPr}^2-\pen/7}
+\tfrac{2}{7}\sum_{\Di\in\nsetlo{\pDi,\ssY}}\pen\Ind{\{\VnormLp{\txdfPr-\xdfPr}^2\geq\pen/7\}}.
\end{multline}
\end{lem}

\begin{rem}\label{co:agg:rem}Keeping \eqref{co:agg:e1} in mind let us
  briefly outline
  the principal arguments of our aggregation strategy.
  Selecting the values  $\pDi$ and $\mDi$
 close to the oracle dimension $\noDio$ the first two terms
  in the upper bound of \eqref{co:agg:e1} are of the order of the
  oracle rate. On the other hand  the
  weights  are in the sequel selected
   such that
  the third and fourth are
  negligible with respect to the oracle rate, while the choice of the
  penalties  allows as usual to bound the deviation of the last two terms
by concentration inequalities.\remEnd
\end{rem}

\begin{te}
  For some constant $\rWc\in\Nz$, we consider  either \pcw
  \begin{equation}\label{ak:de:rWe}
      \rWe:=\frac{\exp(-\rWn\{-\VnormLp{\txdfPr}^2+\penSv\})}{\sum\nolimits_{l \in \nset{n}}\exp(-\rWn\{-\VnormLp{\txdfPr[l]}^2+\penSv[l]\})},\quad \Di\in\nset{1,\ssY}.
    \end{equation}
    or \msw
    \begin{equation}\label{ak:de:msWe}
    \msWe := \Ind{\{\Di = \tDi\}}, \quad \Di \in \nset{\ssY}, \quad \text{ with } \tDi:=\argmin_{\Di\in\nset{\ssY}}\{-\VnormLp{\txdfPr}^2+\penSv\}
  \end{equation}
  respectively similar to the ones defined in \eqref{au:de:erWe} and \eqref{au:de:msWe}.
  Until now we have not  specified the sequence of penalty
  terms. For  a sequence
  $x_{\mbullet}\in\Rz^{\Nz_0}$ and  $\Di\in\Nz_0$ define
  \begin{equation}\label{ak:de:LiSy}
  x_{(\Di)}:=\max\{x_j:j\in\nset{0,\Di}\},\quad
  \liSy{x}:=\tfrac{|\log (\Di x_{(\Di)}\vee(\Di+2))|^2}{|\log(\Di+2)|^2},\;\text{ and }
  \LiSy{x}:=\liSy{x}x_{(\Di)}.
  \end{equation}
  Given $\iSvS \in \pRz^{\Nz_0}$ as in \eqref{iSv} and a numerical  constant $\cpen>0$ we use
  \begin{equation}\label{ak:de:penSv}
    \penSv:= \penD,\quad\Di\in\Nz.
  \end{equation}
  as penalty terms. For the theoretical results below we need that the numerical
  constant satisfies $\cpen\geq84$. However, for a practical  application
  this values is generally too large and a suitable constant might be chosen by
  preliminary calibration experiments see \cite{BaudryMaugisMichel2012}.\\
\end{te}

\begin{te}
  We derive bounds for the risk of the weighted sum estimator $\txdfAg$ with Bayesian weights
  and the model selected estimator $\txdf[\tDi] = \txdf[{\msWe[]}]$ by applying
  \cref{co:agg}.
  From definition \eqref{oo:de:doRao} replacing $x_{\mbullet}$ and
  $y_{\mbullet}$, respectively, by $\sbFS=\sbFxdfS$ and $\LiSvS$ we consider $\doRaL{\Di}$
  for each $\ssY,\Di\in\Nz$. Note that by construction \eqref{iSv} and \eqref{ak:de:LiSy}, we have $\noRaL\geq\noRao$
for all $\ssY\in\Nz$.  We denote in the sequel by $\cst{}$ an universal finite numerical constant with value
changing possibly from line to line.
\end{te}

\begin{pr}\label{ak:ag:ub:pnp} Consider an aggregation $\txdfAg[{\We[]}] = \sum_{\Di \in \nset{n}} \We \txdfPr$  using either \pcw $\We[]:=\rWe[]$
  as in \eqref{ak:de:rWe} or \msw $\We[]:=\msWe[]$
  as in \eqref{ak:de:msWe} and  penalties $(\penSv)_{\Di\in\nset[1,]{\ssY}}$ as in \eqref{ak:de:penSv} with numerical constant $\cpen\geq84$.  Let $\Di_{\ydf}:=\floor{3(400)^2\Vnormlp[1]{\fydf}^2}$.
\begin{Liste}[]
\item[\mylabel{ak:ag:ub:pnp:p}{\dgrau\bfseries{(p)}}]Assume there is $K\in\Nz_0$
  with $1\geq \bFxdf[(K-1)]>0$ and $\bFxdf[K]=0$. If $K=0$ we set $c_{\xdf}:=0$ and
  $c_{\xdf}:=\tfrac{4\cpen}{\VnormLp{\ProjC[0]\xdf}^2\sbFxdf[(K-1)]}$,
   otherwise.  For $\ssY\in\Nz$ let
  $\sDi{\ssY}:=\max\{\Di\in\nset[1,]{\ssY}:\ssY>c_{\xdf}\DipenSv\}$, if
  the defining set is not empty, and
  $\sDi{\ssY}:=\ceil{\Di_{\ydf}\log(2+\ssY)}$ otherwise.
There is a finite constant $\cst{\xdf,\edf}$
given in \eqref{ak:ag:ub:p:e1} depending only on $\xdf$ and $\edf$ such that for all $n\in\Nz$ holds%
\begin{equation}\label{ak:ag:ub:pnp:e1}
  \noRi{\txdfAg[{\We[]}]}{\xdf}{\iSv}
  \leq
  \cst{}\VnormLp{\ProjC[0]\xdf}^2\big[
  \ssY^{-1}\vee\exp\big(\tfrac{-\cmiSv[\sDi{\ssY}]\sDi{\ssY}}{\Di_{\ydf}}\big)\big]
  + \cst{\xdf,\edf}\ssY^{-1}.
\end{equation}
\item[\mylabel{ak:ag:ub:pnp:np}{\dgrau\bfseries{(np)}}] Assume that
  $\bias(\xdf)>0$ for all  $\Di\in\Nz$.
There is a finite constant $\cst{\xdf,\edf}$ given in
\eqref{ak:ag:ub:pnp:p8} depending only on $\xdf$ and $\edf$ such that for all
$\ssY\in\Nz$  holds
 \begin{multline}\label{ak:ag:ub:pnp:e2}
   \noRi{\txdfAg[{\We[]}]}{\xdf}{\iSv}
    \leq
   \cst{}\,[\VnormLp{\ProjC[0]\xdf}^2\vee1]\,\noRaA  +\cst{\xdf,\edf}\,\ssY^{-1}\\
  \text{with }\quad \noRaA:=\min_{\Di\in\nset[1,]{\ssY}}\setB{\big[\doRaL{\Di}\vee\exp\big(\tfrac{-\cmiSv\Di}{\Di_{\ydf}}\big)\big]}.
\end{multline}
\end{Liste}
\end{pr}

\begin{co}\label{ak:ag:ub2:pnp}
  Let the assumptions of \cref{ak:ag:ub:pnp} be satisfied.
  \begin{Liste}[]
  \item[\mylabel{ak:ag:ub2:pnp:p}{\dgrau\bfseries{(p)}}]
    If  in addition
    \begin{inparaenum}\item[\mylabel{ak:ag:ub2:pnp:pc}{\normalfont\rmfamily{\dgrau\bfseries(A1)}}]
      there is $\ssY_{\xdf,\edf}\in\Nz$ such that for all
      $\ssY\geq \ssY_{\xdf,\edf}$ holds
      $\cmiSv[\sDi{\ssY}]\sDi{\ssY}\geq
      \Di_{\ydf}\log\ssY$, \end{inparaenum} then there is a constant $\cst{\xdf,\edf}$ depending
    only on $\xdf$ and $\edf$ such that for all $n\in\Nz$ holds
    $\noRi{\txdfAg[{\We[]}]}{\xdf}{\iSv} \leq
    \cst{\xdf,\edf}\;\ssY^{-1}$.
  \item[\mylabel{ak:ag:ub2:pnp:np}{\dgrau\bfseries{(np)}}]
    If in addition
    \begin{inparaenum}\item[\mylabel{ak:ag:ub2:pnp:npc}{{\dgrau\bfseries(A2)}}:]
      there is  $\ssY_{\xdf,\edf}\in\Nz$ such that $\oDi{\ssY}:=\noDiL$
      as in \eqref{oo:de:doRao} for all $\ssY\geq \ssY_{\xdf,\edf}$ satisfies
      $\oDi{\ssY}\cmSv[\oDi{\ssY}]\geq\Di_{\ydf}|\log\noRaL|$,
    \end{inparaenum}  then there is a constant $\cst{\xdf,\edf}$ depending
    only on $\xdf$ and $\edf$ such that $\noRi{\txdfAg[{\We[]}]}{\xdf}{\iSv}
    \leq \cst{\xdf,\edf}\noRaL$ for all $n\in\Nz$.
  \end{Liste}
\end{co}

\begin{il}\label{ak:il}Here and subsequently, we use for two strictly positive sequences
 $\Nsuite[n]{a_{n}}$ and $\Nsuite[n]{b_{n}}$ the notation
 $a_{n}\sim b_{n}$  if the sequence $\Nsuite[n]{a_{n}/b_{n}}$
 is bounded away both from zero and infinity.
  We  illustrate the last results considering usual
  behaviours for the sequences $\bFS$ and
  $\iSvS$.  Regarding the error density $\edf$ we consider for $a>0$ the
  following two cases \begin{inparaenum}[i]
\item[\mylabel{il:edf:o}{\dg\bfseries{(o)}}]   $\iSv[\Di]\sim
  \Di^{2a}$ and \item[\mylabel{il:edf:s}{\dg\bfseries{(s)}}]  $\iSv[\Di]\sim
  \exp(\Di^{2a})$.
\end{inparaenum}
The error density $\edf$ is called ordinary smooth in case \ref{il:edf:o} and
super smooth in case \ref{il:edf:s},
and it holds, respectively,
  \begin{inparaenum}[i]
\item[\mylabel{ak:il:edf:o}{\dg\bfseries{(o)}}]
  $\sDi{\ssY}\sim\ssY^{1/(2a+1)}$ and
    $\sDi{\ssY}\cmSv[\sDi{\ssY}]\sim\ssY^{1/(2a+1)}$, and
\item[\mylabel{ak:il:edf:s}{\dg\bfseries{(s)}}]
    $\sDi{\ssY}\sim(\log \ssY)^{1/(2a)}$ and
    $\sDi{\ssY}\cmSv[\sDi{\ssY}]\sim (\log \ssY)^{2+1/(2a)}$.
  \end{inparaenum}
  Clearly in both cases \ref{ak:ag:ub2:pnp:pc}
  holds true and hence employing \cref{ak:ag:ub2:pnp}
  \ref{ak:ag:ub2:pnp:p} the aggregated estimator  attains the
  parametric rate. On the other hand, for
  \ref{ak:ag:ub2:pnp:np}
  we use for the deconvolution density $\xdf$  as particular
specifications \begin{inparaenum}[i]
\item[\mylabel{il:xdf:o}{\dg\bfseries{(o)}}] $|\fxdf[\Di]|^2\sim
  \Di^{-2p-1}$ and
\item[\mylabel{il:xdf:s}{\dg\bfseries{(s)}}] $|\fxdf[\Di]|^2\sim
 \Di^{2p-1} \exp(-\Di^{2p})$
\end{inparaenum} with $p>0$.
\\[2ex]
\centerline{\begin{inparaenum}[i]
\begin{tabular}{@{}llllll@{}}
  \toprule
 Order  & $\sbF$
  &$\iSv[\Di]$ &$\noRao$ &  $\noRaL$& $\noRaA$\\
  \midrule
 \item[\mylabel{ak:oo:oo}{\dg\bfseries{[o-o]}}]
        &$\Di^{-2p}$&$\Di^{2a}$
               & $\ssY^{\tfrac{-2p}{2p+2a+1}}$ & $\ssY^{\tfrac{-2p}{2p+2a+1}}$ & $\ssY^{\tfrac{-2p}{2p+2a+1}}$ \\
 \item[\mylabel{ak:oo:os}{\dg\bfseries{[o-s]}}] &$\Di^{-2p}$& $e^{\Di^{2a}}$ & $(\log\ssY)^{\tfrac{-p}{a}}$& $(\log\ssY)^{\tfrac{-p}{a}}$  &$(\log \ssY)^{\tfrac{-p}{a}}$\\
 \item[\mylabel{ak:oo:so}{\dg\bfseries{[s-o]}}]  &$e^{-\Di^{2p}}$&$\Di^{2a}$
   & $(\log \ssY)^{\tfrac{2a+1}{2p}}\ssY^{-1}$&  $(\log\ssY)^{\tfrac{2a+1}{2p}}\ssY^{-1}$
   & $\left\{\begin{array}{@{}ll@{}}(\log\ssY)^{\tfrac{2a+1}{2p}}\ssY^{-1}&:\text{\smaller[1]$p<1/2$},\\
                              (\log\ssY)^{2a+1}\ssY^{-1}&:\text{\smaller[1]$p\geq1/2$}.\end{array}\right.$\\
  \bottomrule
\end{tabular}
\end{inparaenum}}\\[2ex]
To calculate the order in the last table we used that the dimension parameter $\oDi{\ssY}:=\noDiL$ satisfies   \begin{inparaenum}[i]
  \item[\mylabel{ak:il:oo}{\dg\bfseries{[o-o]}}]
  $\oDi{\ssY}\sim\ssY^{1/(2p+2a+1)}$ and
  $\cmiSv[\oDi{\ssY}]\oDi{\ssY}\sim\ssY^{1/(2p+2a+1)}$,
  \item[\mylabel{ak:il:os}{\dg\bfseries{[o-s]}}]
   $\oDi{\ssY}\sim(\log\ssY)^{1/(2a)}$ and
   $\cmiSv[\oDi{\ssY}]\oDi{\ssY}\sim(\log\ssY)^{2+1/(2a)}$, and
  \item[\mylabel{ak:il:so}{\dg\bfseries{[s-o]}}]
   $\oDi{\ssY}\sim(\log\ssY)^{(a+1/2)/p}$ and
   $\cmiSv[\oDi{\ssY}]\oDi{\ssY}\sim(\log\ssY)^{1/(2p)}$.
 \end{inparaenum}
We note that in each of the three cases
the order of $\noRaL$ and the order of the oracle rate $\noRao$
coincide. Moreover, the additional assumption
\ref{ak:ag:ub2:pnp:npc} in \cref{ak:ag:ub2:pnp} \ref{ak:ag:ub2:pnp:np}
is satisfied in case \ref{ak:oo:oo} and \ref{ak:oo:os}, but in \ref{ak:oo:so}
only with $p<1/2$. Consequently, in this situations due to
\cref{ak:ag:ub2:pnp} \ref{ak:ag:ub2:pnp:np} the partially data-driven
aggregation is oracle optimal (up to a constant). Otherwise,
the upper bound $\noRaA$ in \cref{ak:ag:ub:pnp} \eqref{ak:ag:ub:pnp:e2}  faces a
detoriation compared to the rate $\noRaL$.
  In case \ref{ak:oo:so}
  with $p\geq1/2$
  setting $\sDi{\ssY}:=\Di_{\ydf}|\log\noRaL|\sim(\log\ssY)$ the upper bound
  $\noRaA\leq\doRaL{\sDi{\ssY}}\sim(\log\ssY)^{2a+1}\ssY^{-1}$
  features a deterioration at most by a logarithmic factor
  $(\log\ssY)^{(2a+1)(1-1/(2p))}$ compared to the oracle rate $(\log \ssY)^{(2a+1)/(2p)}\ssY^{-1}$.\ilEnd
\end{il}

\begin{te}\paragraph{Minimax optimality.} Rather  than considering for each
  $\Di\in\Nz$ the risk of the OSE $\txdfPr$ for given $\xdf$ and
  $\edf$ we shall measure now its accuracy by a maximal risk over
  pre-specified classes of densities determining
  a priori conditions on $\xdf$ and $\edf$, respectively. For an
  arbitrary positive sequence $x_{\mbullet}\in\pRz^{\Nz_0}$ and
  $\He\in\Lp[2]$ we write shortly
  $\Vnorm[x]{\He}^2:=\sum_{j\in\Zz}x_{|j|}|\fHe{j}|^2$.  Given
  strictly positive sequences $\xdfCwS=\Ksuite[\Di\in\Nz_0]{\xdfCw}$
  and $\edfCwS=\Ksuite[\Di\in\Nz_0]{\edfCw}$, and constants
  $\xdfCr,\edfCr\geq1$ we define
  \begin{equation*}
    \rwCxdf:=\{p\in\cD:
    \Vnorm[1/{\xdfCw[]}]{p}^2\leq\xdfCr\}\quad\text{ and }\quad
    \rwCedf:=\{p\in\cD:\edfCr^{-1}\leq\edfCw[j]|\fou[|j|]{p}|^{2}\leq \edfCr,\;\forall j\in\Zz\}.
  \end{equation*}
  Here and subsequently, we suppose  the following minimal regularity
  conditions are satisfied. \setcounter{ass}{2}
  \begin{ass}\label{ass:Cw}The sequences
    $\xdfCwS$, $\IedfCwS$   are monotonically non-increasing  with
  $\xdfCw[0]=1=\edfCw[0]$,
  $\lim_{\Di\to\infty}\xdfCw=0=\lim_{\Di\to\infty}\edfCw^{-1}$ and
  $\sum_{\Di\in\Nz_0}\xdfCw/\edfCw=\Vnormlp[1]{\xdfCwS/\edfCwS}<\infty$.
  \end{ass}%
  We shall emphasize that for $\Di\in\Nz_0$,
  $\xdf\in\rwCxdf$ and $\edf\in\rwCedf$ hold $\VnormLp{\ProjC[0]\xdf}^2\sbFxdf\leq
\xdfCr\xdfCw[\Di]$ and $1/\edfCr\leq\oiSv[\Di]/\edfCwo\leq\edfCr$ with
$\edfCwo=\Di^{-1}\sum_{j\in\nset[1,]{\Di}}\edfCw[j]$  which we use in the sequel
without further reference.
Exploiting again the identity \eqref{oo:e1}  and  the definition
\eqref{oo:de:doRao} with $x_{\mbullet}$ and $y_{\mbullet}$, respectively,
replaced by $\xdfCwS$ and
$\edfCwoS$ it follows for all $\Di,\ssY\in\Nz$
\begin{equation}\label{oo:e4}
    \nmRi{\txdfPr}{\rwCxdf}{\rwCedf}
    \leq
(2\edfCr+\xdfCr)\dmRao[\xdfCwS,\edfCwoS]{\Di}.
\end{equation}
The upper bound in the last display depends on the dimension parameter
$\Di$ and hence by choosing an optimal value $\nmDio[{\xdfCwS},\edfCwoS]$ the upper bound
will be minimised.
From \eqref{oo:e4} we deduce that
$\nmRi{\txdfPr[{\nmDio[\xdfCwS,\edfCwS]}]}{\rwCxdf}{\rwCedf}\leq(2\edfCr+\xdfCr)\nmRao$ for
all $n\in\Nz$. On the other
  hand  \cite{JohannesSchwarz2013a} have shown  that for all $\ssY\in\Nz$
  \begin{equation}\label{mm:lbn}
  \inf\nolimits_{\txdf}\nmRi{\txdf}{\rwCxdf}{\rwCedf}\geq \cst{}\nmRao,
\end{equation}
 where  $\cst{}>0$ and the infimum is taken over all
  possible estimators $\txdf$ of $\xdf$.  Consequently,
  $\Nsuite[n]{\nmRao}$,  $\Nsuite[n]{\nmDio}$
  and  $\Nsuite[n]{\txdfPr[\nmDio]}$, respectively, is a
  minimax rate, a minimax dimension and minimax optimal estimator (up to a
  constant).
\end{te}

\begin{te}\paragraph{Aggregation.}
  Exploiting \cref{co:agg} we derive now bounds for the
  maximal risk
  of
  the aggregated estimator $\txdfAg[{\We[]}]$ using either \pcw
   $\We[]:=\rWe[]$ as in \eqref{ak:de:rWe} or \msw
   $\We[]:=\msWe[]$ as in \eqref{ak:de:msWe}.
  Keeping the definition \eqref{ak:de:LiSy} in mind
  we use in the sequel that for any $\edf\in\rwCedf$  and $\Di\in\Nz_0$  hold
  \begin{equation}\label{a:ak:mrb:LiSy}
(1+\log\edfCr)^{-2}\leq\liSv/\liCw\leq (1+\log\edfCr)^2\text{ and }
 \zeta_{\edfCr}:=\edfCr(1+\log\edfCr)^2\geq\liSv\miSv/(\liCw\edfCw)\geq\zeta_{\edfCr}^{-1}.
  \end{equation}
It follows for all
$\Di,\ssY\in\Nz$, $\xdf\in\rwCxdf$ and $\edf\in\rwCedf$ immediately
  \begin{equation}\label{ak:ass:pen:mm:c}
    \xdfCr\dmRaL{\Di}\geq\VnormLp{\ProjC[0]\xdf}^2\sbFxdf\quad\text{and}\quad\cpen\zeta_{\edfCr}\dmRaL{\Di}\geq\penSv.
  \end{equation}
  Note that by construction
  $\nmRaL\geq\nmRao[\xdfCwS,\edfCwoS]$ for all $\ssY\in\Nz$.
\end{te}

\begin{pr}\label{ak:ag:ub:pnp:mm} Consider an aggregation $\txdfAg[{\We[]}]$  using either
  Bayesian weights $\We[]:=\rWe[]$
  as in \eqref{ak:de:rWe} or model selection weights $\We[]:=\msWe[]$
  as in \eqref{ak:de:msWe} and  penalties $(\penSv)_{\Di\in\nset[1,]{\ssY}}$ as in \eqref{ak:de:penSv} with numerical constant $\cpen\geq84$.
 Let \ref{ass:Cw} be satisfied and set $\Di_{\xdfCw[]\edfCw[]}:=\floor{3(400)^2\xdfCr\zeta_{\edfCr}\,\Vnormlp[1]{\xdfCwS/\edfCwS}}$.
There is a constant $\cst[\xdfCr\edfCr]{\xdfCw[]\edfCw[]}$ given in
\eqref{ak:ag:ub:pnp:mm:p1} depending only on  $\rwCxdf$  and $\rwCedf$ such that for all
$\ssY\in\Nz$
 \begin{multline}\label{ak:ag:ub:pnp:mm:e1}
 \nmRi{\txdfAg[{\We[]}]}{\rwCxdf}{\rwCedf}
    \leq \cst{}\,(\xdfCr+\zeta_{\edfCr})\,\nmRaA  +\cst[\xdfCr\edfCr]{\xdfCw[]\edfCw[]}\;\ssY^{-1}\\
    \text{with }\quad \nmRaA:=\min_{\Di\in\nset[1,]{\ssY}}\setB{\big[\dmRaL{\Di}\vee\exp\big(\tfrac{-\liCw[\Di]\Di}{\Di_{\xdfCw[]\edfCw[]}}\big)\big]}.
\end{multline}
\end{pr}

\begin{co}\label{ak:ag:ub2:pnp:mm}
  Let the assumptions of \cref{ak:ag:ub:pnp:mm} be satisfied.  If in
  addition
  \begin{inparaenum}\item[\mylabel{ak:ag:ub2:pnp:mm:c}{{\normalfont\rmfamily\dgrau\bfseries(A2')}}]
    there is $\ssY_{\xdfCw[]\edfCw[]}\in\Nz$  such that
    $\oDi{\ssY}:=\nmDiL$ as in  \eqref{oo:de:doRao} satisfies
    $\oDi{\ssY}\liCw[\oDi{\ssY}]\geq\Di_{\xdfCw[]\edfCw[]}|\log\nmRaL|$
    for all $\ssY\geq \ssY_{\xdfCw[]\edfCw[]}$,
  \end{inparaenum}
  then there is a constant $\cst[\xdfCr\edfCr]{\xdfCw[]\edfCw[]}$
  depending only on the classes $\rwCxdf$  and $\rwCedf$ such that
  $ \nmRi{\txdfAg[{\We[]}]}{\rwCxdf}{\rwCedf}\leq
  \cst[\xdfCr\edfCr]{\xdfCw[]\edfCw[]}\;\nmRaL$ for all $n\in\Nz$.
\end{co}

\begin{il}\label{ak:il:mm} We illustrate the last results  considering usual
configurations for  $\xdfCwS$ and $\edfCwS$. \\[2ex]
\centerline{\begin{inparaenum}[i]
\begin{tabular}{@{}llllll@{}}
  \toprule
  & $\xdfCw$
  &$\edfCw$ &$\nmRao[\xdfCwS,\edfCwoS]$
   & $\nmRaL$& $\nmRaA$\\
  \midrule
 \item[\mylabel{ak:mm:oo}{\dg\bfseries{[o-o]}}]
        &$\Di^{-2p}$&$\Di^{2a}$
               & $\ssY^{\tfrac{-2p}{2p+2a+1}}$ & $\ssY^{\tfrac{-2p}{2p+2a+1}}$ & $\ssY^{\tfrac{-2p}{2p+2a+1}}$ \\
 \item[\mylabel{ak:mm:os}{\dg\bfseries{[o-s]}}] &$\Di^{-2p}$& $e^{\Di^{2a}}$ & $(\log\ssY)^{\tfrac{-p}{a}}$& $(\log\ssY)^{\tfrac{-p}{a}}$  &$(\log \ssY)^{\tfrac{-p}{a}}$\\
 \item[\mylabel{ak:mm:so}{\dg\bfseries{[s-o]}}]  &$e^{-\Di^{2p}}$&$\Di^{2a}$
   & $(\log \ssY)^{\tfrac{2a+1}{2p}}\ssY^{-1}$&  $(\log\ssY)^{\tfrac{2a+1}{2p}}\ssY^{-1}$
   & $\left\{\begin{array}{@{}ll@{}}(\log\ssY)^{\tfrac{2a+1}{2p}}\ssY^{-1}&:\text{\smaller[1]$p<1/2$},\\
                              (\log\ssY)^{2a+1}\ssY^{-1}&:\text{\smaller[1]$p\geq1/2$}.\end{array}\right.$\\
  \bottomrule
\end{tabular}
\end{inparaenum}}\\[2ex]
  We note that in each of the three cases the order of $\nmRaL$
  coincide with the order of the minimax rate
  $\nmRao[\xdfCwS,\edfCwoS]$.  Moreover, the additional assumption
\ref{ak:ag:ub2:pnp:mm:c} in \cref{ak:ag:ub2:pnp:mm}
is satisfied in case \ref{ak:mm:oo} and \ref{ak:mm:os}, but in \ref{ak:mm:so}
only with $p<1/2$. Consequently, in this situations due to
\cref{ak:ag:ub2:pnp:mm}  the
partially data-driven aggregation is minimax optimal (up to a constant). Otherwise,
the upper bound $\nmRaA$  in \cref{ak:ag:ub:pnp:mm} \eqref{ak:ag:ub:pnp:mm:e1} faces a
detoriation compared to $\nmRaL$., e.g. in case \ref{ak:mm:so} with
$p\geq1/2$   by a logarithmic factor
  $(\log\ssY)^{(2a+1)(1-1/(2p))}$.\ilEnd
\end{il}

\section{Data-driven aggregation: unknown error density}\label{au}

\begin{te}In this section we dispense with any knowledge about the
  error density $\edf$. Instead we assume two independent sample
  $(\rY_i)_{i\in\nset[1,]{\ssY}}$ and $(\rE_i)_{i\in\nset[1,]{\ssE}}$ as in \eqref{i:obs}.
We denote by $\FuEx[\ssY,\ssE]{\xdf,\edf}$, $\FuEx[\ssY]{\xdf,\edf}$ and $\FuEx[\ssE]{\edf}$ the expectation with respect to their
joint distribution $\FuVg[\ssY,\ssE]{\xdf,\edf}$, and marginals
$\FuVg[\ssY]{\xdf,\edf}$, and $\FuVg[\ssE]{\edf}$, respectively.
\paragraph{Risk bound.}
Exploiting the independence assumption, the risk of the orthogonal
series estimators $\hxdfPr$
 given in \eqref{i:hxdfPr}
 can be decomposed for each $\ssY,\ssE,\Di\in\Nz$ as follows
  \begin{multline}\label{oo:e5}
    \mnoRi{\hxdfPr}{\xdf}{\iSv}
    =  \ssY^{-1}\sum_{|j|\in\nset[1,]{\Di}}\iSv[j](1-|\fydf[j]|^2)\FuEx[\ssE]{\rE}\big(\Vabs{\hfedfmpI[j]\fedf[j]}^2\big)+\VnormLp{\ProjC[0]\xdf}^2\bias^2(\xdf)\\
+\sum_{|j|\in\nset[1,]{\Di}}|\fxdf[j]|^2\FuEx[\ssE]{\rE}\big(\Vabs{\hfedf[j]-\fedf[j]}^2\Vabs{\hfedfmpI[j]}^2\big)
+\sum_{|j|\in\nset[1,]{\Di}}\fxdf[j]^2\FuVg[\ssE]{\rE}(|\hfedf[j]|^2<1/\ssE).
\end{multline}
Exploiting \cref{oSv:re} in the appendix we
control the deviations of the additional terms estimating the
error density. Therewith,
setting $\Vnorm[1\wedge\iSv/\ssE]{\ProjC[0]\xdf}^2=2\sum_{j\in\Nz}|\fxdf[j]|^2[1\wedge
\iSv[j]/\ssE]$, selecting $\oDi{\ssY}:=\noDio$ as in \eqref{oo:de:doRao}
with $\noRao=\doRao{\oDi{\ssY}}$
it follows for all $\ssY,\ssE\in\Nz$
\begin{equation}\label{oo:e6}
  \mnoRi{\hxdfPr[\oDi{\ssY}]}{\xdf}{\iSv}
  \leq
(\VnormLp{\ProjC[0]\xdf}^2+8) \noRao
+8(\cst{}+1)\moRa.
\end{equation}
\end{te}

\begin{rem}\label{oo:rem:nm}
  Note that $\VnormLp{\ProjC[0]\xdf}^2=0$ implies
  $\moRa=0$, while for $\VnormLp{\ProjC[0]\xdf}^2>0$ holds
  $\moRa\geq
  \VnormLp{\ProjC[0]\xdf}^2
  \ssE^{-1}$, and hence
  any additional term of order $\ssY^{-1}+\ssE^{-1}$
  is negligible with respect to
  $\noRa+\moRa$, since
  $\noRa\geq \ssY^{-1}$. On the other hand if
  $\Vnorm[\iSv]{\xdf}^2<\infty$ then  $\moRa\leq \ssE^{-1}\Vnorm[\iSv]{\xdf}^2$.
  Consequently, in case \ref{oo:xdf:p}  the order of the upper bound is parametric in both
sample sizes, i.e.,  $\mnoRi{\hxdfPr[\oDi{\ssY}]}{\xdf}{\iSv}\leq
\cst{\xdf,\edf} (\ssY\wedge\ssE)^{-1}$ for all $\ssY,\ssE\in\Nz$ and a
finite constant $\cst{\xdf,\edf}>0$ depending on $\xdf$ and $\edf$ only.
  We shall further
  emphasise that in case $\ssY=\ssE$ for any density $\xdf$ and $\edf$  it holds
  \begin{multline}\label{oo:e7}
    \moRa[\ssY]=\sum_{|j|\in\nset[1,]{\noDi}}|\fxdf[j]|^2[1\wedge
    \ssY^{-1}\iSv[j]]+\sum_{|j|>\noDi}|\fxdf[j]|^2[1\wedge
    \ssY^{-1}\iSv[j]]\\\leq
    \VnormLp{\ProjC[0]\xdf}^2 \oDi{\ssY}
    \oiSv[\oDi{\ssY}]/\ssY +
    \VnormLp{\ProjC[0]\xdf}^2\bias[\oDi{\ssY}]^2\leq
    2\VnormLp{\ProjC[0]\xdf}^2\doRao{\oDi{\ssY}}
  \end{multline}
  which in turn implies
  $\mnoRi{\hxdfPr[\oDi{\ssY}]}{\xdf}{\iSv}\leq\cst{}(1\vee\VnormLp{\ProjC[0]\xdf}^2)\noRao.$
  In other words, the estimation of the unknown error density $\edf$
  is negligible whenever $\ssY\leq\ssE$.\remEnd
\end{rem}

\begin{te}\paragraph{Aggregation.} Introducing aggregation weights
$\We[]$ consider an aggregation
  $\hxdf[{\We[]}]=\sum_{\Di\in\nset[1,]{\ssY}}\We\hxdfPr$
  of the orthogonal series estimators $\hxdfPr$, $\Di\in\Nz$, defined
  in \eqref{i:hxdfPr} with coefficients $\Zsuite{\fhxdfPr[{\We[]}]{j}}$ satisfying $\fhxdfPr[{\We[]}]{j}=0$ for
  $|j|>\ssY$, and  $\fhxdfPr[{\We[]}]{j}=\FuVg{\We[]}(\nset{|j|,\ssY})\hfedfmpI[j]\hfydf[j]$
 for any  $|j|\in\nset[1,]{\ssY}$.
  We note that again by construction $\fhxdfPr[{\We[]}]{0}=1$,
    $\fhxdfPr[{\We[]}]{-j}=\overline{\fhxdfPr[{\We[]}]{j}}$ and $1\geq |\fhxdfPr[{\We[]}]{j}|$. Hence,
    $\hxdf[{\We[]}]$ is real and integrates to one, however, it is not
    necessarily non-negative.  Our aim is to prove an upper bound for its risk
  $\mnoRi{\hxdfPr[{\We[]}]}{\xdf}{\iSv}$ and its maximal risk
  $\mnmRi{\hxdfPr[{\We[]}]}{\rwCxdf}{\rwCedf}$. Here and subsequently, we denote
  $\pxdfPr:=\sum_{j=-\Di}^{\Di}\hfedfmpI[j]\fydf[j]\bas_j=\sum_{j=-\Di}^{\Di}\hfedfmpI[j]\fedf[j]\fxdf[j]\bas_j$
  for   $\Di\in\Nz$. For
  arbitrary aggregation weights and penalties, the next lemma establishes an
  upper bound for the loss of the aggregated estimator. Selecting the weights and penalties suitably, it provides in the sequel our key argument.
\end{te}

\begin{lem}\label{co:agg:au}
  Consider an weighted sum  $\hxdf[{\We[]}]$ with arbitrary aggregation weights $\We[]$ and non-negative penalty terms
  $(\pen)_{\Di\in\nset[1,]{\ssY}}$.   For any
  $\mDi,\pDi\in\nset[1,]{\ssY}$ holds
  \begin{multline}\label{co:agg:au:e1}
    \VnormLp{\hxdf[{\We[]}]-\xdf}^2\leq
    3\VnormLp{\hxdfPr[\pDi]-\pxdfPr[\pDi]}^2
    +3 \VnormLp{\ProjC[0]\xdf}^2\bias[\mDi]^2(\xdf)
    \\\hfill
    +3 \VnormLp{\ProjC[0]\xdf}^2\FuVg{\We[]}(\nsetro[1,]{\mDi})
    +\tfrac{3}{7}\sum_{l\in\nsetlo{\pDi,\ssY}}\pen[l]\We[l]
    \Ind{\{\VnormLp{\hxdfPr[l]-\pxdfPr[l]}^2<\pen[l]\}}\\\hfill
    +3\sum_{l\in\nsetlo{\pDi,\ssY}}\vectp{\VnormLp{\hxdfPr[l]-\pxdfPr[l]}^2-\pen[l]/7}
    +\tfrac{3}{7}\sum_{l\in\nsetlo{\pDi,\ssY}}\pen[l]
    \Ind{\{\VnormLp{\hxdfPr[l]-\pxdfPr[l]}^2\geq\pen[l]/7\}}\\
    +6\sum_{j\in\nset[1,]{\ssY}}|\hfedfmpI[j]|^2|\fedf[j]-\hfedf[j]|^2|\fxdf[j]|^2
    +2\sum_{j\in\nset[1,]{\ssY}}\Ind{\{|\hfedf[j]|^2<1/\ssE\}}|\fxdf[j]|^2
 \end{multline}
\end{lem}

\begin{rem}The upper bound in \eqref{co:agg:au:e1} is similar to
  \eqref{co:agg:e1} apart of the last two terms which are controled again
  by \cref{oSv:re}. However, in oder to control
  the third and fourth term we replace in both weights and penalties
   the quantities  $\iSv$, which are
  not anymore known in advance, by their natural estimators.
  \remEnd
\end{rem}

\begin{te}We consider  either for some constant
    $\rWc\in\Nz$  \pcw
    $(\erWe)_{\Di\in\nset[1,]{\ssY}}$ as in \eqref{au:de:erWe} or \msw   $(\msWe)_{\Di\in\nset[1,]{\ssY}}$ as in
    \eqref{au:de:msWe}. Given $\eiSvS=\Ksuite[j\in\Nz_0]{\eiSv[j]}$, with
  \begin{equation}\label{eiSv}
    \eiSv[j]:=|\hfedfmpI[j]|^2 = |\hfedfmpI[-j]|^2 = \vert \hfedfI[j] \vert^{2}\Ind{\{|\hfedf[j]|^2\geq1/\ssE\}},
  \end{equation}
    $\LiSy{\eiSv}$ as in \eqref{ak:de:LiSy} with $x_{\mbullet}$
    replaced by $\eiSvS$ and a numerical  constant $\cpen>0$ we use
  \begin{equation}\label{au:de:peneSv}
    \peneSv:=\peneD,\quad\Di\in\Nz
  \end{equation}
  as penalty terms.\
\end{te}
\begin{te}\end{te}

\begin{thm}\label{au:ag:ub:pnp}Consider a weighted sum $\hxdf[{\We[]}] = \sum\nolimits_{\Di \in \nset{\ssY}} \We \hxdfPr$  using either
  \pcw $\We[]:=\erWe[]$ as in \eqref{au:de:erWe} or \msw $\We[]:=\msWe[]$ as in \eqref{au:de:msWe} and
  penalties $(\peneSv)_{\Di\in\nset[1,]{n}}$ as in
  \eqref{au:de:peneSv}  with numerical constant $\cpen\geq84$.  Let
  $\Di_{\ydf}:=\floor{3(400)^2\Vnormlp[1]{\fydf}^2}$ and for
    $\ssE\in\Nz$ set
    $\sDi{\ssE}:=\max\{\Di\in\nset[1,]{\ssE}:
    289\log(\Di+2)\cmiSv[\Di]\miSv[\Di]\leq\ssE\}$,
    if the defining set is not empty, and
    $\sDi{\ssE}:=\ceil{\Di_{\ydf}\log(2+\ssE)}$ otherwise. 
  \begin{Liste}[]
  \item[\mylabel{au:ag:ub:pnp:p}{\dgrau\bfseries{(p)}}]Assume there is
    $K\in\Nz_0$ with $1\geq \bFxdf[(K-1)]>0$ and
    $\bFxdf[K]=0$. If $K=0$ we set  $c_{\xdf}:=0$ and
    $c_{\xdf}:=\tfrac{104\cpen}
    {\VnormLp{\ProjC[0]\xdf}^2\sbFxdf[(K-1)]}$, otherwise. For
    $\ssY\in\Nz$ let
    $\sDi{\ssY}:=\max\{\Di\in\nset[1,]{\ssY}:\ssY>c_{\xdf}\DipenSv\}$,
    if the defining set is not empty, and
    $\sDi{\ssY}:=\ceil{\Di_{\ydf}\log(2+\ssY)}$ otherwise.
    There is a constant
    $\cst{\xdf\edf}$ given in \eqref{au:ag:ub:p:e1} depending only
    on $\xdf$ and $\edf$ such that for all $\ssY,\ssE\in\Nz$ holds%
    \begin{equation}\label{au:ag:ub:pnp:e1}
       \mnoRi{\hxdfAg[{\We[]}]}{\xdf}{\iSv}
       \leq
       \cst{}\VnormLp{\ProjC[0]\xdf}^2\big[(\ssY\wedge\ssE)^{-1} \vee
       \exp\big(\tfrac{-\cmiSv[{(\sDi{\ssY}\wedge\sDi{\ssE})}]
         (\sDi{\ssY}\wedge\sDi{\ssE})}{\Di_{\ydf}}\big)\big]\\
       + \cst{\xdf\edf}(\ssY\wedge\ssE)^{-1}.
     \end{equation}
   \item[\mylabel{au:ag:ub:pnp:np}{\dgrau\bfseries{(np)}}] Assume
     $\bFxdf>0$ for all $\Di\in\Nz$ and consider $ \noRaA$ as in
     \eqref{ak:ag:ub:pnp:e2}.
   There is a constant
    $\cst{\xdf\edf}$ given in \eqref{au:ag:ub:pnp:p13} depending only
    on $\xdf$ and $\edf$ such that for all $\ssY,\ssE\in\Nz$ holds
    \begin{multline}\label{au:ag:ub:pnp:e2}
     \mnoRi{\hxdfAg[{\We[]}]}{\xdf}{\iSv}
     \leq\cst{}\{[1\vee\VnormLp{\ProjC[0]\xdf}^2]\;\noRaA
    +\moRaA\}   + \cst{\xdf\edf}(\ssY\wedge\ssE)^{-1} \\\hfill
     \text{with }\quad \moRaA:=\moRa\vee\VnormLp{\ProjC[0]\xdf}^2\big[\bias[\sDi{\ssE}]^2(\xdf)
     \vee\exp\big(\tfrac{-\cmiSv[\sDi{\ssE}]\sDi{\ssE}}{\Di_{\ydf}}\big)\big].
  \end{multline}
\end{Liste}
\end{thm}

\begin{co}\label{au:ag:ub2:pnp}
  Let the assumptions of \cref{au:ag:ub:pnp} be satisfied and
    in addition
    \begin{inparaenum}
    \item[\mylabel{au:ag:ub2:pnp:pc:b}{{\normalfont\rmfamily\dgrau\bfseries(A4)}}]
            there is $\ssE_{\xdf\edf}\in\Nz$ such that
      $\cmiSv[\sDi{\ssE}]\sDi{\ssE}\geq \Di_{\ydf}\log\ssE$ for all
      $\ssE\geq \ssE_{\xdf\edf}$.
    \end{inparaenum}
  \begin{Liste}[]
  \item[\mylabel{au:ag:ub2:pnp:p}{\dgrau\bfseries{(p)}}]
    If \ref{ak:ag:ub2:pnp:pc} as in \cref{ak:ag:ub2:pnp} and \ref{au:ag:ub2:pnp:pc:b}
    hold true, then there is a constant $\cst{\xdf\edf}$ depending
    only on $\xdf$ and $\edf$ such that for all $\ssY,\ssE\in\Nz$ holds
    $\mnoRi{\hxdfAg[{\We[]}]}{\xdf}{\iSv} \leq
    \cst{\xdf\edf}(\ssY\wedge\ssE)^{-1}$.
  \item[\mylabel{au:ag:ub2:pnp:np}{\dgrau\bfseries{(np)}}]
    If  \ref{ak:ag:ub2:pnp:npc} as in \cref{ak:ag:ub2:pnp} and  \ref{au:ag:ub2:pnp:pc:b}
    hold true, then there is a constant $\cst{\xdf\edf}$ depending
    only on $\xdf$ and $\edf$ such that $\mnoRi{\hxdfAg[{\We[]}]}{\xdf}{\iSv}
    \leq \cst{\xdf\edf}\big(\noRaL+\moRa+\sbFxdf[\sDi{\ssE}]\big)$ for all $\ssY,\ssE\in\Nz$ holds true.
  \end{Liste}
\end{co}

\begin{il}\label{au:il:oo}
  Consider the cases \ref{il:edf:o} and \ref{il:edf:s} for the
  error density $\edf$ as in \cref{ak:il:oo}, where in both cases
 \cref{ak:ag:ub2:pnp} \ref{ak:ag:ub2:pnp:pc} holds true (cf. \cref{ak:il:oo}
 \ref{ak:il:edf:o} and \ref{ak:il:edf:s}). Moreover \cref{au:ag:ub2:pnp} \ref{au:ag:ub2:pnp:pc:b} is  satisfied, since
   \begin{inparaenum}[i]
  \item[\mylabel{au:il:edf:o}{\dg\bfseries{(o)}}]
    $\sDi{\ssE}\sim(\ssE/\log\ssE)^{1/(2a)}$ and
    $\sDi{\ssE}\cmSv[\sDi{\ssE}]\sim (\ssE/\log\ssE)^{1/(2a)}$, and
  \item[\mylabel{au:il:edf:s}{\dg\bfseries{(s)}}]
    $\sDi{\ssE}\sim(\log\ssE)^{1/(2a)}$
    and $\sDi{\ssE}\cmSv[\sDi{\ssE}]\sim (\log \ssE)^{2+1/(2a)}$.
  \end{inparaenum}
  Therefore, employing \cref{au:ag:ub2:pnp}
  \ref{au:ag:ub2:pnp:p} the fully data-driven aggregation  attains the
  parametric rate. For  \ref{au:ag:ub2:pnp:np} due to \cref{au:ag:ub2:pnp}
  the risk of the fully data-driven aggregated estimator has the order  $\noRaL+\moRa$,
  if \ref{ak:ag:ub2:pnp:npc} and
  \ref{au:ag:ub2:pnp:pc:b} are satisfied and
  $\bias[\sDi{\ssE}]^2(\xdf)$ is negligible with
  respect to $\moRa$.   The upper
  bound  $\moRaA$  in \cref{au:ag:ub:pnp} \eqref{au:ag:ub:pnp:e2}
  faces otherwise a detoriation compared to $\moRa$ which we illustrate
  considering  the cases
  in \cref{ak:il:oo}. Note that  the other upper bound  $\noRaA$ in \cref{au:ag:ub:pnp} \eqref{au:ag:ub:pnp:e2} already
  appears in    \cref{ak:ag:ub:pnp} and has been discussed in
  \cref{ak:il:oo}. Therefore, we
  state below the order of the additional term $\moRaA$ only.\\[2ex]
  \centerline{\begin{inparaenum}[i]
    \begin{tabular}{@{}lllll@{}}
      \toprule
      Order &$|\fxdf[\Di]|^2$ &$\iSv[\Di]$    &$\moRa$ & $\moRaA$ \\
      \midrule
      \item[\mylabel{au:oo:oo}{\dg\bfseries{[o-o]}}]
       & $\Di^{-2p-1}$&$\Di^{2a}$& $\left\{\begin{array}{l}\ssE^{-p/a}\\
          (\ssE/\log\ssE)^{-1},\\
                                             \ssE^{-1}\end{array}\right.$
            & $\begin{array}{@{}ll@{}}(\ssE/\log\ssE)^{-p/a}&:\text{\smaller[1]$p<a$},\\
          (\ssE/\log\ssE)^{-1}&:\text{\smaller[1]$p=a$}\\
          \ssE^{-1}&:\text{\smaller[1]$p>a$}.
        \end{array}$ \\
      \item[\mylabel{au:oo:os}{\dg\bfseries{[o-s]}}]
            & $\Di^{-2p-1}$ & $e^{\Di^{2a}}$ & $|\log \ssE|^{-p/a}$
                                                       & $|\log \ssE|^{-p/a}$\\
      \item[\mylabel{au:oo:so}{\dg\bfseries{[s-o]}}]
            &$\Di^{2p-1}e^{-\Di^{2p}}$& $\Di^{2a}$ & $\ssE^{-1}$& $\ssE^{-1}$\\
      \bottomrule
    \end{tabular}
  \end{inparaenum}}
  \\[2ex]
  Combining the \cref{ak:il:oo,au:il:oo} the fully data-driven aggregation
  attains the rate $\noRao+\moRa$   in case \ref{au:oo:os}, \ref{au:oo:oo} with $p\geq a$,   and \ref{au:oo:so} with $p\leq1/2$. In  case \ref{au:oo:oo} with $p< a$ and
  \ref{ak:oo:so} with $p>1/2$ its rate
   features a detoriation compared to $\noRao+\moRa$  by a logarithmic factor $(\log\ssE)^{p/a}$  and $(\log\ssY)^{(2a+1)(1-1/(2p))}$, respectively.\ilEnd
\end{il}

\begin{te}
\paragraph{Minimax optimality.}
For  $\ssE\in\Nz$ setting
$\mmRa:=\supF{\xdfCw[j](1\wedge \edfCw[j]/\ssE)}$
 it holds $\moRa\leq   \edfCr\xdfCr \mmRa$
  for all $\xdf\in\rwCxdf$ and $\edf\in\rwCxdf$.
Exploiting again the upper bound
\eqref{oo:e6} and  the definition
\eqref{oo:de:doRao}
for all
 $\Di,\ssY,\ssE\in\Nz$ follows immediately
\begin{equation*}
  \mnmRi{\hxdfPr}{\rwCxdf}{\rwCedf}
  \leq \cst{}\,\xdfCr\edfCr\,(\dmRao{\Di}+\mmRa).
\end{equation*}
The upper bound in the last display depends on the dimension parameter
$\Di$ and hence by choosing an optimal value $\nmDio$ as in
\eqref{oo:de:doRao} the upper bound
will be minimised and it holds
$ \mnmRi{\hxdfPr[\nmDio]}{\rwCxdf}{\rwCedf}
\leq \cst{}\xdfCr\edfCr\,(\nmRao+\mmRa)$. On the other
  hand  \cite{JohannesSchwarz2013a} have shown  that for all
 $\ssY,\ssE\in\Nz$
  \begin{equation*}
  \inf\nolimits_{\hxdf}\nmRi{\hxdf}{\rwCxdf}{\rwCedf}\geq \cst{}\mmRa,
\end{equation*}
 where $\cst{}>0$  and the infimum is taken over all
  possible estimators $\hxdf$ of $\xdf$.  Consequently, combining
  \eqref{mm:lbn} and  the
  last lower bound   $\Nsuite[\ssY,\ssE]{\nmRao+\mmRa}$,  $\Nsuite[\ssY]{\nmDio}$
  and  $\Nsuite[\ssY]{\hxdfPr[\nmDio]}$, respectively, is a
  minimax rate, a minimax dimension and minimax optimal estimator (up to a
  constant).
\end{te}

\begin{te}\paragraph{Aggregation.}
    By applying
  \cref{co:agg:au} we derive bounds for the maximal risk defined
  of the fully data-driven aggregation.
  \end{te}

\begin{thm}\label{au:mrb:ag:ub:pnp}
  Consider an aggregation
  $\hxdf[{\We[]}]$ using either \pcw $\We[]:=\erWe[]$
  as in \eqref{au:de:erWe} or \msw $\We[]:=\msWe[]$
  as in \eqref{au:de:msWe} and penalties $(\peneSv)_{\Di\in\nset[1,]{\ssY}}$ as in \eqref{au:de:peneSv}  with numerical constant $\cpen\geq84$.
   Let \ref{ass:Cw} be satisfied and $\nmRaA$ as in \eqref{ak:ag:ub:pnp:mm:e1}. Set
   $\Di_{\xdfCw[]\edfCw[]}:=\floor{3(400)^2\xdfCr\zeta_{\edfCr}\,\Vnormlp[1]{\xdfCwS/\edfCwS}}$
    and for
    $\ssE\in\Nz$,
    $\sDi{\ssE}:=\max\{\Di\in\nset[1,]{\ssE}:
    289\log(\Di+2)\zeta_{\edfCr}\liCw[\Di]\edfCw[\Di]\leq\ssE\}$,
    if the defining set is not empty, and
    $\sDi{\ssE}:=\ceil{\Di_{\xdfCw[]\edfCw[]}\log(2+\ssE)}$ otherwise.
Then there is a constant $\cst[\xdfCr\edfCr]{\xdfCw[]\edfCw[]}$ given in
\eqref{au:mrb:ag:ub:pnp:p4} depending only on the classes $\rwCxdf$  and $\rwCedf$ such that for all
$\ssY,\ssE\in\Nz$
  \begin{multline}\label{au:mrb:ag:ub:pnp:e1}
    \mnmRi{\hxdfAg[{\We[]}]}{\rwCxdf}{\rwCedf}\leq
    \cst{}(\xdfCr+\zeta_{\edfCr})\big(\nmRaA + \mmRaA\big) + \cst[\xdfCr\edfCr]{\xdfCw[]\edfCw[]}(\ssY\wedge\ssE)^{-1} \\\hfill     \text{with }\quad \mmRaA:=\mmRa\vee\xdfCw[\sDi{\ssE}]\vee\exp\big(\tfrac{-\liCw[\sDi{\ssE}]\sDi{\ssE}}{\Di_{\xdfCw[]\edfCw[]}}\big)\big].
  \end{multline}
\end{thm}

\begin{co}\label{au:mrb:ag:ub2:pnp}
  Let the assumptions of \cref{au:mrb:ag:ub:pnp} be satisfied.  If
  \ref{ak:ag:ub2:pnp:mm:c} as in \cref{ak:ag:ub2:pnp:mm} and in
  addition
  \begin{inparaenum}
  \item[\mylabel{au:mrb:ag:ub2:pnp:pc:b}{{\normalfont\rmfamily\dgrau\bfseries(A4')}}]
    there is $\ssE_{\xdfCw[]\edfCw[]}\in\Nz$ such that
    $\liCw[\sDi{\ssE}]\sDi{\ssE}\geq \Di_{\xdfCw[]\edfCw[]}\log\ssE$ for all
    $\ssE\geq \ssE_{\xdfCw[]\edfCw[]}$,
  \end{inparaenum}
  then there  is a constant $\cst[\xdfCr\edfCr]{\xdfCw[]\edfCw[]}$
  depending only on the classes $\rwCxdf$  and $\rwCedf$ such that   for
  all $\ssY,\ssE\in\Nz$ holds
  $\mnmRi{\hxdfAg[{\We[]}]}{\rwCxdf}{\rwCedf} \leq
  \cst[\xdfCr\edfCr]{\xdfCw[]\edfCw[]}\big(\nmRaL+\mmRa+\xdfCw[\sDi{\ssE}]\big)$.
\end{co}

\begin{il}\label{au:il:mm}We have
  discussed the order of $\nmRaA$ appearing in \cref{au:mrb:ag:ub:pnp}
  for typical configurations in \cref{ak:il:mm}, thus we
  state below the order of the additional term only.\\[2ex]
  \centerline{\begin{inparaenum}[i]
    \begin{tabular}{@{}lllll@{}}
      \toprule
      Order &$\xdfCw$ &$\edfCw$    &$\mmRa$ & $\mmRaA$ \\
      \midrule
      \item[\mylabel{au:mm:oo}{\dg\bfseries{[o-o]}}]
       & $\Di^{-2p}$&$\Di^{2a}$& $\left\{\begin{array}{l}\ssE^{-p/a}\\
                                             \ssE^{-1}\end{array}\right.$
            &
              $\begin{array}{@{}ll@{}}(\ssE/\log\ssE)^{-p/a}&:\text{\smaller[1]$p\leq a$},\\
          \ssE^{-1}&:\text{\smaller[1]$p>a$}.
        \end{array}$ \\
      \item[\mylabel{au:mm:os}{\dg\bfseries{[o-s]}}]
            & $\Di^{-2p}$ & $e^{\Di^{2a}}$ & $|\log \ssE|^{-p/a}$
                                                       & $|\log \ssE|^{-p/a}$\\
      \item[\mylabel{au:mm:so}{\dg\bfseries{[s-o]}}]
            &$e^{-\Di^{2p}}$& $\Di^{2a}$ & $\ssE^{-1}$& $\ssE^{-1}$\\
      \bottomrule
    \end{tabular}
  \end{inparaenum}}\\[2ex]
  Note that in all cases the additional assumption  \ref{au:mrb:ag:ub2:pnp:pc:b} in \cref{au:mrb:ag:ub2:pnp} is satisfied (as in
  \cref{au:il:oo}), and hence $\mmRaA$ is of order
  $\mmRa+\xdfCw[\sDi{\ssE}]$. Moreover,  in case \ref{au:mm:oo},
  \ref{au:mm:os} and \ref{au:mm:so} holds $\xdfCw[\sDi{\ssE}]\sim(\ssE/\log\ssE)^{-p/a}$,
  $\xdfCw[\sDi{\ssE}]\sim(\log \ssE)^{-p/a}$ and
   $\xdfCw[\sDi{\ssE}]\sim\exp(-(\ssE/\log \ssE)^{p/a})$,
   respectively. Consequently, $\xdfCw[\sDi{\ssE}]$ is negligible compared
   to $\mmRa$ in  case \ref{au:mm:os} and \ref{au:mm:so}, but in \ref{au:mm:oo} for $p> a$ only.
   Combining \cref{ak:il:mm,au:il:mm} the fully data-driven aggregation
  attains the minimax  rate in case \ref{au:mm:os} with $p>0$, \ref{au:mm:oo} with $p>a$ and \ref{au:mm:so} with
  $p\leq1/2$, while in  case \ref{au:mm:oo} with $p\leq a$ and
  \ref{au:mm:so} with $p>1/2$  its rate features a detoriation
   by a logarithmic factor
  $(\log\ssE)^{p/a}$  and $(\log\ssY)^{(2a+1)(1-1/(2p))}$,
  respectively, compared to the minimax  rate.\ilEnd
\end{il}

\section{Simulation study}\label{si}

\begin{figure}
\centering
\begin{subfigure}{.45\textwidth}
  \centering
  \includegraphics[width=\linewidth]{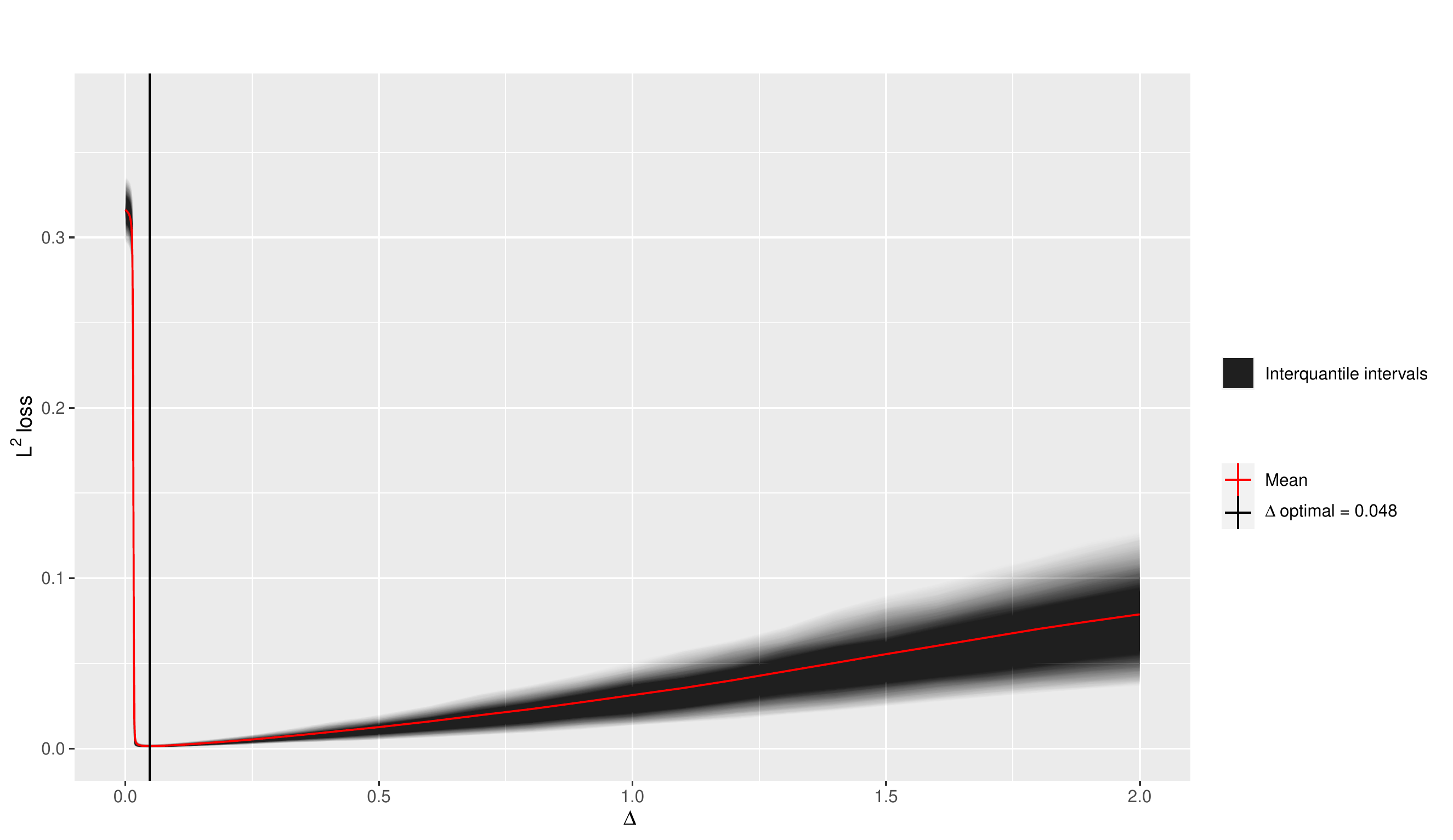}
  \label{si:delta:erWe}
\end{subfigure}%
\begin{subfigure}{.45\textwidth}
  \centering
  \includegraphics[width=\linewidth]{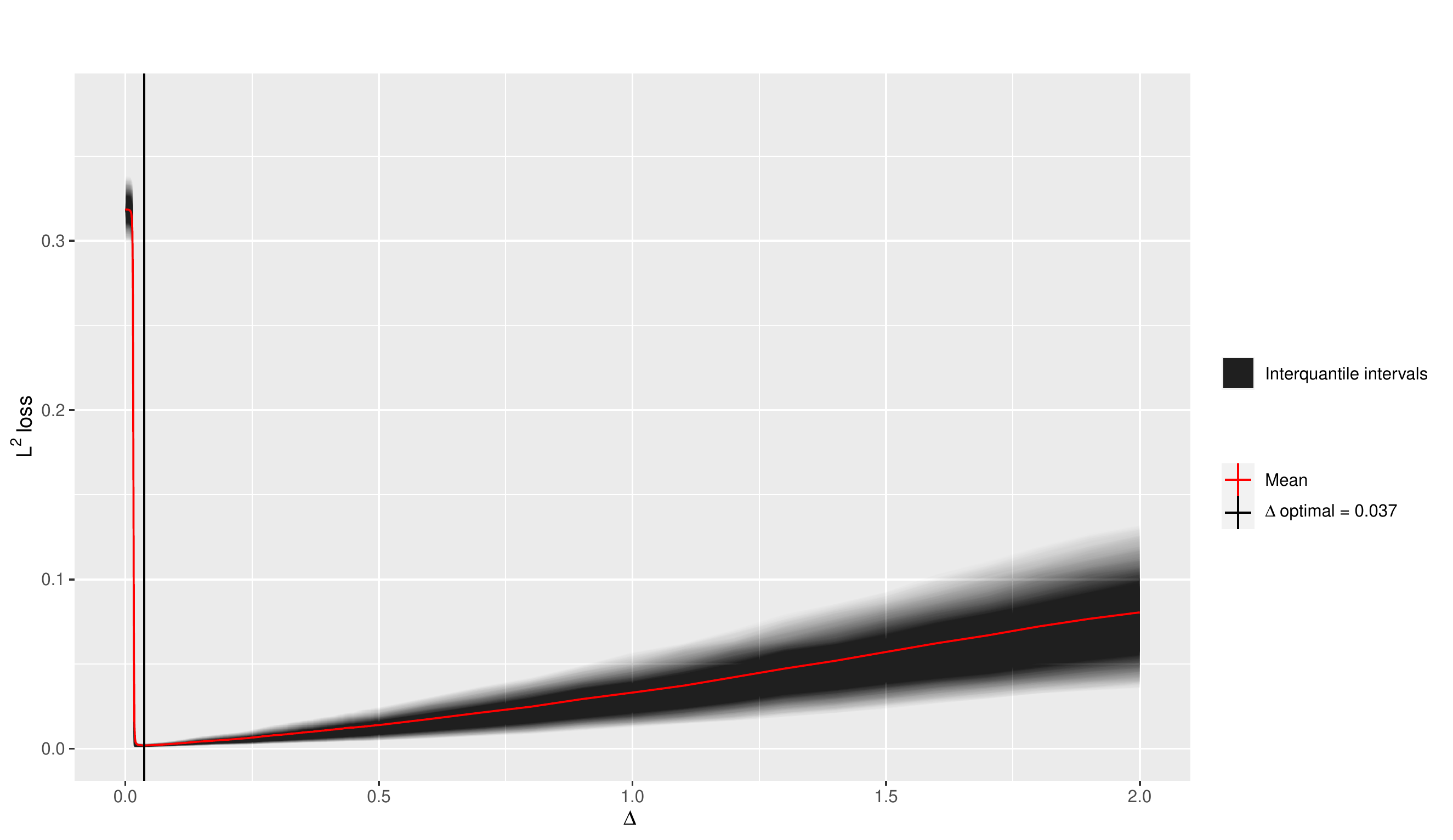}
  \label{si:delta:msWe}
\end{subfigure}
\caption{Empirical Bayesian risk for Bayesian (left) and model selection (right) weights over $1000$ replicates
as a function of the constant $\cpen$ and a minimal value (black vertical)}
\label{si:tune:res}
\end{figure}

\begin{figure}
\centering
\begin{subfigure}{.5\textwidth}
  \centering
  \includegraphics[width=\linewidth]{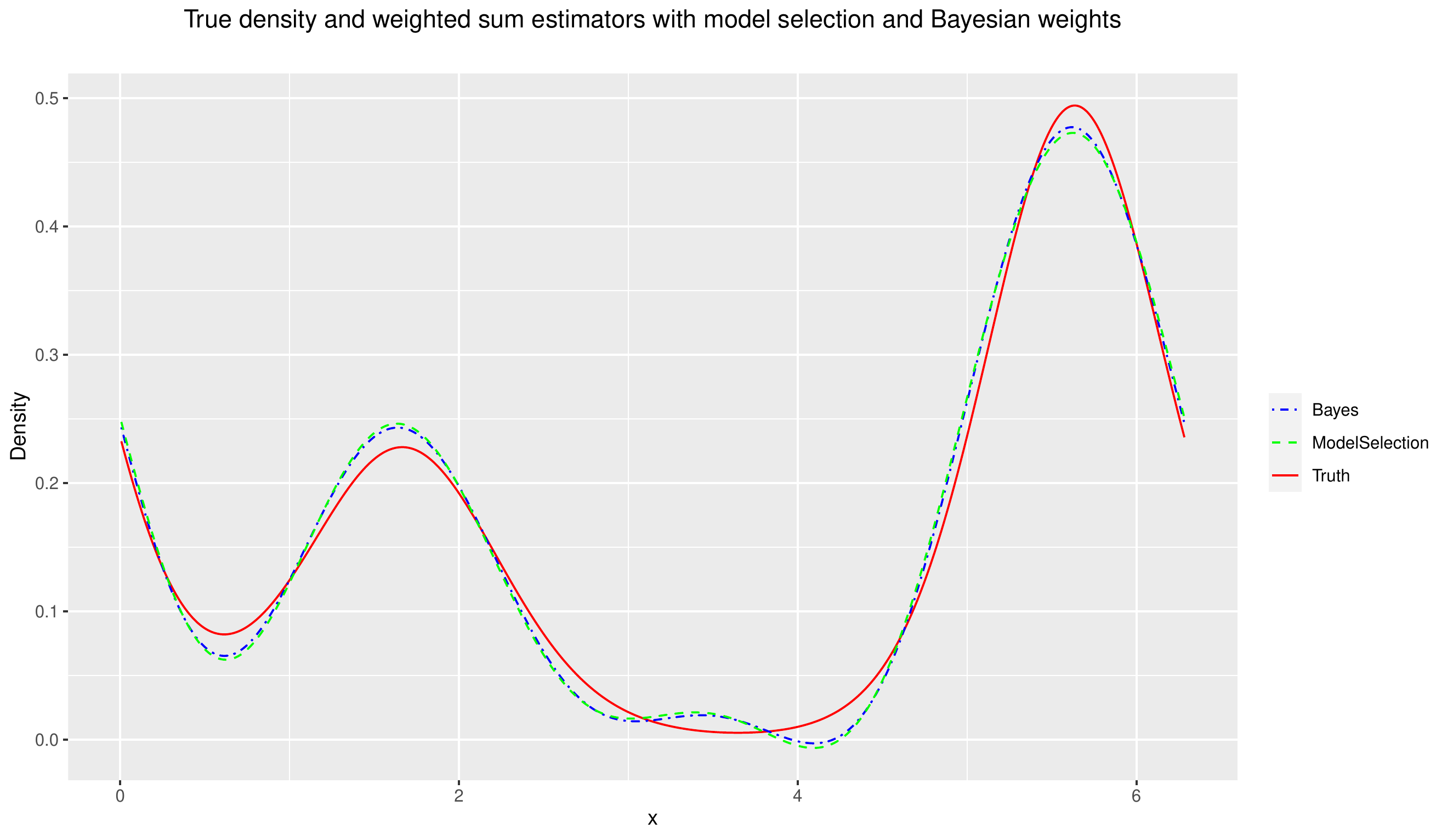}
  \caption{}
  \label{si:il:ex}
\end{subfigure}%
\begin{subfigure}{.5\textwidth}
  \centering
  \includegraphics[width=\linewidth]{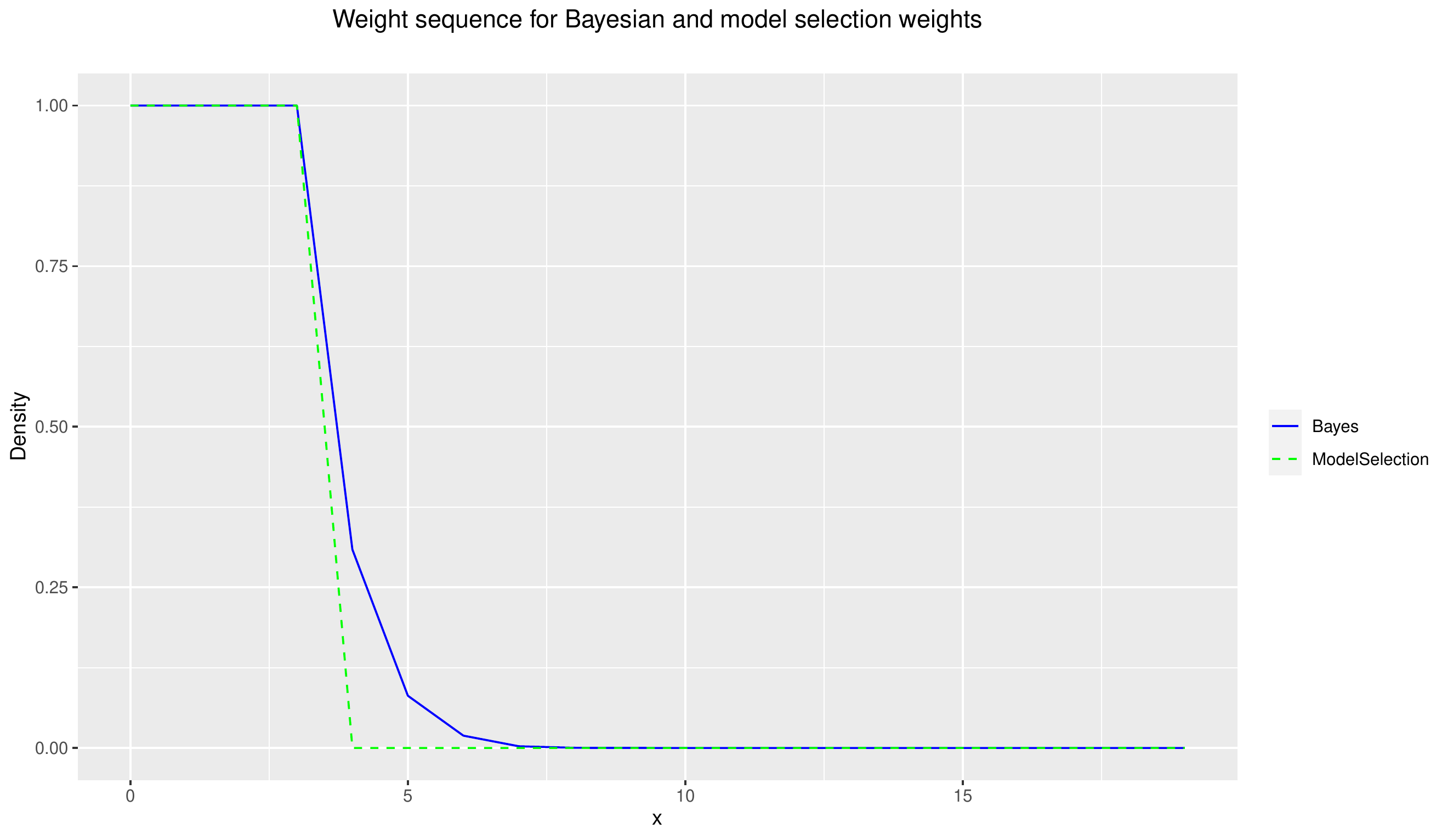}
  \caption{}
  \label{si:il:exWe}
\end{subfigure}
\caption{Weighted sum estimators using Bayesian weights or model selection weights (a), and the associated random weights (b) with a sample of size $n = m = 1000$}
\label{si:il}
\end{figure}

\begin{te}Let us illustrate the performance of the fully data-driven weighted sum of OSE’s either with  model
  selection  \eqref{au:de:msWe} or  Bayesian
   \eqref{au:de:erWe}  ($\rWc = 1$)  weights or by a simulation study.  As a first step, we
  calibrate the constant $\cpen$ appearing in the penalty  \eqref{au:de:peneSv}.  Indeed,
  \cref{au:ag:ub:pnp,au:mrb:ag:ub:pnp} stipulate that any choice $\cpen \geq 84$ ensures optimal rates but this
  is not a necessary condition and the constant obtained this way is often too large.  Hence, we select a value
  minimising a Bayesian empirical risk obtained by repeating $1000$ times the procedure as described
  hereafter.  We randomly pick a noise density $\edf$ and a density of interest $\xdf$, respectively, from a
  family of wrapped asymmetric Laplace distributions and a family of wrapped normal distributions. For the
  noise density the location parameter is uniformly-distributed in $[0, 1]$ and both the scale, and the
  asymmetry parameter follow a $\Gamma$ distribution with shape $0.5$ and scale $1$. For the density of interest the
  mean is again uniformly-distributed in $[0, 1]$, while the standard deviation has a $\Gamma$ distribution with
  shape $9$ and scale $0.5$. Next we generate a sample $(\rE_i)_{i\in\nset{\ssE}}$
  of size $\ssE = 5000$ from $\edf$, and a sample $(\rY_i)_{i\in\nset{\ssY}}$ of
  size $\ssY = 500$ from $\ydf = \xdf \oast \edf$. We use them to construct the estimators $\hxdfAg$, and
  $\hxdfAg[{\msWe[]}]$  as in \cref{au:ag:ub:pnp} for a  range of values of $\cpen$.  Finally, we compute and store
  the $\Lp[2]$-loss of each estimator obtained this way. Given the result of the  $1000$  repetitions,  for each
  value of the constant $\cpen$ we use the sample of $\Lp[2]$-losses to compute estimators of the mean squared
  error and  the quantiles of the distribution of the $\Lp[2]$-loss.
  Finally, we select and fix from now on  a value of $\cpen$ that minimises the empirical mean squared error.
  The results of this procedure are reported in \cref{si:tune:res}.  Using the calibrated constants and
  samples of size $\ssY = \ssE = 1000$   in \cref{si:il} we depict a realisation of the weighted sum estimators
  with Bayesian or model selection weights. In this example, $\xdf$ is a mixture of two von Mises distributions and $\edf$ is a wrapped asymetric Laplace distribution.
  The two estimators estimate the true density properly and behave similarly, we investigate next if there can be a significant performance difference.  \end{te}

\begin{figure}
\centering
\begin{subfigure}{.3\textwidth}
  \centering
  \includegraphics[width=\linewidth]{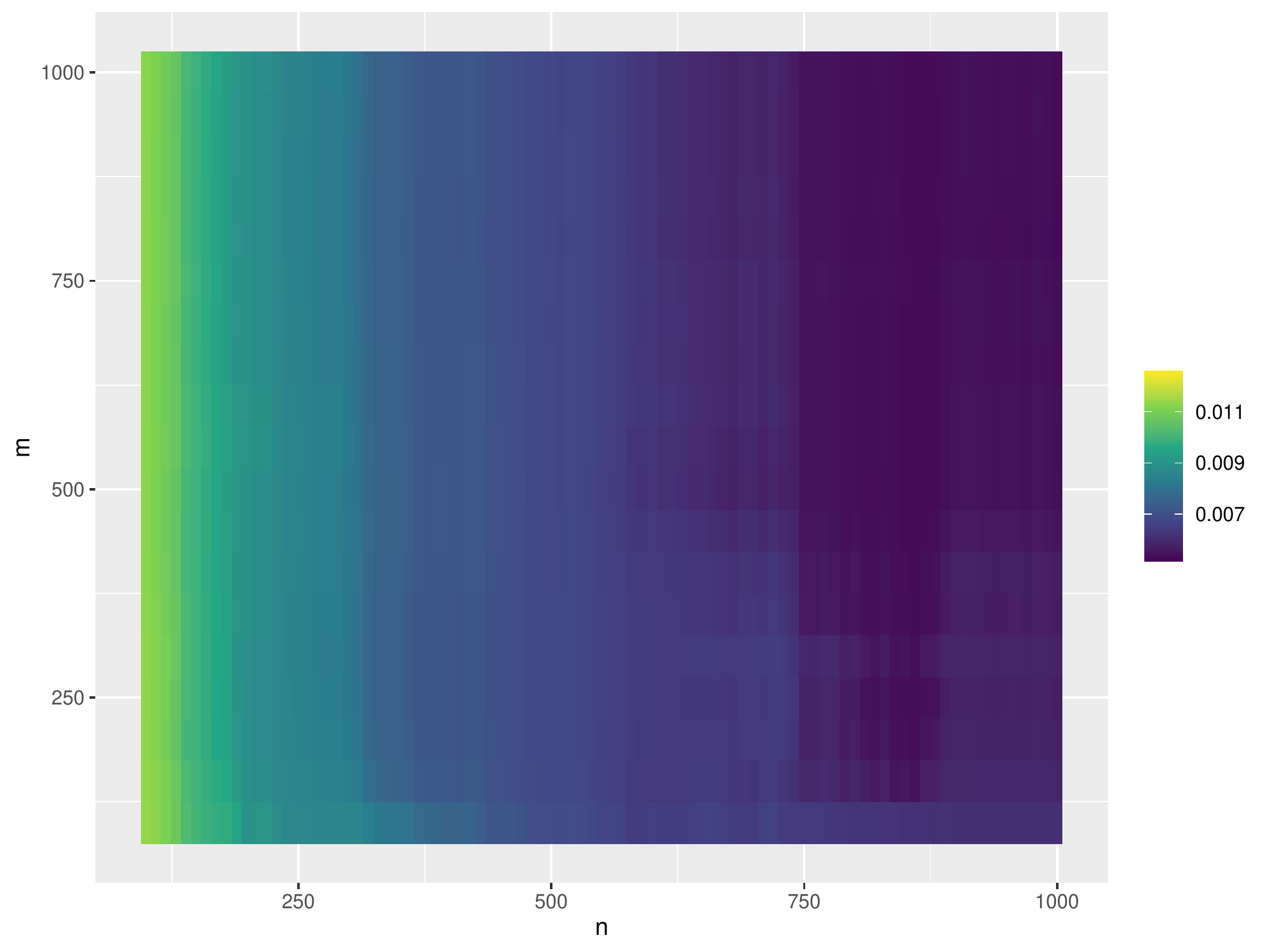}
  \caption{}
  \label{si:matrices:ag}
\end{subfigure}%
\begin{subfigure}{.3\textwidth}
  \centering
  \includegraphics[width=\linewidth]{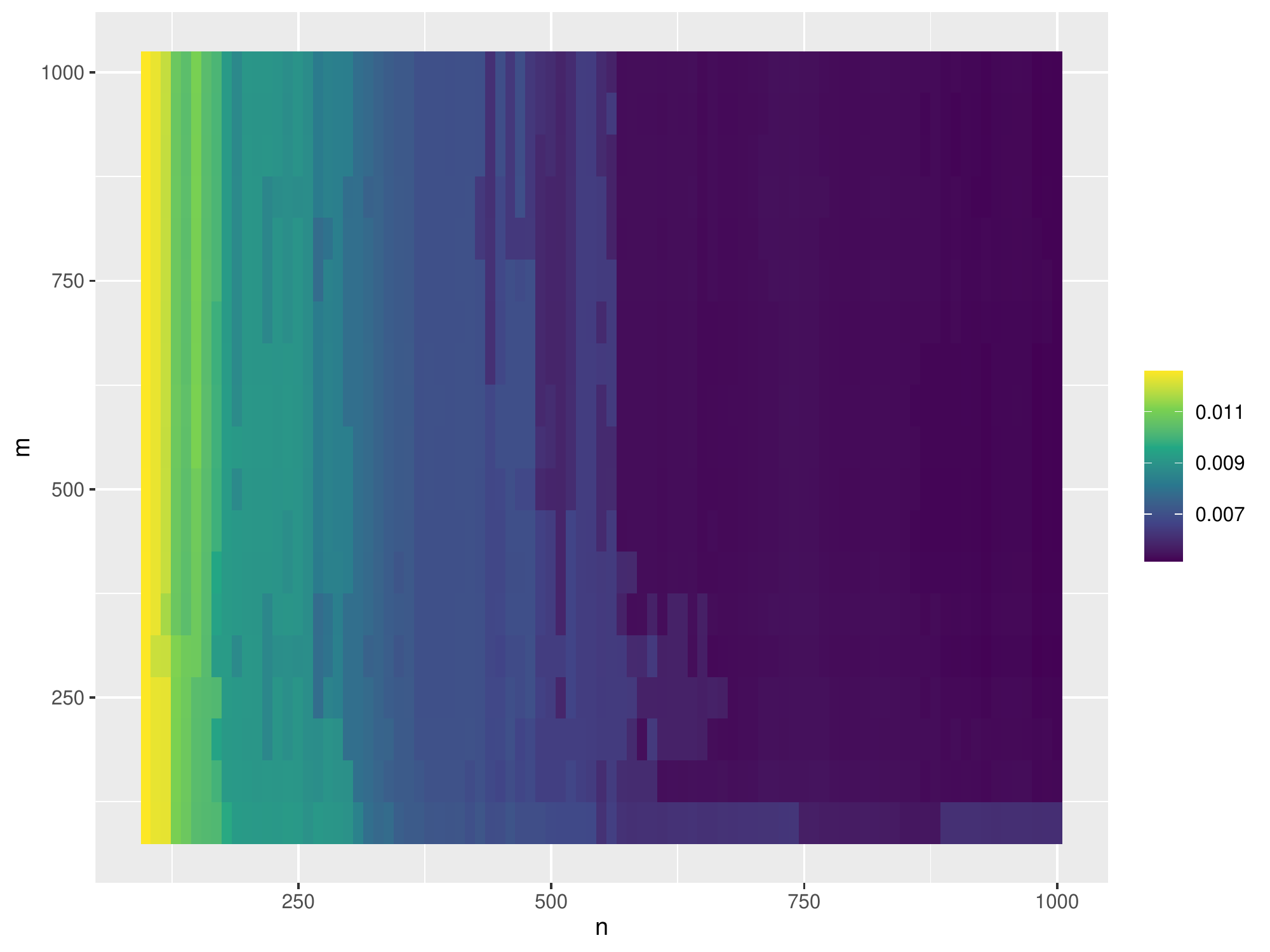}
  \caption{}
  \label{si:matrices:ms}
\end{subfigure}%
\begin{subfigure}{.3\textwidth}
  \centering
  \includegraphics[width=\linewidth]{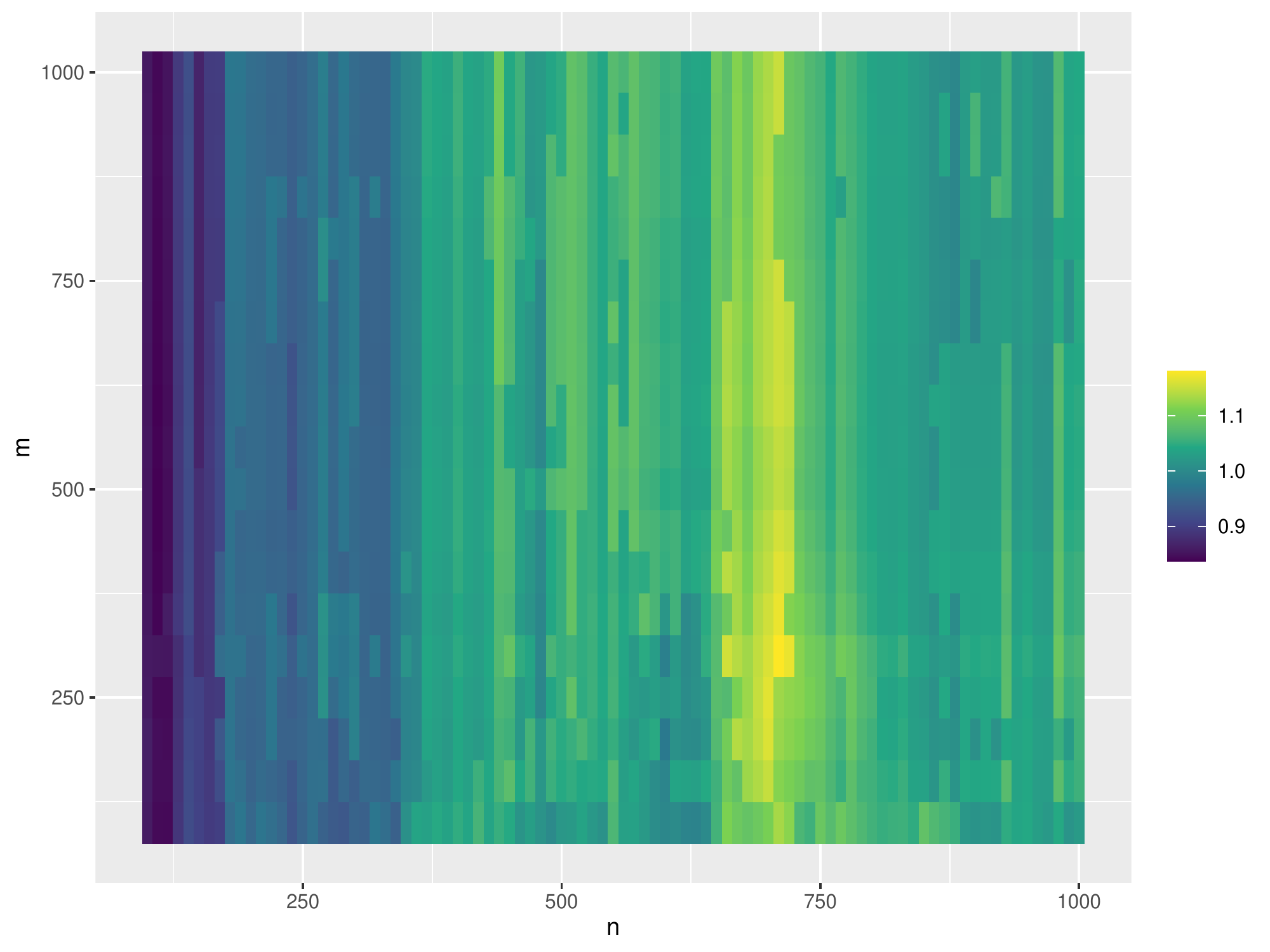}
  \caption{}
  \label{si:matrices:ra}
\end{subfigure}%
\caption{Empirical Bayesian risk for Bayesian (a) or model selection (b) weights and their ratio (c) over
  $100$ replicates as a function of the sample size $n$ (abscissa) and $m$ (ordinate).}
\label{si:matrices}
\end{figure}

\begin{figure}
\centering
\begin{subfigure}{.3\textwidth}
  \centering
  \includegraphics[width=\linewidth]{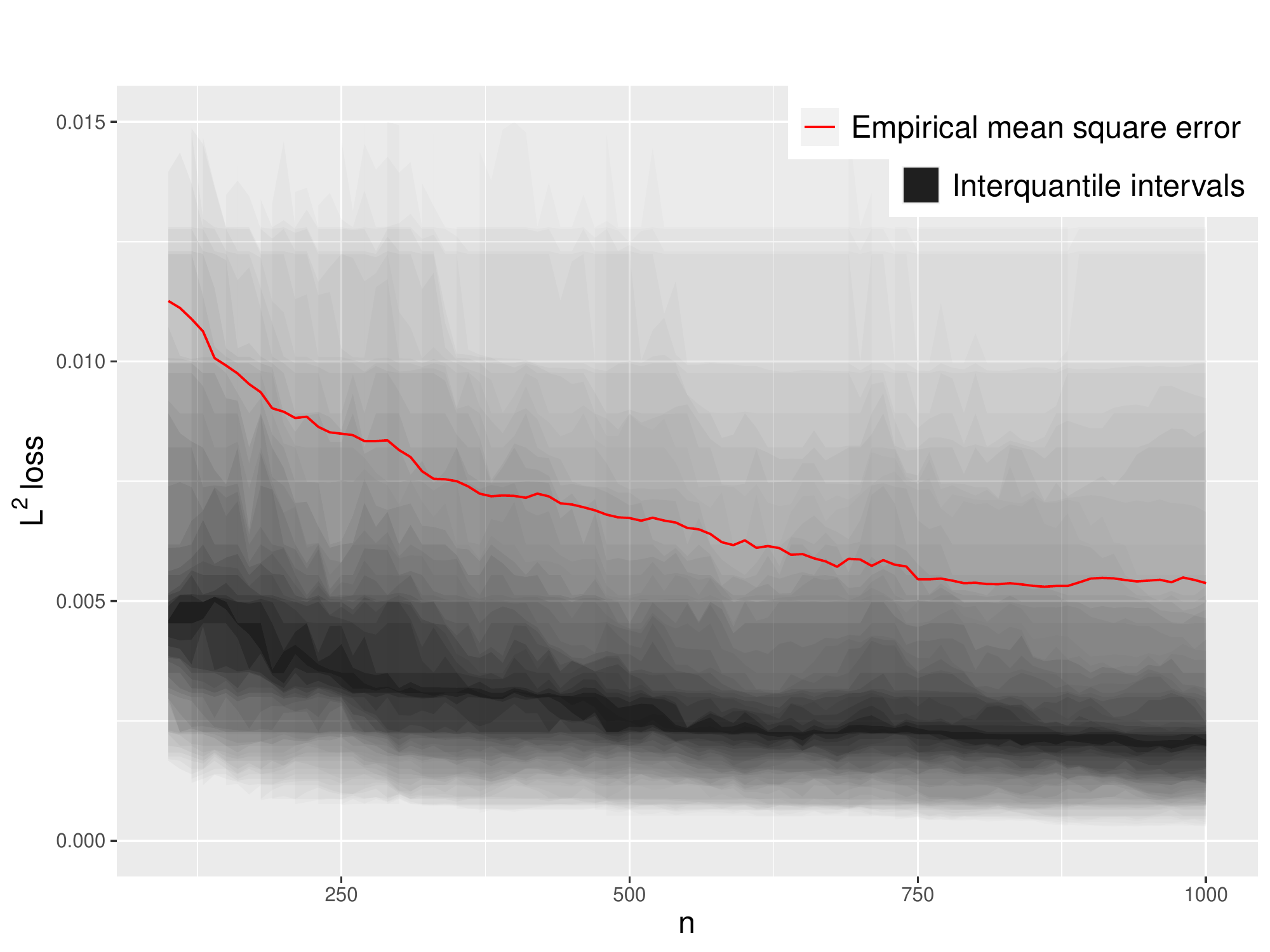}
  \caption{}
  \label{si:rate:agn}
\end{subfigure}%
\begin{subfigure}{.3\textwidth}
  \centering
  \includegraphics[width=\linewidth]{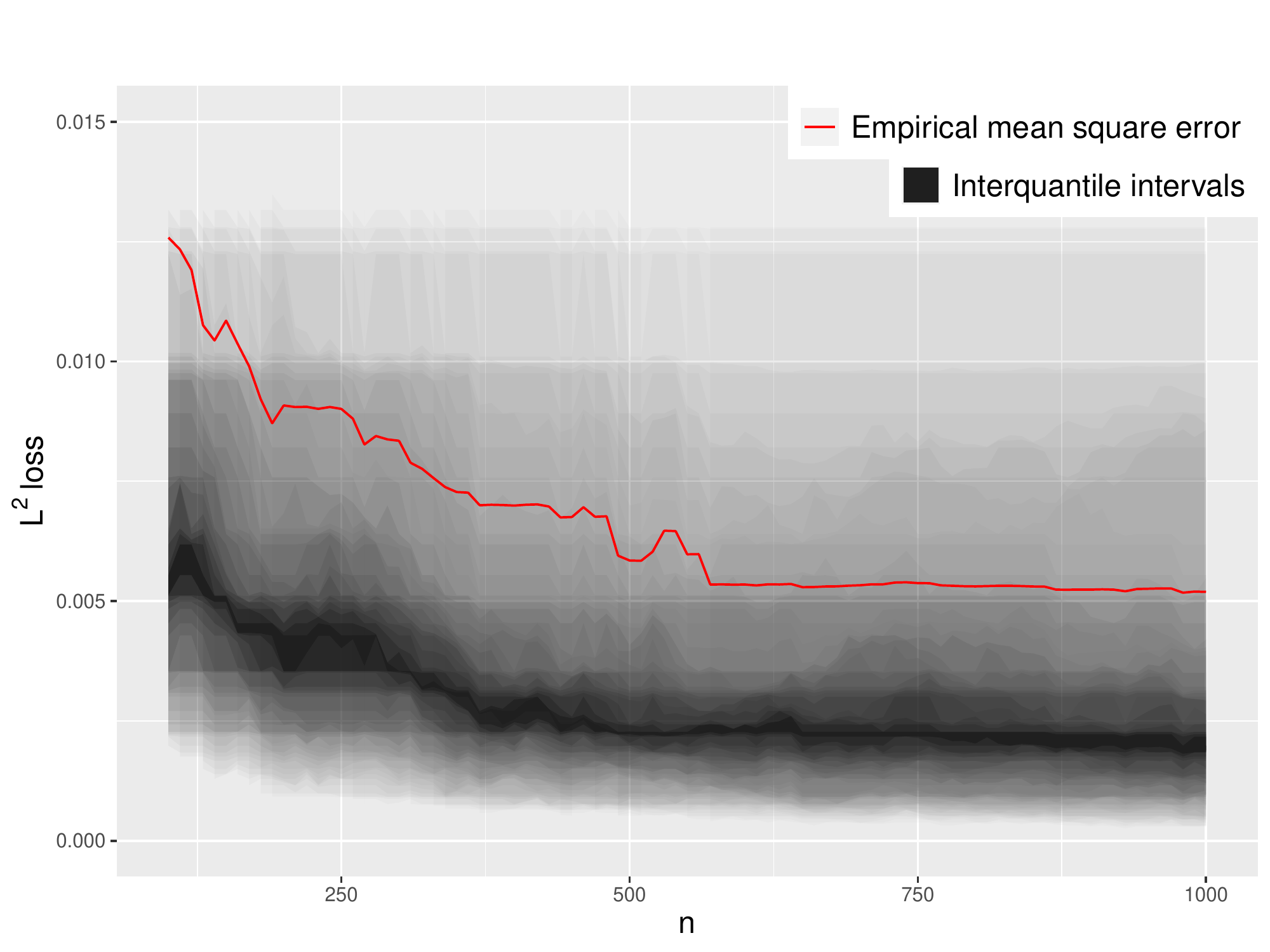}
  \caption{}
  \label{si:rate:msn}
\end{subfigure}%
\begin{subfigure}{.3\textwidth}
  \centering
  \includegraphics[width=\linewidth]{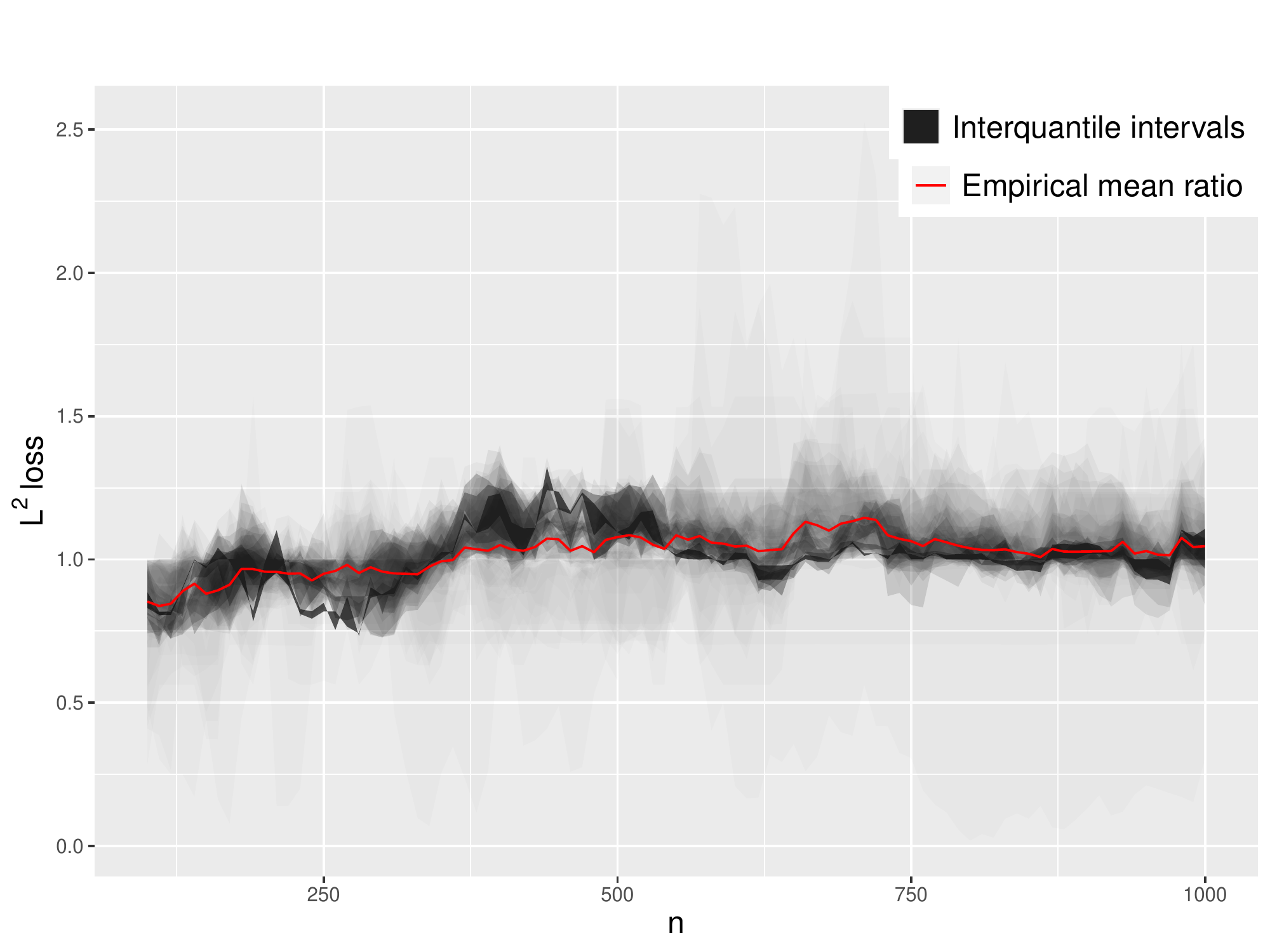}
  \caption{}
  \label{si:rate:ran}
\end{subfigure}%
\\
\begin{subfigure}{.3\textwidth}
  \centering
  \includegraphics[width=\linewidth]{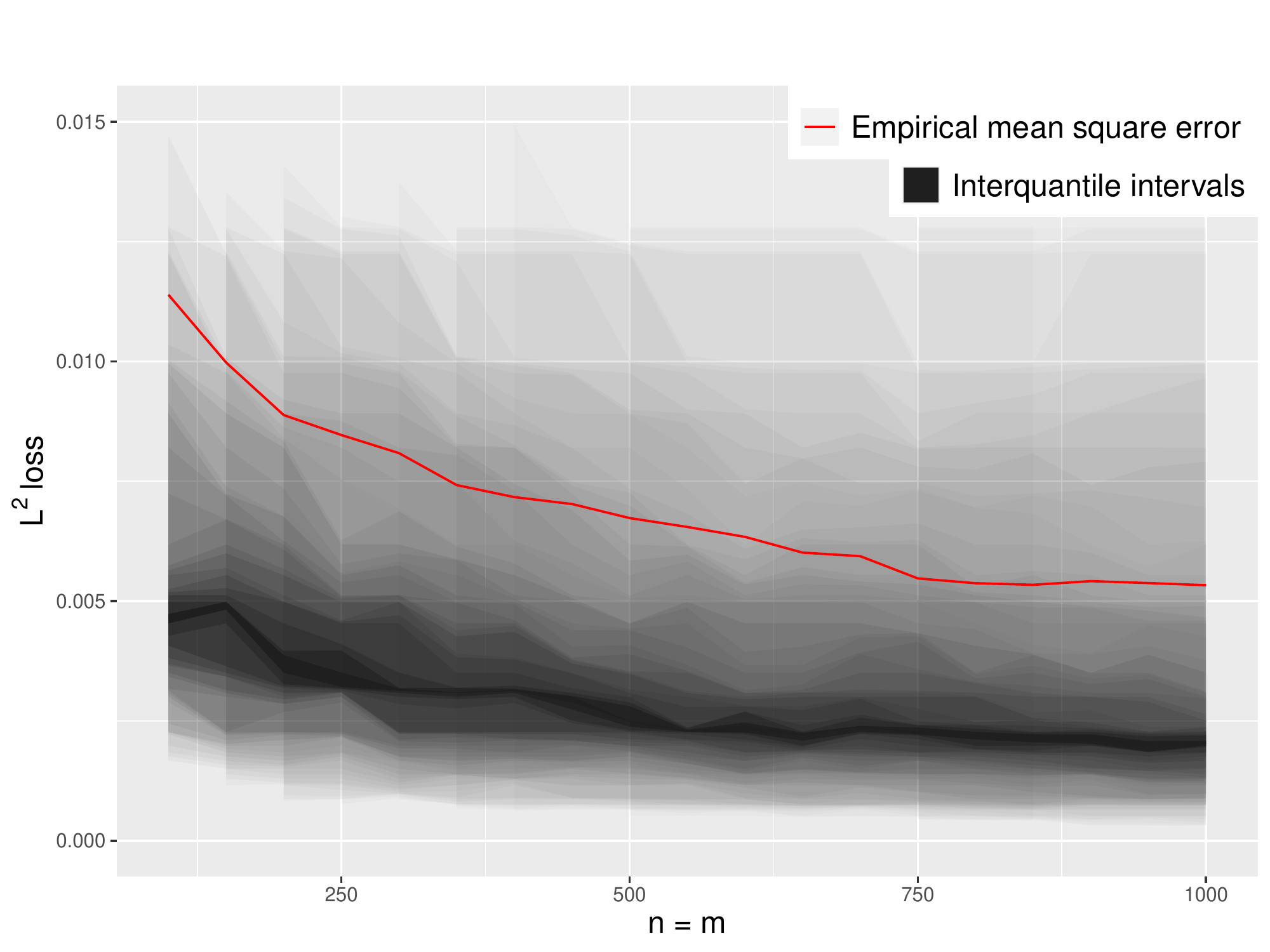}
  \caption{}
  \label{si:rate:agmn}
\end{subfigure}%
\begin{subfigure}{.3\textwidth}
  \centering
  \includegraphics[width=\linewidth]{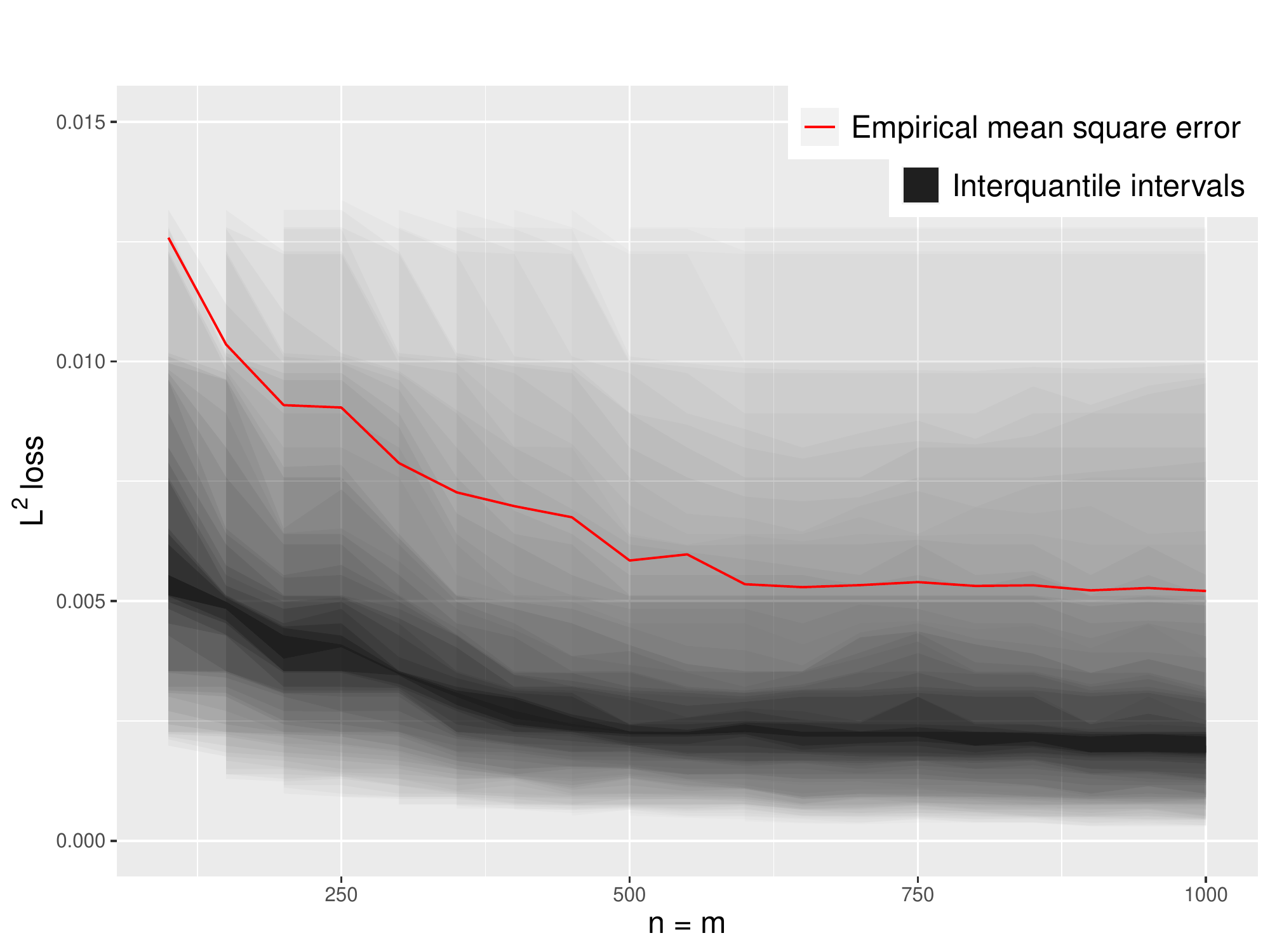}
  \caption{}
  \label{si:rate:msmn}
\end{subfigure}%
\begin{subfigure}{.3\textwidth}
  \centering
  \includegraphics[width=\linewidth]{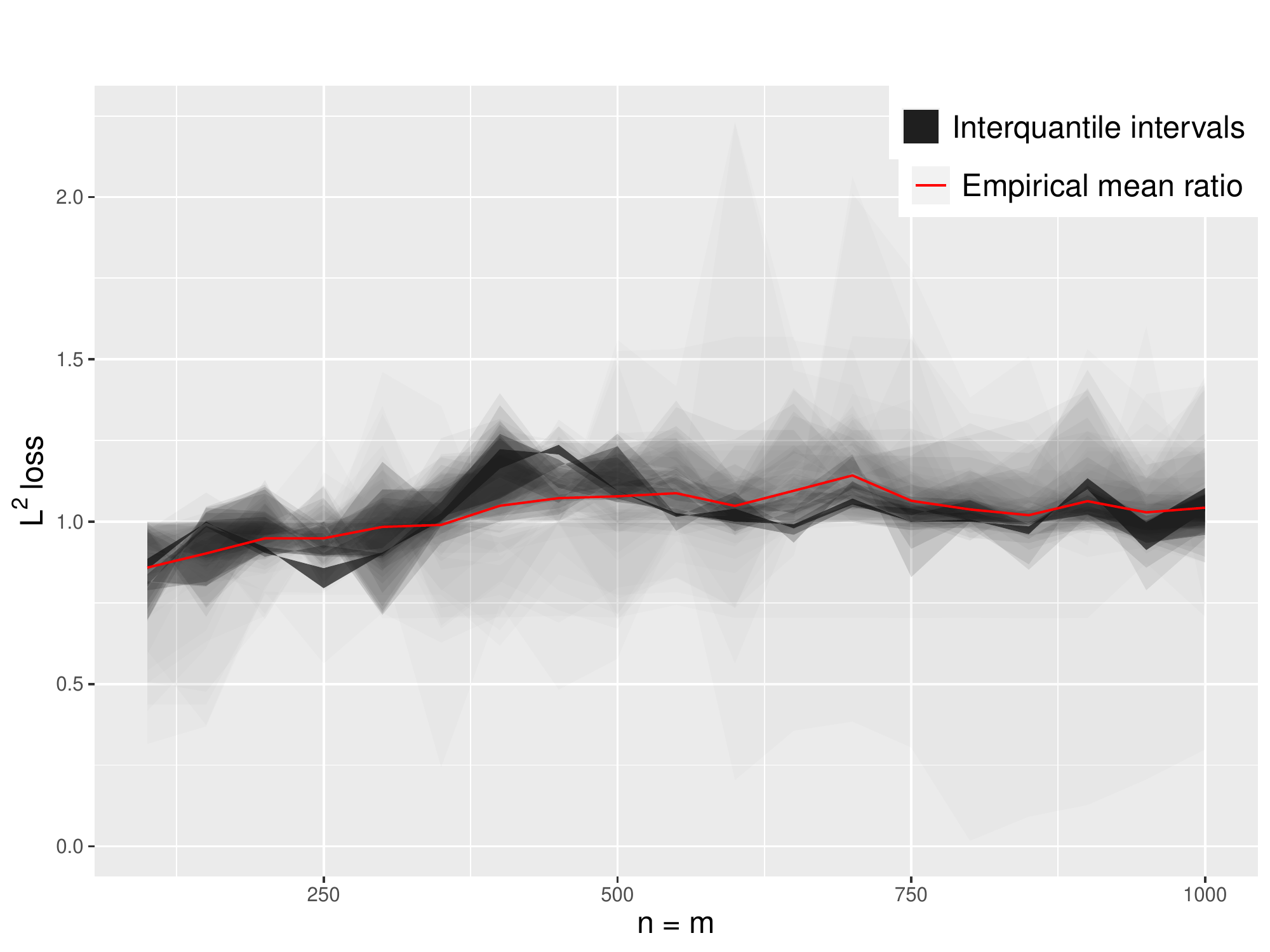}
  \caption{}
  \label{si:rate:ramn}
\end{subfigure}%
\\
\begin{subfigure}{.3\textwidth}
  \centering
  \includegraphics[width=\linewidth]{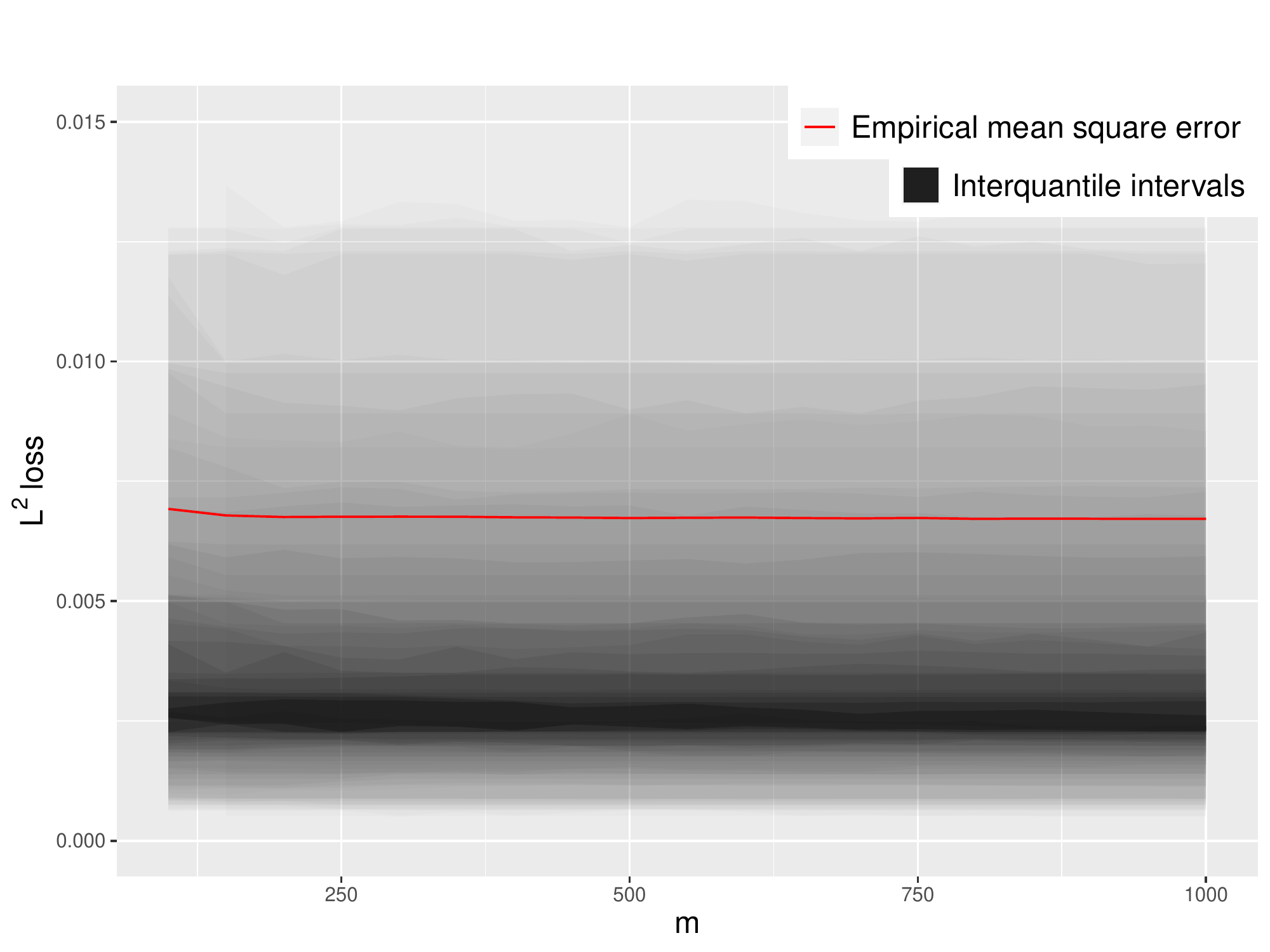}
  \caption{}
  \label{si:rate:agm}
\end{subfigure}%
\begin{subfigure}{.3\textwidth}
  \centering
  \includegraphics[width=\linewidth]{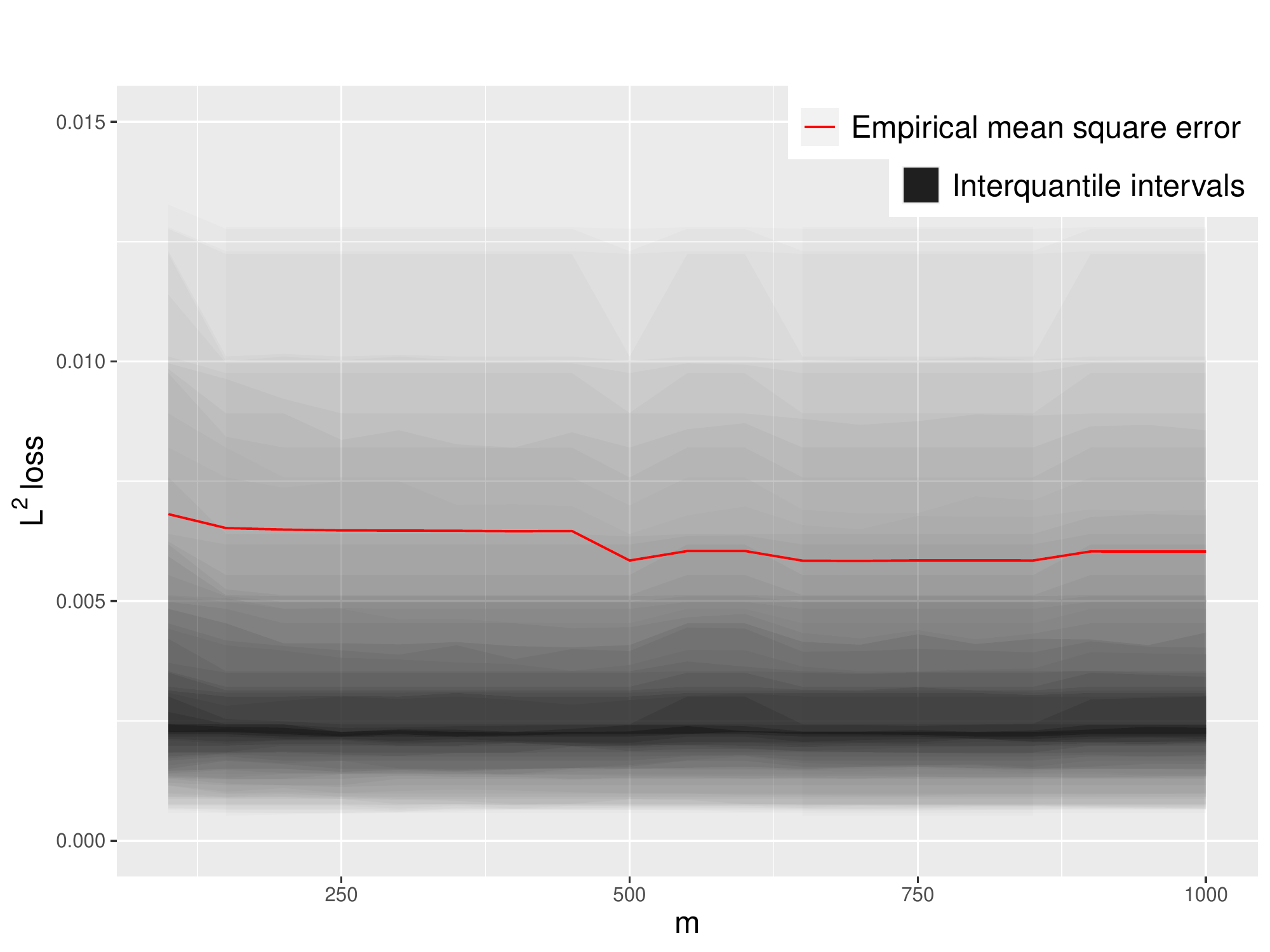}
  \caption{}
  \label{si:rate:msm}
\end{subfigure}%
\begin{subfigure}{.3\textwidth}
  \centering
  \includegraphics[width=\linewidth]{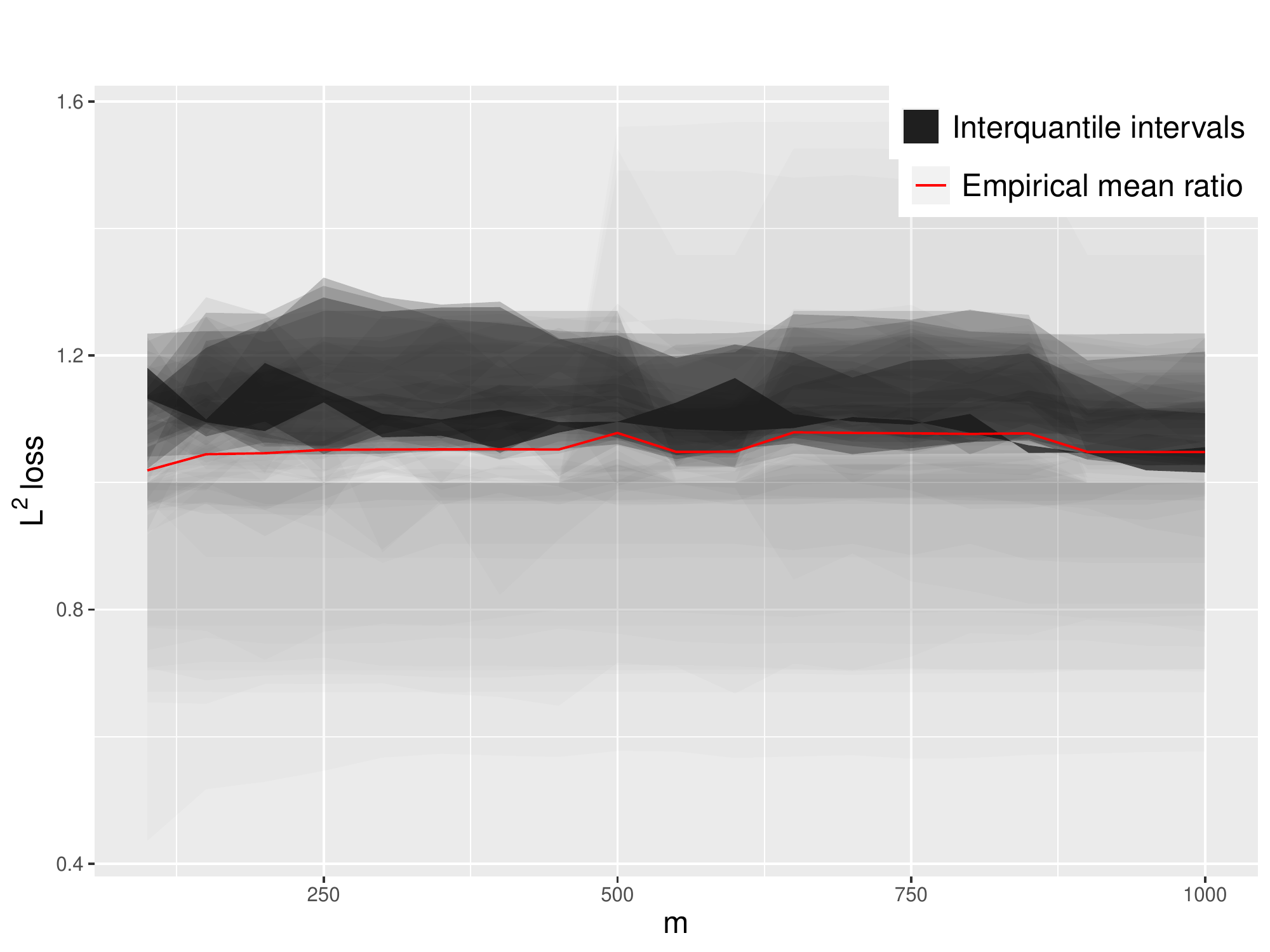}
  \caption{}
  \label{si:rate:ram}
\end{subfigure}%
\caption{Empirical convergence rate for weighted sum estimators with Bayesian and model selection weights, and their ratio}
\label{si:rate}
\end{figure}

\begin{te}
  In the remaining part of this section we illustrate the numerical performance of the weighted sum estimators and their dependence on the sample sizes $n$ and $m$ by reporting the Bayesian empirical risk obtained by repeating 100 times a procedure described next.
  In opposite to above we randomly pick a noise density $\edf$ and a density of interest $\xdf$, respectively, from a family of wrapped normal distributions and a family of wrapped asymetric Laplace distributions.
  For the noise density the mean and concentration parameters are uniformly distributed in $[0, 1]$.
  For the density of interest, location parameter is uniformly-distributed in $[0, 1]$, and both the scale, and the asymmetry parameters follow a $\Gamma$ distribution with shape $1$ and scale $5$.
  Note that the families differ from the ones used to calibrate the constant $\cpen$.
  Next we generate a sample $(\rE_i)_{i\in\nset{\ssE}}$ of size $\ssE = 1000$ from $\edf$, and a sample $(\rY_i)_{i\in\nset{\ssY}}$ of size $\ssY = 1000$ from $\ydf = \xdf \oast \edf$. For a range of subsamples
  with different samples sizes we construct the estimators $\hxdfAg$, and $\hxdfAg[{\msWe[]}]$ as in
  \cref{au:ag:ub:pnp} and compute their $\Lp[2]$-losses. Given the results of the $100$ repetitions, for the
  different values of $n$ and $m$ we use the sample of $\Lp[2]$-losses to compute estimators of the mean
  squared error and the quantiles of the distribution of the $\Lp[2]$-loss. The evolution of the $\Lp[2]$-loss for the weighted sum
  estimator with Bayesian weights or model selection weights, and their ratio, is represented in
  \cref{si:matrices}, when the sample sizes $\ssY$ and $\ssE$ vary.  In \cref{si:matrices:ag,si:matrices:ms}, both empirical errors decrease nicely as $\ssY$ and $\ssE$ increase.
  In \cref{si:matrices:ra}, it seems like, on smaller sample sizes the estimator with Bayesian weights performs
  better than the one with model selection weights, while the opposite happens for larger samples.
\end{te}
\begin{te}
  In \cref{si:rate} more attention is given to the spread of the $\Lp[2]$-loss around its empirical mean.
  The three columns (from left to right) refer to the estimator with Bayesian or  model selection weights, and their ratio.
  In the first row (\cref{si:rate:agn,si:rate:msn,si:rate:ran}) the noise sample size is fixed at $\ssE = 500$
  and in each graph  the sample-size $n$ increases from $100$ to $1000$.
  In the second row (\cref{si:rate:agmn,si:rate:msmn,si:rate:ramn}) both sample have the same size $\ssE =
  \ssY$ which again in each graph  increases from $100$ to $1000$.
  In the last row (\cref{si:rate:agm,si:rate:msm,si:rate:ram}) the size of the noisy sample is fixed at $\ssY =
  500$  and in each graph  the sample-size $m$ of the noise  increases from $100$ to $1000$.
  These graphics show that the distribution of the $\Lp[2]$-losses is skewed. However, in all cases, both estimators behave reasonably.
\end{te}

\begin{te}
  The simulations were performed with the R software, using the libraries 'circular', 'ggplot2', 'reshape2', 'foreach', and 'doParallel'.
  (see \cite{R-Core-Team:2018aa, Agostinelli:2017aa, Wickham:2016aa, Wickham:2007aa, foreach, doParallel}).
  All the scripts are available upon request to the authors.
\end{te}

\appendix
\setcounter{subsection}{0}
\section*{Appendix}
\numberwithin{equation}{subsection}
\renewcommand{\thesubsection}{\Alph{subsection}}
\renewcommand{\theco}{\Alph{subsection}.\arabic{co}}
\numberwithin{co}{subsection}
\renewcommand{\thelem}{\Alph{subsection}.\arabic{lem}}
\numberwithin{lem}{subsection}
\renewcommand{\therem}{\Alph{subsection}.\arabic{rem}}
\numberwithin{rem}{subsection}
\renewcommand{\thepr}{\Alph{subsection}.\arabic{pr}}
\numberwithin{pr}{subsection}

\subsection{Preliminaries}\label{a:prel}
\begin{te}
This section gathers technical results. The next result
is due to \cite{johannes2020adaptive}.
\end{te}

\begin{lem}\label{re:contr}Given $\ssY\in\Nz$ and
  $\pxdfPr[],\dxdfPr[]\in\Lp[2]$ consider the
  families of  orthogonal projections
  $\setB{\dxdfPr=\dProj{\Di}{}\dxdfPr[],\Di\in\nset{\ssY}}$ and $\setB{\pxdfPr=\dProj{\Di}{}\pxdf,\Di\in\nset{\ssY}}$. If $\VnormLp{\dProj{\Di}{}^\perp\pxdf}^2=\VnormLp{\ProjC[0]\pxdf}^2\bias^2(\pxdf)$ for all
  $\Di\in\nset{\ssY}$, then for any $l\in\nset{n}$ holds
 \begin{resListeN}[]
\item\label{re:contr:e1}
$\VnormLp{\dxdfPr[k]}^2-\VnormLp{\dxdfPr[l]}^2\leq
\tfrac{11}{2}\VnormLp{\dxdfPr[l]-\pxdfPr[l]}^2-\tfrac{1}{2}\VnormLp{\ProjC[0]\pxdf}^2\{\bias[k]^2(\pxdf)-\bias[l]^2(\pxdf)\}$,
for all $k\in\nsetro{l}$;
\item\label{re:contr:e2}
$\VnormLp{\dxdfPr[k]}^2-\VnormLp{\dxdfPr[l]}^2\leq \tfrac{7}{2}\VnormLp{\dxdfPr[k]-\pxdfPr[k]}^2+\tfrac{3}{2}\VnormLp{\ProjC[0]\pxdf}^2
\{\bias[l]^2(\pxdf)-\bias[k]^2(\pxdf)\}$, for all $k\in\nsetlo{l,n}$.
\end{resListeN}
\end{lem}

\begin{te}
 The next assertion
provides our key arguments in order to control the deviations of the
reminder terms.  Both inequalities are due to
\cite{Talagrand1996}, the formulation of the first part \cref{re:tal:e4} can be found for
example in \cite{KleinRio2005}, while the second part \cref{re:tal:e5} is  based on
equation (5.13) in Corollary 2 in \cite{BirgeMassart1995} and stated
in this form for example in  \cite{ComteMerlevede2002}.
\end{te}

\begin{lem}(Talagrand's inequalities)\label{re:tal} Let
  $(\rOb_i)_{i\in\nset[1,]{\ssY}}$ be independent $\cZ$-valued random variables
  and let
  $\overline{\nu_{\He}}=n^{-1}\sum_{i \in \nset{\ssY}}\left[\nu_{\He}(\rOb_i)-\Ex\left(\nu_{\He}(\rOb_i)\right)
  \right]$ for $\nu_{\He}$ belonging to a countable class
  $\{\nu_{\He},\He\in\cH\}$ of measurable functions. If the following
  conditions are satisfied
\begin{equation}\label{re:tal:e3}
	\sup_{\He\in\cH}\sup_{z\in\cZ}|\nu_{\He}(z)|\leq \Tch,\qquad \Ex(\sup_{\He\in\cH}|\overline{\nu_{\He}}|)\leq \TcH,\qquad \sup_{\He\in\cH}\frac{1}{n}\sum_{i\in\nset{\ssY}} \Var(\nu_{\He}(\rOb_i))\leq \Tcv,
\end{equation}
 then there is an universal numerical constant $\cst{}>0$ such that
\begin{align}
	 &\Ex\vectp{\sup_{\He\in\cH}|\overline{\nu_{\He}}|^2-6\TcH^2}\leq \cst{} \left[\frac{\Tcv}{n}\exp\left(\frac{-n\TcH^2}{6\Tcv}\right)+\frac{\Tch^2}{n^2}\exp\left(\frac{-n \TcH}{100\Tch}\right) \right]\label{re:tal:e4} \\
	&\pM\big(\sup_{\He\in\cH}|\overline{\nu_{\He}}|^2\geq6\TcH^2\big)\leq 3\big[\exp\big(\frac{-n\TcH^2}{400\Tcv}\big)+\exp\big(\frac{-n\TcH}{200\Tch}\big)\big].\label{re:tal:e5}
\end{align}
\end{lem}

\begin{rem}\label{rem:re:tal} Introduce  the unit ball $\mBaH:=\set{\He\in\mHiH:\VnormLp{\He}\leq1}$
 contained in the linear subspace
 $\mHiH$.
 Setting
 $\nu_{\He}(\rY)=\sum_{|j|\in\nset{\Di}}\ofHe{j}\fedfI[j]\bas_j(-\rY)$
 we have
\begin{equation*}
	\VnormLp{\txdfPr-\xdfPr}^2
=\sup_{\He\in\mBaH}|\sum_{|j|\in\nset{\Di}}\fedfI[j]\{\tfrac{1}{\ssY}\sum_{i \in \nset{\ssY}}(\bas_j(-\rY_i)-\fydf[j])\}\ofHe{j}|^2
=\sup_{\He\in\mBaH}|\overline{\nu_{\He}}|^2.
\end{equation*}
The last identity provides the necessary argument to link the
next  \cref{re:conc,re:cconc} and Talagrand's inequalities in \cref{re:tal}.  Note that, the unit ball $\mBaH$ is not a countable set of functions, however, it contains a countable dense subset, say $\cH$, since $\Lp[2]$ is separable, and it is straightforward to see that $\sup_{\He\in\mBaH}|\overline{\nu_{\He}}|^2=\sup_{\He\in\cH}|\overline{\nu_{\He}}|^2$.
\remEnd
\end{rem}

\begin{te}
 The proof of \cref{re:conc} given in
  \cite{JohannesSchwarz2013a}  makes use of  \cref{re:tal} by computing the quantities
  $\Tch$, $\TcH$, and $\Tcv$ which verify the three inequalities \eqref{re:tal:e3}.   We provide in \cref{re:cconc} a slight
  modification of this result following along the lines of
  the proof of  \cref{re:conc} in \cite{JohannesSchwarz2013a}.
\end{te}

\begin{lem}\label{re:conc}
  Let
  $\miSv= \max\{\iSv[j],j\in\nset{\Di}\}$, $\cmSv\geq1$ and
  $\DipenSv=\cmSv\Di \miSv$, then there is a numerical constant
  $\cst{}$ such that for all $\ssY\in\Nz$ and
  $\Di\in\nset{\ssY}$ holds
  \begin{resListeN}[]
  \item\label{re:conc:i}
    $\nEx \vectp{\VnormLp{\txdfPr-\xdfPr}^2 - 12\DipenSv\ssY^{-1}}  \leq
    \cst{} \bigg[\tfrac{\Vnormlp[1]{\fydf}\,\miSv}{\ssY}\exp\big(\tfrac{-\cmSv\Di}{3\Vnormlp[1]{\fydf}}\big)+\tfrac{2\Di\miSv}{n^2}\exp\big(\tfrac{-\sqrt{n\cmSv}}{200}\big) \bigg]$
  \item\label{re:conc:ii}
    $\nVg\big(\VnormLp{\txdfPr-\xdfPr}^2 \geq 12\DipenSv\ssY^{-1}\big)\leq
    3 \bigg[\exp\big(\tfrac{-\cmSv\Di}{200\Vnormlp[1]{\fydf}}\big)
    +\exp\big(\tfrac{-\sqrt{\ssY\cmSv}}{200}\big)\bigg]$
  \end{resListeN}
\end{lem}

\begin{lem}\label{re:cconc}
  Consider
  $\hxdfPr-\pxdfPr=\sum_{|j|\in\nset{\Di}}\hfedfmpI[j](\hfydf[j]-\fydf[j])\bas_j$.
  Denote by
  $\FuVg[\ssY]{\rY|\rE}$ and $\FuEx[\ssY]{\rY|\rE}$ the conditional
  distribution and expectation, respectively, of
  $(\rY_i)_{i\in\nset{\ssY}}$ given $(\rE_i)_{i\in\nset{\ssE}}$.  Let
  $\eiSv[j]=|\hfedfmpI[j]|^2$,
  $\oeiSv=\tfrac{1}{\Di}\sum_{j\in\nset{\Di}}\eiSv[j]$,
  $\meiSv= \max_{j\in\nset{\Di}}\eiSv[j]$,
  $\DiepenSv=\cmeiSv\Di \meiSv$ and $\cmeiSv\geq1$.  Then there is a
  numerical constant $\cst{}$ such that for all $\ssY\in\Nz$ and
  $\Di\in\nset{\ssY}$ holds
  \begin{resListeN}[]
  \item\label{re:cconc:i}
    $\FuEx[\ssY]{\rY|\rE}
    \vectp{\VnormLp{\hxdfPr-\pxdfPr}^2 - 12\DiepenSv\ssY^{-1}}  \leq
    \cst{} \bigg[\tfrac{\Vnormlp[1]{\fydf}\,\meiSv}{\ssY}
    \exp\big(\tfrac{-\cmeiSv\Di}{3\Vnormlp[1]{\fydf}}\big)
    +\tfrac{2\Di\meiSv}{n^2}\exp\big(\tfrac{-\sqrt{\ssY\cmeiSv}}{200}\big) \bigg]$
  \item\label{re:cconc:ii}
$\FuVg[\ssY]{\rY|\rE}\big(\VnormLp{\hxdfPr-\pxdfPr}^2 \geq 12\DiepenSv\ssY^{-1}\big)\leq
3 \bigg[\exp\big(\tfrac{-\cmeiSv\Di}{200\Vnormlp[1]{\fydf}}\big)
+\exp\big(\tfrac{-\sqrt{\ssY\cmeiSv}}{200}\big)\bigg]$
\end{resListeN}
\end{lem}

\begin{pro}[Proof of \cref{re:cconc}.]\label{pro:cconc}
For  $\He\in\mBaH$ set
$\nu_{\He}(\rY)=\sum_{|j|\in\nset{\Di}}\ofHe{j}\hfedfmpI[j]\bas_j(-\rY)$ where
 $\FuEx[\ssY]{\rY|\rE}\nu_{\He}(\rY)=\sum_{|j|\in\nset{\Di}}\ofHe{j}\hfedfmpI[j]\fydf[j]$ and
 $\VnormLp{\hxdfPr-\pxdfPr}^2=\sup_{\He\in\mBaH}|\overline{\nu_{\He}}|^2$ (see \cref{rem:re:tal}).
 We
 intent to apply \cref{re:tal}. Therefore,  we compute  next quantities $\Tch$, $\TcH$,
and $\Tcv$  verifying the three  inequalities required in
\cref{re:tal}. First, we have $\sup_{\He\in\mBaH}\sup_{y\in[0,1]}|\nu_{\He}(y)|^2 = 2\sum_{j\in\nset{\Di}} \eiSv[j]
  \leq 2\Di\meiSv=:\Tch^2$.
Next, find  $\TcH$.
Exploiting $\sup_{\He\in\mBaH}
|\VskalarLp{\hxdfPr-\pxdfPr,\He}|^2=
\sum_{|j|\in\nset{\Di}}\eiSv[|j|]\,|\hfydf[j]-\fydf[j]|^2$ and
$\FuEx[\ssY]{\rY|\rE}|\hfydf[j]-\fydf[j]|^2\leq
\tfrac{1}{\ssY}$, it holds
$\FuEx[\ssY]{\rY|\rE}\big(\sup_{\He\in\mBaH} |\VskalarLp{\hxdfPr-\pxdfPr,\He}^2\big)
\leq 2\sum_{|j|\in\nset{\Di}}\eiSv[|j|]/\ssY\leq
2\DiepenSv/\ssY  =: \TcH^2$.
\noindent Finally, consider $\Tcv$.  Using
$\FuEx[]{\rY|\rE}\big(\bas_j(Y_1)\bas_{j'}(-Y_1)\big)
=\fydf[j'-j]$ for each $\He\in \mBaH$  holds
\begin{equation*}
  \FuEx[]{\rY|\rE}|\nu_{\He}(Y_1)|^2
=\sum_{|j|,|j'|\in\nset{\Di}}\fHe{j}\ohfedfmpI[j]\fydf[j'-j] \hfedfmpI[j']\,\ofHe{j'}
=\Vskalarlp{\DiPro[k]\hA\DiPro[k]\;\fHe{},\fHe{}}
\end{equation*}
defining the Hermitian and positive semi-definite matrix $\hA:=
\Zsuite[j,j']{\ohfedfmpI[j] \hfedfmpI[j'] \fydf[j'-j]}$ and the mapping
$\DiPro[k]:\Cz^\Zz\to\Cz^\Zz$ with
$z\mapsto\DiPro[k]z=(z_l\Ind{\{|l|\in\nset{\Di}\}})_{l\in\Zz}$. Obviously,
$\DiPro[k]$ is an orthogonal projection from $\lp^2$ onto the linear
subspace spanned by all $\lp^2$-sequences with support on the
index-set $\nset{-\Di,-1}\cup\nset{\Di}$. Straightforward
   algebra shows
\begin{multline*}
  \sup_{\He\in\mBaH} \tfrac{1}{\ssY}\sum_{i \in \nset{\ssY}}\Var_{\rY|\rE}(\nu_{\He}(Y_i))
\leq \sup_{\He\in\mBaH} \Vskalarlp{\DiPro[k]\hA\DiPro[k]\fHe{},\fHe{}}= \sup_{\He\in\mBaH} \Vnormlp{\DiPro[k]\hA\DiPro[k]\fHe{}}\leq\Vnorm[s]{\DiPro[k]\hA\DiPro[k]}.
\end{multline*}
where $\Vnorm[s]{M}:=\sup_{\Vnormlp{x}\leq 1}\Vnormlp{Mx}$ denotes the
spectral-norm of a linear $M:\lp^2\to\lp^2$. For  a
sequence $z\in\Cz^\Zz$ let $\Diag{z}$   be the multiplication operator
given by $\Diag{z}x:=\Zsuite{z_j x_j}$.   Clearly, we have
 $\DiPro[k]\hA\DiPro[k] = \DiPro[k]\Diag{\hfedfmpI[]} \DiPro[k]  \cC_{\fydf}
 \DiPro[k] \Diag  {\ohfedfmpI[]}\DiPro[k],$
where $  \cC_{\fydf} := \Zsuite[j,j']{\fou[j-j']{\ydf}}.$
Consequently,
\begin{multline*}
   \sup_{\He\in\mBaH}\tfrac{1}{\ssY}\sum_{i \in \nset{\ssY}}\Var_{\rY|\rE}(\nu_{\He}(Y_i))
\leq \Vnorm[s]{\DiPro[k]\Diag{\hfedfmpI[]}\DiPro[k]}\; \Vnorm[s]{ \cC_{\fydf}}\;\Vnorm[s]{\DiPro[k]\Diag{\ohfedfmpI[]}\DiPro[k]}\\=
\Vnorm[s]{\DiPro[k]\Diag{\hfedfmpI[]} \DiPro[k]}^2\; \Vnorm[s]{\cC_{\fydf}},
\end{multline*}
where $\Vnorm[s]{\DiPro[k]\Diag{\hfedfmpI[]}\DiPro[k]}^2=
\max\{\eiSv[j],j\in\nset{\Di}\}=\meiSv$.
For $(\cC_{\fydf}z)_k:=
\sum_{j\in\Zz}\fydf[j-k]z_j$, $k\in\Zz$ it is
easily verified that $\Vnormlp{\cC_{\fydf} z}^2\leq
\Vnormlp[1]{\fydf}^2\Vnormlp{z}^2$  and hence $ \Vnorm[s]{\cC_{\fydf}}\leq\Vnormlp[1]{\fydf}$, which
  implies
  \begin{equation*}
    \sup_{\He\in\mBaH}\tfrac{1}{\ssY}\sum_{i \in \nset{\ssY}}\Var_{\rY|\rE}(\nu_{\He}(Y_i))\leq\Vnormlp[1]{\fydf}\,\meiSv=:\Tcv.
  \end{equation*}
  Replacing in \cref{rem:re:tal} \eqref{re:tal:e3} and \eqref{re:tal:e4}
  the quantities $\Tch,\TcH$ and $\Tcv$ together with
  $\DipeneSv=\cmeiSv\Di\meiSv$ gives the assertion \ref{re:cconc:i} and
  \ref{re:cconc:ii}
  , which completes the
  proof.\proEnd
\end{pro}

\begin{lem}\label{oSv:re}
There is a finite numerical constant $\cst{}>0$ such that
for all $j\in\Zz$ hold
\begin{inparaenum}[\dgrau\upshape(i)]
\item[] $\ssE^2\FuEx[\ssE]{\rE}\Vabs{\fedf[j]-\hfedf[j]}^4\leq\cst{}$,
\item[\mylabel{oSv:re:i}{{\dr\upshape(i)}}]
$\FuEx[\ssE]{\rE}\big(\Vabs{\fedf[j]\hfedfmpI[j]}^2\big)\leq 4$;
\item[\mylabel{oSv:re:ii}{{\dr\upshape(ii)}}]
$\FuVg[\ssE]{\rE}(\Vabs{\hfedfmpI[j]}^2<1/\ssE)\leq4(1\wedge \iSv[j]/\ssE)$,
\item[\mylabel{oSv:re:iii}{{\dr\upshape(iii)}}] $\FuEx[\ssE]{\rE}\big(\Vabs{\fedf[j]-\hfedf[j]}^2\Vabs{\hfedfmpI[j]}^2\big)\leq
  4\cst{}(1\wedge \iSv[j]/\ssE)$. Given $\Di\in\Nz$ for all
  $j\in\nset{\Di}$ we have
\item[\mylabel{re:xevent:e1}{{\dr\upshape(iv)}}] $\FuVg[\ssE]{\rE}\big(|\hfedf[j]/\fedf[j]-1|>1/3\big)\leq 2\exp\big(-\tfrac{\ssE|\fedf[j]|^2}{72}\big)\leq 2\exp\big(-\tfrac{\ssE}{72\miSv}\big)$.
\end{inparaenum}
\end{lem}

\begin{pro}[Proof of \cref{oSv:re}.]
The elementary properties \ref{oSv:re:i}-\ref{oSv:re:iii} are shown,
for example, in
\cite{JohannesSchwarz2013a} and the assertion \ref{re:xevent:e1} follows directly from Hoeffding's inequality.\proEnd
\end{pro}

\begin{lem}\label{re:evrest}Let $\ssE,\Di\in\Nz$ and set
  $\aixEv[\Di]:=\{1/4\leq\eiSv[j]/\iSv[j]\leq9/4:\;\forall\;j\in\nset{\Di}\}$.
\begin{resListeN}[]
\item\label{re:evrest:i}
If $\miSv[k]\leq(4/9)\ssE$ then $\FuVg[\ssE]{\rE}(\aixEv^c)\leq 2\Di\exp\big(-\tfrac{\ssE}{72\miSv}\big)$.
\item\label{re:evrest:ii}
For  $\ssE_{\Di}:=\floor{9\miSv/4}$  holds
$\FuVg[\ssE]{\rE}(\aixEv^c)\leq(555\Di\ssE_{k}^2\ssE^{-2})\wedge(12\Di\ssE_{\Di}\ssE^{-1})$
for all
$\ssE\in\Nz$.
\item\label{re:evrest:iii}
If $\ssE\geq289\log(\Di+2)\cmiSv\miSv$ then $\FuVg[\ssE]{\rE}(\aixEv^c)\leq(11226\ssE^{-2})\wedge(53\ssE^{-1})$.
\end{resListeN}
\end{lem}

\begin{pro}[Proof of \cref{re:evrest}.]
We start our proof with the observation that  for each $j\in\Zz$ with
$\iSv[j]\leq (4/9)\ssE$ holds
$\{|\hfedf[j]/\fedf[j]-1|\leq1/3\}\subseteq
\{1/2\leq |\fedf[j]\hfedfmpI[j]|\leq3/2\}$
Consequently, if $\miSv\leq
(4/9)\ssE$ then
  $\aixEv^c
  \subset\bigcup_{j \in \nset{\Di}}\{|\hfedf[j]/\fedf[j]-1|>1/3\}$
and hence \ref{re:evrest:i} follows from \cref{oSv:re}
\ref{re:xevent:e1}.
Consider \ref{re:evrest:ii}. Given  $\Di\in\Nz$ and
$\ssE_{\Di}:=\floor{9\miSv/4}\in\Nz$ we distinguish for $\ssE\in\Nz$ the cases \begin{inparaenum}[i]\renewcommand{\theenumi}{\dgrau\rm(\alph{enumi})}\item\label{pro:evrest:ii:c1}
$\ssE>\ssE_{\Di}$
and \item\label{pro:evrest:ii:c2}$\ssE\in\nset{\ssE_{\Di}}$.
\end{inparaenum}
In case \ref{pro:evrest:ii:c1} it holds $ \miSv\leq (4/9)\ssE$, and hence \ref{re:evrest:i} implies \ref{re:evrest:ii}.
In case \ref{pro:evrest:ii:c2}  \ref{re:evrest:ii} holds trivially, since $\FuVg[\ssE]{\rE}(\aixEv^c)\leq \ssE_{\Di}^2\ssE^{-2}\wedge\ssE_{\Di}\ssE^{-1}$.
Consider \ref{re:evrest:iii}.
Since $\ssE\geq289\log(\Di+2)\cmiSv\miSv\geq
(9/4)\miSv[\Di]$ from \ref{re:evrest:i} follows
\begin{equation*}
\ssE^2\FuVg[\ssE]{\rE}(\aixEv^c)\leq \Di\ssE^2\exp\big(-\tfrac{\ssE}{72\miSv}\big)
\leq 11226 \Di\miSv^2\exp\big(-\tfrac{\ssE}{144\miSv}\big) \leq 11226
\end{equation*}
and analogously
$\ssE\FuVg[\ssE]{\rE}(\aixEv^c) \leq 53$, which completes the proof.
\proEnd\end{pro}

\begin{lem}\label{re:aixEv} Consider for any
$l\in\Nz$  the event
$\aixEv[l]:=\setB{\tfrac{1}{4}\leq\iSv[j]^{-1}\eiSv[j]\leq\tfrac{9}{4},\,\forall\,j\in\nset{l}}$.
For each $l\in\nset{\ssY}$ and  $k\in\nsetro{l}$ setting
$\VnormLp{\dProj{\Di}{l}\pxdfPr[\ssY]}^2:=\sum_{|j|\in\nsetlo{k,l}}\iSv[j]^{-1}\eiSv[j]|\fxdf[j]|^2$ hold
  \begin{resListeN}[]
  \item\label{re:aixEv:i}
    $\VnormLp{\dProj{\Di}{l}\pxdfPr[\ssY]}^2\leq\VnormLp{\dProj{\Di}{\ssY}\pxdfPr[\ssY]}^2=\VnormLp{\ProjC[\Di]\pxdfPr[\ssY]}^2$
    and
    $\VnormLp{\dProj{\Di}{l}\pxdfPr[\ssY]}^2\Ind{\aixEv[l]}\geq
    \tfrac{1}{4}\VnormLp{\ProjC[0]\xdf}^2(\bias^2(\xdf)-\bias[l]^2(\xdf))$.
  \end{resListeN}
  Moreover, for any $l\in\Nz$ and $\Di\in\nset{l}$ hold
  \begin{resListeN}[\addtocounter{ListeN}{1}]
  \item\label{re:aixEv:ii}
    $\meiSv[l]\leq\ssE$,
    $\tfrac{1}{4}\miSv[l]\leq
    \meiSv[l]\Ind{\aixEv[l]}\leq\tfrac{9}{4}\miSv[l]$,
    $\cmeiSv[l]\geq1$,
    $\tfrac{9}{100}\cmiSv[l]\leq\cmeiSv[l]\Ind{\aixEv[l]}\leq\tfrac{49}{16}\cmiSv[l]$,
    and hence
  $\tfrac{1}{50}\penSv\leq\peneSv\Ind{\aixEv[l]}\leq7\penSv$;
  \item\label{re:aixEv:iv} $\setB{\meiSv[l]<1}=\setB{\meiSv[l]=0}$, and hence
    $\peneSv[l]=\peneSv[l]\Ind{\{\meiSv[l]\geq1\}}$.
  \end{resListeN}
\end{lem}
\begin{pro}[Proof of \cref{re:aixEv}.]The assertions
  \ref{re:aixEv:i} and \ref{re:aixEv:ii} follow by elementary calculations from the definition of  the event
$\aixEv[l]$, and we omit the details.
Consider \ref{re:aixEv:iv}. For each  $j\in\Zz$ holds
$\eiSv[j]=|\hfedfmpI[j]|^2=0$ on the event
$\setB{|\hfedf[j]|^2<1/\ssE}$ and
$\eiSv[j]\geq1$ on the complement $\setB{|\hfedf[j]|^2\geq1/\ssE}$,
since  $|\hfedf[j]|^2\leq1$. Consequently,
$\setB{\eiSv[j]<1}=\setB{|\hfedf[j]|^2<1/\ssE}=\setB{\eiSv[j]=0}$,
which implies \ref{re:aixEv:iv},
and completes the proof.\proEnd
\end{pro}

\subsection{Proofs of \cref{ak}}\label{a:ak}
\begin{pro}[Proof of \cref{co:agg}.]
We start the proof with the observation that $\ftxdfPr[{\We[]}]{0}-\fxdf[0]=0$,
and for each $j\in\Zz$ holds
$\oftxdfPr[{\We[]}]{j}-\ofxdf[j]=\ftxdfPr[{\We[]}]{-j}-\fxdf[-j]$,  where
$\ftxdfPr[{\We[]}]{j}-\fxdf[j]=-\fxdf[j]$ for all $|j|>\ssY$ and
\begin{equation*}
  \ftxdfPr[{\We[]}]{j}-\fxdf[j]=\fedfI[j](\hfydf[j]-\fydf[j])\FuVg{\We[]}(\nset{|j|,\ssY})-\fxdf[j]\FuVg{\We[]}(\nsetro{|j|})\text{ for all }|j|\in\nset{\ssY}.
\end{equation*}
Consequently, (keep in mind that $|\fedfI[j]|^2=\iSv[j]$)  we  have
  \begin{multline}\label{co:agg:pro1}
    \VnormLp{\txdfPr[{\We[]}]-\xdf}^2
\leq
   \sum_{|j|\in\nset{\ssY}}2\{\iSv[j]|\hfydf[j]-\fydf[j]|^2 \FuVg{\We[]}(\nset{|j|,\ssY})\} \\+ \sum_{|j|\in\nset{\ssY}}2|\fxdf[j]|^2\FuVg{\We[]}(\nsetro{|j|})+\sum_{|j|>n}|\fxdf[j]|^2,
 \end{multline}
where we consider the first  and the two other terms on the right hand
side separately. Considering the first term we split the sum into two parts.  Precisely,
\begin{multline}\label{co:agg:pro2}
\sum_{|j|\in\nset{\ssY}}\iSv[j]|\hfydf[j]-\fydf[j]|^2
\FuVg{\We[]}(\nset{|j|,\ssY})
\leq \VnormLp{\txdfPr[\pDi]-\xdfPr[\pDi]}^2
+\sum_{l\in\nsetlo{\pDi,\ssY}}\We[l]\VnormLp{\txdfPr[l]-\xdfPr[l]}^2\\
\hfill\leq\tfrac{1}{7}\pen[\pDi]
+\sum_{l\in\nset{\pDi,\ssY}}\vectp{\VnormLp{\txdfPr[l]-\xdfPr[l]}^2-\pen[l]/7}\\
+\tfrac{1}{7}\sum_{l\in\nsetlo{\pDi,\ssY}}\We[l]\pen[l]\Ind{\{\VnormLp{\txdfPr[l]-\xdfPr[l]}^2\geq\pen[l]/7\}}
+\tfrac{1}{7}\sum_{l\in\nsetlo{\pDi,\ssY}}\pen[l]\We[l]\Ind{\{\VnormLp{\txdfPr[l]-\xdfPr[l]}^2<\pen[l]/7\}}
\end{multline}
Considering the second and third term we split the first sum into two
parts and obtain
\begin{multline}\label{co:agg:pro3}
\sum_{|j|\in\nset{\ssY}}|\fxdf[j]|^2\FuVg{\We[]}(\nsetro{|j|})+\sum_{|j|>\ssY}|\fxdf[j]|^2\\
\hspace*{5ex}\leq  \sum_{|j|\in\nset{\mDi}}|\fxdf[j]|^2\FuVg{\We[]}(\nsetro{|j|})+ \sum_{|j|\in\nsetlo{\mDi,n}}|\fxdf[j]^2+
  \sum_{|j|>n}|\fxdf[j]|^2\\\hfill
\leq \VnormLp{\Proj[{\mHiH[0]^\perp}]\xdf}^2\{\FuVg{\We[]}(\nsetro{\mDi})+\sbF[\mDi]\}
\end{multline}
Combining  \eqref{co:agg:pro1} and  \eqref{co:agg:pro2}, \eqref{co:agg:pro3} we obtain   the assertion, which completes the proof.\proEnd
\end{pro}
\subsubsection{Proof of \cref{ak:ag:ub:pnp} and \cref{ak:ag:ub2:pnp}}\label{a:ak:rb}
\begin{pro}[Proof of \cref{ak:ag:ub:pnp}.]

\begin{te}We present the main arguments to prove \cref{ak:ag:ub:pnp}. The technical details are gathered in \cref{ak:re:SrWe:ag,ak:re:SrWe:ms,ak:re:nd:rest,ak:ag:ub:p}
  in the end of this section. Keeping in mind the definition
  \eqref{oo:de:doRao} and \eqref{ak:de:penSv} here and subsequently we use that
    \begin{equation}\label{ak:ass:pen:oo:c}
    \doRaL{\Di}\geq\bias^2(\xdf)\quad\text{and}\quad\cpen\doRaL{\Di}\geq\penSv
    \quad\text{for all }\Di\in\nset{\ssY}.
  \end{equation}
\end{te}

\begin{te}
   For arbitrary
  $\pdDi,\mdDi\in\nset{\ssY}$ (to be choosen suitable below) let us define
  \begin{multline}\label{ak:de:*Di:ag}
    \mDi:=\min\set{\Di\in\nset{\mdDi}: \VnormLp{\ProjC[0]\xdf}^2\bias^2(\xdf)\leq
      \VnormLp{\ProjC[0]\xdf}^2\bias[\mdDi]^2(\xdf)+4\penSv[\mdDi]}\quad\text{and}\\
    \pDi:=\max\set{\Di\in\nset{\pdDi,\ssY}:
      \penSv \leq 6\VnormLp{\ProjC[0]\xdf}^2\bias[\pdDi]^2(\xdf)+ 4\penSv[\pdDi]}
  \end{multline}
  where the defining set obviously contains $\mdDi$ and $\pdDi$,
  respectively, and hence, it is not empty.
\end{te}

\begin{te}
  We intend to combine the upper bound in \eqref{co:agg:e1} and the
  bounds considering  \pcw $\We[]=\rWe[]$
  as in \eqref{ak:de:rWe} and \msw $\We[]=\msWe[]$
  as in \eqref{ak:de:msWe} given
  in \cref{ak:re:SrWe:ag} and  \cref{ak:re:SrWe:ms}, respectively. First note, that due to \cref{ak:re:SrWe:ag} \ref{ak:re:SrWe:ag:i} we have
  \begin{multline*}
    \nEx\FuVg{\rWe[]}(\nsetro{\mDi})\leq\Ind{\{\mDi>1\}}
    \tfrac{1}{\rWc\cpen}\exp\big(-\tfrac{3\rWc\cpen}{14}\DipenSv[\mdDi]\big)\\+\Ind{\{\mDi>1\}}
    \nVg\big(\VnormLp{\txdfPr[\mdDi]-\xdfPr[\mdDi]}^2
    \geq\penSv[\mdDi]/7\big)
  \end{multline*}
  and, hence from \eqref{co:agg:e1} for  \pcw $\We[]=\rWe[]$
  as in \eqref{ak:de:rWe} follows immediately
  \begin{multline}\label{co:agg:ag}
    \noRi{\txdfAg}{\xdf}{\iSv}
    \leq \tfrac{2}{7}\penSv[\pDi]
    +2\VnormLp{\ProjC[0]\xdf}^2\bias[\mDi]^2(\xdf)
    \\\hfill
    +\ssY^{-1}\tfrac{32}{7\rWc}
    +\tfrac{2}{\rWc\cpen}\VnormLp{\ProjC[0]\xdf}^2\Ind{\{\mDi>1\}}
    \exp\big(-\tfrac{3\rWc\cpen}{14}\DipenSv[\mdDi]\big)\\
    +2\sum_{\Di\in\nset{\pdDi,n}}\nEx\vectp{\VnormLp{\txdfPr-\xdfPr}^2-\tfrac{1}{7}\penSv}
    +\tfrac{2}{7}\sum_{\Di\in\nset{\pdDi,\ssY}}\penSv
    \nVg\big(\VnormLp{\txdfPr-\xdfPr}^2\geq\tfrac{1}{7}\penSv\big)\\
    \hfill+2\VnormLp{\ProjC[0]\xdf}^2\Ind{\{\mDi>1\}}\nVg
    \big(\VnormLp{\txdfPr[\mdDi]-\xdfPr[\mdDi]}^2\geq\tfrac{1}{7}
    \penSv[\mdDi]\big)
  \end{multline}
\end{te}

\begin{te}
On the other hand for \msw $\We[]=\msWe[]$ we combine again the upper bound in \eqref{co:agg:e1} and the bounds given in \cref{ak:re:SrWe:ms}.
Clearly, due to  \cref{ak:re:SrWe:ms} we have
$\nEx\FuVg{\msWe[]}(\nsetro{\mDi})=\nVg\big(\VnormLp{\txdfPr[\mdDi]-\xdfPr[\mdDi]}^2\geq\penSv[\mdDi]/7\big)$
and, hence from \eqref{co:agg:e1} follows immediately
\begin{multline}\label{co:agg:ms}
  \noRi{\txdfPr[\hDi]}{\xdf}{\iSv}\leq \tfrac{2}{7}\penSv[\pDi]
+2\VnormLp{\ProjC[0]\xdf}^2\bias[\mDi]^2(\xdf)
\\
+2\sum_{\Di\in\nset{\pdDi,\ssY}}\nEx\vectp{\VnormLp{\txdfPr-\xdfPr}^2-\tfrac{1}{7}\penSv}
+\tfrac{2}{7}\sum_{\Di\in\nset{\pdDi,\ssY}}\penSv\nVg\big(\VnormLp{\txdfPr-\xdfPr}^2\geq\tfrac{1}{7}\penSv\big)\\\hfill
+2\VnormLp{\ProjC[0]\xdf}^2\Ind{\{\mDi>1\}}\nVg\big(\VnormLp{\txdfPr[\mdDi]-\xdfPr[\mdDi]}^2\geq\tfrac{1}{7}\penSv[\mdDi]\big)
\end{multline}
The deviations of the last three terms in \eqref{co:agg:ag} and \eqref{co:agg:ms}
 we bound in \cref{ak:re:nd:rest} by exploiting usual
concentration inequalities. Precisly,  we obtain
  \begin{multline}\label{ak:ag:ub:p3}
    \nEx\VnormLp{\txdfAg[{\We[]}]-\xdf}^2\leq \tfrac{2}{7}\penSv[\pDi]
    +2\VnormLp{\ProjC[0]\xdf}^2\sbF[\mDi]
 + \cst{}\VnormLp{\ProjC[0]\xdf}^2\Ind{\{\mDi>1\}}
 \exp\big(\tfrac{-\cmiSv[\mdDi]\mdDi}{\Di_{\ydf}}\big)
 \\
 +\cst{}\big(
 \VnormLp{\ProjC[0]\xdf}^2\Ind{\{\mDi>1\}}
+\miSv[\Di_{\ydf}]^2\Di_{\ydf}^3+\miSv[\ssY_{o}]^2 \big)\ssY^{-1}.
  \end{multline}
  Indeed, combining
  \cref{ak:re:nd:rest} and  \eqref{co:agg:ag} for  \pcw   we obtain
  \begin{multline}\label{ak:ag:ub:p1}
    \noRi{\txdfAg}{\xdf}{\iSv}\leq \tfrac{2}{7}\penSv[\pDi]
    +2\VnormLp{\ProjC[0]\xdf}^2\sbF[\mDi]
    \\\hfill
 + \cst{}\VnormLp{\ProjC[0]\xdf}^2\Ind{\{\mDi>1\}}\big( \tfrac{1}{\rWc}
     \exp\big(-\tfrac{3\rWc\cpen}{14}\DipenSv[\mdDi]\big)+
 \exp\big(\tfrac{-\cmiSv[\mdDi]\mdDi}{200\Vnormlp[1]{\fydf}}\big)\big)
 \\
 +\cst{}\big(\tfrac{1}{\rWc}+
 \VnormLp{\ProjC[0]\xdf}^2\Ind{\{\mDi>1\}}
+\miSv[\Di_{\ydf}]^2\Di_{\ydf}^3+\miSv[\ssY_{o}]^2 \big)\ssY^{-1}
  \end{multline}
 Therewith, be using that $\LiSv[\mdDi]\geq\liSv[\mdDi]$,
 $\tfrac{3\rWc\cpen}{14}>\tfrac{1}{200\Vnormlp[1]{\fydf}}>\tfrac{1}{\Di_{\ydf}}$
(since $\rWc\geq1$ and $\Vnormlp[1]{\fydf}\geq|\fydf[0]|=1$) from
\eqref{ak:ag:ub:p1}  follows the upper bound \eqref{ak:ag:ub:p3}.
Consider secondly  model selection weights $\We[]=\msWe[]$
  as in \eqref{ak:de:msWe}. Combining
  \cref{ak:re:nd:rest}, $200\Vnormlp[1]{\fydf}\leq \Di_{\ydf}$ and the upper bound given in \eqref{co:agg:ms}
  we obtain \eqref{ak:ag:ub:p3}.\
\end{te}

\begin{te} From the upper bound \eqref{ak:ag:ub:p3} for a suitable
coice of the dimension parameters $\mdDi,\pdDi\in\nset{n}$ we derive
  separately the risk bound in the two cases
  \ref{oo:xdf:p} and \ref{oo:xdf:np} considered in \cref{ak:ag:ub:pnp}.
The tedious
case-by-case analysis   for  \ref{ak:ag:ub:pnp:p} is deferred to \cref{ak:ag:ub:p} in the end of this section.
In case \ref{ak:ag:ub:pnp:np} with
$\oDi{\ssY}:=\noDiL\in\nset{n}$
and $\doRaL{\Di}$ as in \eqref{oo:de:doRao} we  set
$\pdDi:=\oDi{\ssY}$ and let $\mdDi\in\nset{\ssY}$. Keeping
\eqref{ak:ass:pen:oo:c} in mind
 the definition
\eqref{ak:de:*Di:ag} of $\pDi$ and $\mDi$ implies
$\penSv[\pDi] \leq 2(3\VnormLp{\ProjC[0]\xdf}^2+ 2\cpen)\doRaL{\oDi{\ssY}}$ and
$\VnormLp{\ProjC[0]\xdf}^2\sbF[\mDi]\leq
      (\VnormLp{\ProjC[0]\xdf}^2+4\cpen)\doRaL{\mdDi}$
which together with
$\doRaL{\mdDi}\geq\doRaL{\oDi{\ssY}}=\noRaL=\min\Nset[\Di\in\Nz]{\doRaL{\Di}}\geq\ssY^{-1}$
and exploiting
\eqref{ak:ag:ub:p3} implies the assertion
\eqref{ak:ag:ub:pnp:e2}, that is for all $\mdDi\in\nset{\ssY}$ holds%
 \begin{multline}\label{ak:ag:ub:pnp:p8}
 \noRi{\txdfAg[{\We[]}]}{\xdf}{\iSv}
   \leq
   \cst{}(\VnormLp{\ProjC[0]\xdf}^2\vee1)\big[\doRaL{\mdDi}\vee\exp\big(\tfrac{-\cmiSv[\mdDi]\mdDi}{\Di_{\ydf}}\big)\big]\\
   +\cst{}\big[\miSv[\Di_{\ydf}]^2\Di_{\ydf}^3+\miSv[\ssY_{o}]^2 \big]\ssY^{-1},
\end{multline}
with $\ssY_{o}=15(600)^4$,  which completes the proof of \cref{ak:ag:ub:pnp}.
\end{te}
\proEnd\end{pro}

\begin{pro}[Proof of \cref{ak:ag:ub2:pnp}.]
  Consider the case \ref{ak:ag:ub2:pnp:p}. Under \ref{ak:ag:ub2:pnp:pc}  for all
  $\ssY> \ssY_{\xdf,\iSv}$ we have trivially
  $\exp\big({-\cmiSv[\sDi{\ssY}]\sDi{\ssY}}/{\Di_{\ydf}}\big)\leq\ssY^{-1}$,
  while for $\ssY\in\nset{\ssY_{\xdf,\iSv}}$ holds
  $\exp\big({-\cmiSv[\sDi{\ssY}]\sDi{\ssY}}/{\Di_{\ydf}}\big)\leq1\leq
  \ssY_{\xdf,\iSv} \ssY^{-1}$. Thereby, from \eqref{ak:ag:ub:pnp:e1}
  in \cref{ak:ag:ub:pnp} follows immediately the assertion
  \ref{oo:xdf:p}.
  In case
  \ref{ak:ag:ub2:pnp:np} due to
  \ref{ak:ag:ub2:pnp:npc}  for $\oDi{\ssY}:=\noDiL$
  as in \eqref{oo:de:doRao} we have trivially
  $\exp\big({-\cmiSv[\oDi{\ssY}]\oDi{\ssY}}/{\Di_{\ydf}}\big)\leq
  \noRaL$ while for $\ssY\in\nset{\ssY_{\xdf,\iSv}}$ holds
  $\exp\big({-\cmiSv[\oDi{\ssY}]\oDi{\ssY}}/{\Di_{\ydf}}\big)\leq1\leq
  \ssY\noRaL\leq \ssY_{\xdf,\iSv} \noRaL$. Thereby, from
  \eqref{ak:ag:ub:pnp:e2} in \cref{ak:ag:ub:pnp} with
  $\noRaL=\min_{\Di\in\nset{\ssY}}\doRaL{\Di}$ follows
  \ref{ak:ag:ub2:pnp:np}, which completes the proof of \cref{ak:ag:ub2:pnp}.\proEnd
\end{pro}

\begin{te}
 Below  we state and prove the technical
 \cref{ak:re:SrWe:ag,ak:re:SrWe:ms,ak:re:nd:rest,ak:ag:ub:p} used in the proof of \cref{ak:ag:ub:pnp}. The
  proof of \cref{ak:re:SrWe:ag} is based on \cref{re:rWe} given first.
\end{te}

\begin{lem}\label{re:rWe} Consider  \pcw
  $\rWe[]$ as in \eqref{ak:de:rWe}. Let $l\in\nset{\ssY}$.
  \begin{resListeN}[]
  \item\label{re:rWe:i} For all $k\in\nsetro{l}$ holds
    $\rWe\Ind{\setB{7\VnormLp{\txdfPr[l]-\xdfPr[l]}^2<\penSv[l]}}$
    \\\null\hfill
    $\leq\exp\big(\rWn\big\{\tfrac{25}{14}\penSv[l]+\tfrac{1}{2}\VnormLp{\ProjC[0]\xdf}^2\sbFxdf[l]-\tfrac{1}{2}\VnormLp{\ProjC[0]\xdf}^2\sbFxdf -\penSv\big\}\big).$
  \item\label{re:rWe:ii} For all $\Di\in\nsetlo{l,\ssY}$ holds
    $\rWe\Ind{\setB{7\VnormLp{\txdfPr-\xdfPr}^2<\penSv}}$\\\null\hfill
    $\leq\exp\big(\rWn\big\{-\tfrac{1}{2}\penSv
    +\tfrac{3}{2}\VnormLp{\ProjC[0]\xdf}^2\sbFxdf[l]+\penSv[l]\big\}\big)$.
  \end{resListeN}
\end{lem}

\begin{pro}[Proof of \cref{re:rWe}.]
  Given $\Di,l\in\nset{\ssY}$ and an event $\dmEv{\Di}{l}$ (to be
  specified below) it follows
  \begin{multline}\label{re:rWe:pro1}
    \rWe\Ind{\dmEv{\Di}{l}}
    =\frac{\exp(-\rWn\{-\VnormLp{\txdfPr}^2+\penSv\})}
    {\sum_{l\in\nset{\ssY}}\exp(-\rWn\{-\VnormLp{\txdfPr[l]}^2+\penSv[l]\})}
    \Ind{\dmEv{\Di}{l}}\\
    \leq
    \exp\big(\rWn\big\{\VnormLp{\txdfPr}^2-\VnormLp{\txdfPr[l]}^2
    +(\penSv[l]-\penSv)\big\}\big)\Ind{\dmEv{\Di}{l}}
  \end{multline}
  We distinguish the two cases \ref{re:rWe:i} $\Di\in\nsetro{l}$ and
  \ref{re:rWe:ii} $\Di\in\nsetlo{l,n}$.  Consider first \ref{re:rWe:i}
  $\Di\in\nsetro{l}$. Due to  \cref{re:contr} \ref{re:contr:e1}
  (with  $\dxdf:=\txdfPr[\ssY]$ and $\pxdf:=\xdf$) from \eqref{re:rWe:pro1} we conclude
  \begin{multline*}
    \rWe\Ind{\dmEv{\Di}{l}}
    \leq \exp\big(\rWn\big\{\tfrac{11}{2}\VnormLp{\txdfPr[l]-\xdfPr[l]}^2-\tfrac{1}{2}\VnormLp{\ProjC[0]\xdf}^2(\sbFxdf[k]-\sbFxdf[l])+(\penSv[l]-\penSv[k])\big\}\big)\Ind{\dmEv{k}{l}}
  \end{multline*}
  Considering  $\dmEv{\Di}{l}:=\setB{7\VnormLp{\txdfPr[l]-\xdfPr[l]}^2<\penSv[l]}$
  the last bound implies
  \begin{multline*}
    \rWe\Ind{\setB{7\VnormLp{\txdfPr[l]-\xdfPr[l]}^2<\penSv[l]}}
    \leq \exp\big(\rWn\big\{\tfrac{11}{14}\penSv[l]-\tfrac{1}{2}\VnormLp{\ProjC[0]\xdf}^2(\sbFxdf[k]-\sbFxdf[l])+(\penSv[l]-\penSv[k])\big\}\big).
 \end{multline*}
 Rearranging the arguments of the last upper bound we obtain the
  assertion \ref{re:rWe:i}.
  Consider secondly \ref{re:rWe:ii} $\Di\in\nsetlo{l,n}$. From
  \cref{re:contr} \ref{re:contr:e2} (with  $\dxdf:=\txdfPr[\ssY]$ and $\pxdf:=\xdf$) and
  \eqref{re:rWe:pro1} follows
  \begin{multline*}
    \rWe[k]\Ind{\dmEv{l}{k}}
    \leq
    \exp\big(\rWn\big\{\tfrac{7}{2}\VnormLp{\txdfPr[k]-\xdfPr[k]}^2
    +\tfrac{3}{2}\VnormLp{\ProjC[0]\xdf}^2(\bias[l]^2(\xdf)-\bias^2(\xdf))
    +(\penSv[l]-\penSv)\big\}\big)\Ind{\dmEv{l}{k}}.
  \end{multline*}
  Setting $\dmEv{l}{\Di}:=\{7\VnormLp{\txdfPr-\xdfPr}^2<\penSv\}$ and exploiting $\sbFxdf\geq0$ we obtain  \ref{re:rWe:ii},
  which completes the proof.\proEnd
\end{pro}

\begin{lem}\label{ak:re:SrWe:ag}Consider  \pcw $\rWe[]$
  as in \eqref{ak:de:rWe} and penalties $(\penSv)_{\Di\in\nset{\ssY}}$ as in \eqref{ak:de:penSv}.
  For any
  $\mdDi,\pdDi\in\nset{\ssY}$ and associated $\pDi,\mDi\in\nset{\ssY}$
  as in \eqref{ak:de:*Di:ag} hold
  \begin{resListeN}
  \item\label{ak:re:SrWe:ag:i}
    $\FuVg{\rWe[]}(\nsetro{\mDi})\leq
    \tfrac{1}{\rWc\cpen}\Ind{\{\mDi>1\}}
    \exp\big(-\tfrac{3\rWc\cpen}{14}\DipenSv[\mdDi]\big)+\Ind{\setB{\VnormLp{\txdfPr[\mdDi]-\xdfPr[\mdDi]}^2\geq\penSv[\mdDi]/7}}$;
  \item\label{ak:re:SrWe:ag:ii}
    $\sum_{\Di\in\nsetlo{\pDi,\ssY}}\penSv\rWe
    \Ind{\{\VnormLp{\txdfPr[\Di]-\xdfPr[\Di]}^2<\penSv/7\}}
    \leq\tfrac{16}{\rWc} \ssY^{-1}$.
  \end{resListeN}
\end{lem}

\begin{pro}[Proof of \cref{ak:re:SrWe:ag}.]
  Consider \ref{ak:re:SrWe:ag:i}. Let $\mDi\in\nset{\mdDi}$ as in
  \eqref{ak:de:*Di:ag}. For the non trivial case $\mDi>1$ from
  \cref{re:rWe} \ref{re:rWe:i} with $l=\mdDi$ follows for all
  $\Di<\mDi\leq \mdDi$
  \begin{multline*}
    \rWe\Ind{\setB{\VnormLp{\txdfPr[\mdDi]-\xdfPr[\mdDi]}^2<\penSv[\mdDi]/7}}
    \\\leq
    \exp\big(\rWn\big\{-\tfrac{1}{2}\VnormLp{\ProjC[0]\xdf}^2\sbF
    +(\tfrac{25}{14}\penSv[\mdDi]
    +\tfrac{1}{2}\VnormLp{\ProjC[0]\xdf}^2\sbF[\mdDi])
    -\penSv\big\}\big),
  \end{multline*}
  and hence by exploiting the definition \eqref{ak:de:*Di:ag} of $\mDi$, that is
  $\VnormLp{\ProjC[0]\xdf}^2\sbF\geq \VnormLp{\ProjC[0]\xdf}^2\sbF[{(\mDi-1)}]>
  \VnormLp{\ProjC[0]\xdf}^2\sbF[\mdDi]+4\penSv[\mdDi]$,
  we obtain for each $\Di\in\nsetro{\mDi}$
  \begin{equation*}
    \rWe\Ind{\setB{\VnormLp{\txdfPr[\mdDi]-\xdfPr[\mdDi]}^2<\penSv[\mdDi]/7}}
    \leq\exp\big(-\tfrac{3}{14}\rWn\penSv[\mdDi]-\rWn\penSv\big).
  \end{equation*}
  The last upper bound together with
  $\penSv=\cpen \DipenSv \ssY^{-1}\geq \cpen\Di\ssY^{-1}$,
  $\Di\in\nset{n}$, as   in \eqref{ak:de:LiSy} gives
  \begin{multline*}
    \FuVg{\rWe[]}(\nsetro{\mDi})\leq
    \FuVg{\rWe[]}(\nsetro{\mDi})
    \Ind{\setB{\VnormLp{\txdfPr[\mdDi]-\xdfPr[\mdDi]}^2<\penSv[\mdDi]/7}}
    +\Ind{\setB{\VnormLp{\txdfPr[\mdDi]-\xdfPr[\mdDi]}^2\geq\penSv[\mdDi]/7}}\\
    \hfill\leq\exp\big(-\tfrac{3\rWc}{14}\ssY\penSv[\mdDi]\big)
    \sum_{k\in\nsetro{\mDi}}\exp(-\rWc\cpen\Di)
    +\Ind{\setB{\VnormLp{\txdfPr[\mdDi]-\xdfPr[\mdDi]}^2\geq\penSv[\mdDi]/7}}
  \end{multline*}
  which combined with $\sum_{\Di\in\Nz}\exp(-\mu\Di)\leq \mu^{-1}$ for any $\mu>0$
  implies \ref{ak:re:SrWe:ag:i}.  Consider \ref{ak:re:SrWe:ag:ii}. Let $\pDi\in\nset{\pdDi,\ssY}$ as in \eqref{ak:de:*Di:ag}. For the non trivial case $\pDi<\ssY$ from \cref{re:rWe} \ref{re:rWe:ii}
  with $l=\pdDi$ follows for all $\Di>\pDi\geq \pdDi$
  \begin{equation*}
    \rWe\Ind{\setB{\VnormLp{\txdfPr-\xdfPr}^2<\penSv/7}}
    \leq \exp\big(\rWn\big\{-\tfrac{1}{2}\penSv
    +\tfrac{3}{2}\VnormLp{\ProjC[0]\xdf}^2\sbF[\pdDi]
    +\penSv[\pdDi]\big\}\big),
  \end{equation*}
  and hence by employing the definition \eqref{ak:de:*Di:ag} of $\pDi$, that is,
  $\tfrac{1}{4}\penSv\geq \tfrac{1}{4}\penSv[(\pDi+1)] >
  \penSv[\pdDi]$ $+\tfrac{3}{2}\VnormLp{\ProjC[0]\xdf}^2\sbF[\pdDi]$, we
  obtain for each $\Di\in\nsetlo{\pDi,\ssY}$
  \begin{equation*}
    \rWe\Ind{\setB{\VnormLp{\txdfPr-\xdfPr}^2<\penSv/7}}
    \leq \exp\big(\rWn\big\{-\tfrac{1}{4} \penSv\big\}\big).
  \end{equation*}
   Consequently, using $\penSv=\cpen \Di\cmiSv\miSv \ssY^{-1}$, $\Di\in\nset{\ssY}$, as   in \eqref{ak:de:LiSy} implies
  \begin{equation}\label{ak:re:SrWe:ag:pe1}
    \sum_{\Di\in\nsetlo{\pDi,\ssY}}\penSv\rWe\Ind{\{\VnormLp{\txdfPr-\xdfPr}^2<\pen/7\}}
    \leq \cpen\ssY^{-1}\sum_{\Di\in\nsetlo{\pDi,n}} \Di\cmiSv\miSv\exp\big(-\tfrac{\rWc\cpen}{4}\Di\cmiSv\miSv\big)
  \end{equation}
  Exploiting that
  $(\cmiSv)^{1/2}=\tfrac{\log (\Di\miSv \vee
    (\Di+2))}{\log(\Di+2)}\geq1$, $\Di\miSv\leq \exp
  ((\cmiSv)^{1/2}\log(\Di+2))$ for each $\Di\in\Nz$,
  $\cpen/4\geq2\log(3e)$ and $\rWc\geq1$ for all $k\in\Nz$ holds
  $\tfrac{\rWc\cpen}{4} k-\log(k+2)\geq1$. Making further use of the
  elementary inequality $a\exp(-ab)\leq \exp(-b)$ for $a,b\geq1$ it
  follows
  \begin{multline*}
    \cmiSv\Di \miSv\exp\big(-\tfrac{\rWc\cpen}{4}\cmiSv\Di\miSv\big)
    \leq\cmiSv\exp\big(-\tfrac{\rWc\cpen}{4}\cmiSv\Di\miSv
    + \sqrt{\cmiSv}\log(\Di+2)\big)
    \\\hfill\leq
    \cmiSv\exp\big(-\cmiSv(\tfrac{\rWc\cpen}{4}\Di-\log(\Di+2))\big)
    \leq\exp\big(-(\tfrac{\rWc\cpen}{4}\Di-\log(\Di+2))\big)\\
    =(\Di+2)\exp\big(-\tfrac{\rWc\cpen}{4}\Di\big).
  \end{multline*}
  which  with $\sum_{\Di\in\Nz}\mu\Di\exp(-\mu\Di)\leq 2$ and
  $\sum_{\Di\in\Nz}\mu\exp(-\mu\Di)\leq 1$ for any $\mu>1$ implies
  \begin{displaymath}
    \sum_{\Di\in\nsetlo{\pDi,\ssY}}\cmiSv\Di \miSv\exp\big(-\tfrac{\rWc\cpen}{4}\cmiSv\Di\miSv\big)
    \leq \sum_{\Di=\pDi+1}^\infty(\Di+2)\exp\big(-\tfrac{\rWc\cpen}{4}\Di\big)
    \leq \tfrac{16}{\cpen\rWc}.
  \end{displaymath}
  Combining the last bound and \eqref{ak:re:SrWe:ag:pe1} we obtain
  assertion \ref{ak:re:SrWe:ag:ii}, which completes the proof.\proEnd
\end{pro}

\begin{te}
  The next result can be directly deduced from \cref{ak:re:SrWe:ag} by letting $\rWc\to\infty$. However, we
  think the direct proof given in \cref{ak:re:SrWe:ms} provides an interesting illustration of the values
  $\pDi,\mDi\in\nset{\ssY}$ as defined in \eqref{ak:de:*Di:ag}.
\end{te}

\begin{lem}\label{ak:re:SrWe:ms}Consider  \msw $\msWe[]$
  as in \eqref{ak:de:msWe} and penalties $(\penSv)_{\Di\in\nset{\ssY}}$ as in \eqref{ak:de:penSv}.  For any $\mdDi,\pdDi\in\nset{\ssY}$ and associated
  $\pDi,\mDi\in\nset{n}$ as in \eqref{ak:de:*Di:ag} hold
  \begin{resListeN}[]
  \item\label{ak:re:SrWe:ms:i}
    $\FuVg{\msWe[]}(\nsetro{\mDi})\Ind{\{\VnormLp{\xdfPr[\mdDi]-\xdfPr[\mdDi]}^2
      <\penSv[\mdDi]/7\}}=0$;
  \item\label{ak:re:SrWe:ms:ii}
    $\sum_{\Di\in\nsetlo{\pDi,\ssY}}\penSv\msWe\Ind{\{\VnormLp{\txdfPr-\xdfPr}^2<\penSv/7\}}=0$.
  \end{resListeN}
\end{lem}

\begin{pro}[Proof of \cref{ak:re:SrWe:ms}.]
  By definition of $\hDi$ it holds
  $-\VnormLp{\txdfPr[\hDi]}^2+\penSv[\hDi]\leq
  -\VnormLp{\txdfPr}^2+\penSv$ for all $\Di\in\nset{\ssY}$, and
  hence
  \begin{equation}\label{ak:re:SrWe:ms:pr:e1}
    \VnormLp{\txdfPr[\hDi]}^2-\VnormLp{\txdfPr}^2\geq
    \penSv[\hDi]-\penSv\text{ for all }\Di\in\nset{\ssY}.
  \end{equation}
  Consider \ref{ak:re:SrWe:ms:i}. Let $\mDi\in\nset{\mdDi}$ as in \eqref{ak:de:*Di:ag}. For the non trivial case $\mDi>1$ it is sufficient to show, that
  $\{\hDi\in\nsetro{\mDi}\}\subseteq
  \{\VnormLp{\txdfPr-\xdfPr}^2\geq\penSv[\mdDi]/7\}$ holds.  On the event $\{\hDi\in\nsetro{\mDi}\}$ we have
  $1\leq\hDi<\mDi\leq\mdDi$ and thus the definition
  \eqref{ak:de:*Di:ag} of $\mDi$  implies
  \begin{equation}\label{ak:re:SrWe:ms:pr:e2}
    \VnormLp{\ProjC[0]\xdf}^2\bias[\hDi]^2(\xdf)\geq
    \VnormLp{\ProjC[0]\xdf}^2\bias[(\mDi-1)]^2(\xdf)>
    \VnormLp{\ProjC[0]\xdf}^2\bias[\mdDi]^2(\xdf)+4\penSv[\mdDi].
  \end{equation}
  On the other hand side from \cref{re:contr} \ref{re:contr:e1}  (with  $\dxdf:=\txdfPr[\ssY]$ and $\pxdf:=\xdf$) follows
  \begin{equation}\label{ak:re:SrWe:ms:pr:e3}
    \VnormLp{\txdfPr[\hDi]}^2-\VnormLp{\txdfPr[\mdDi]}^2\leq
    \tfrac{11}{2}\VnormLp{\txdfPr[\mdDi]-\xdfPr[\mdDi]}^2
    -\tfrac{1}{2}\VnormLp{\ProjC[0]\xdf}^2\{\bias[\hDi]^2(\xdf)-\bias[\mdDi]^2(\xdf)\}.
  \end{equation}
  Combining, first \eqref{ak:re:SrWe:ms:pr:e1} and
  \eqref{ak:re:SrWe:ms:pr:e3}, and secondly
  \eqref{ak:re:SrWe:ms:pr:e2} with $\penSv[\hDi]\geq0$ we conclude
  \begin{equation*}
    \tfrac{11}{2}\VnormLp{\txdfPr[\mdDi]-\xdfPr[\mdDi]}^2\geq
    \penSv[\hDi]-\penSv[\mdDi]
    +\tfrac{1}{2}\VnormLp{\ProjC[0]\xdf}^2\{\bias[\hDi]^2(\xdf)-\bias[\mdDi]^2(\xdf)\}>\tfrac{11}{14}\penSv[\mdDi],
  \end{equation*}
  hence
  $\{\hDi\in\nsetro{\mDi}\}\subseteq
  \{\VnormLp{\txdfPr-\xdfPr}^2\geq\penSv[\mdDi]/7\}$, which shows
  \ref{ak:re:SrWe:ms:i}.
  Consider \ref{ak:re:SrWe:ms:ii}.  Let $\pDi\in\nset{\pdDi,\ssY}$ as
  in \eqref{ak:de:*Di:ag}. For the non trivial case $\pDi<\ssY$  it is sufficient to show that,
  $\{\hDi\in\nsetlo{\pDi,\ssY}\}\subseteq
  \{\VnormLp{\txdfPr[\hDi]-\xdfPr[\hDi]}^2\geq\penSv[\hDi]/7\}$.   On the
  event $\{\hDi\in\nsetlo{\pDi,\ssY}\}$ holds $\hDi>\pDi\geq\pdDi$ and
  thus the definition \eqref{ak:de:*Di:ag}  of $\pDi$  implies
  \begin{equation}\label{ak:re:SrWe:ms:pr:e4}
    \penSv[\hDi]\geq \penSv[(\pDi+1)] > 6\VnormLp{\ProjC[0]\xdf}^2\bias[\pdDi]^2(\xdf)+ 4\penSv[\pdDi]
  \end{equation}
  and due to \cref{re:contr} \ref{re:contr:e2}  (with  $\dxdf:=\txdfPr[\ssY]$ and $\pxdf:=\xdf$) also
  \begin{equation}\label{ak:re:SrWe:ms:pr:e5}
    \VnormLp{\txdfPr[\hDi]}^2-\VnormLp{\txdfPr[\pdDi]}^2\leq
    \tfrac{7}{2}\VnormLp{\txdfPr[\hDi]-\xdfPr[\hDi]}^2+\tfrac{3}{2}\VnormLp{\ProjC[0]\xdf}^2
    \{\bias[\pdDi]^2(\xdf)-\bias[\hDi]^2(\xdf)\}.
  \end{equation}
  Combining, first \eqref{ak:re:SrWe:ms:pr:e1} and
  \eqref{ak:re:SrWe:ms:pr:e5}, and secondly
  \eqref{ak:re:SrWe:ms:pr:e4} with $\bias[\hDi]^2(\xdf)\geq0$ it
  follows that
  \begin{multline*}
    \tfrac{7}{2}\VnormLp{\txdfPr[\hDi]-\xdfPr[\hDi]}^2\geq
    \penSv[\hDi]-\penSv[\pdDi]  -\tfrac{3}{2}\VnormLp{\ProjC[0]\xdf}^2
    \{\bias[\pdDi]^2(\xdf)-\bias[\hDi]^2(\xdf)\}>\tfrac{1}{2}\penSv[\hDi]\hfill
  \end{multline*}
  hence $\{\hDi\in\nsetlo{\pDi,\ssY}\}\subseteq
  \{7\VnormLp{\txdfPr[\hDi]-\xdfPr[\hDi]}^2\geq\penSv[\hDi]\}$, which
  shows
 \ref{ak:re:SrWe:ms:ii} and completed the proof.\proEnd
\end{pro}

\begin{lem}\label{ak:re:nd:rest}Consider $(\penSv)_{\Di\in\nset{\ssY}}$ as in \eqref{ak:de:penSv} with $\cpen\geq84$.
Let $\Di_{\ydf}:=\floor{3(400\Vnormlp[1]{\fydf})^2}$ and
    $\ssY_{o}:=15({600})^4$.  There exists a finite numerical constant  $\cst{}>0$ such that for all $\ssY\in\Nz$ and all $\mdDi\in\nset{\ssY}$  hold
\begin{resListeN}
\item\label{ak:re:nd:rest1}
$\sum_{\Di \in \nset{\ssY}}\nEx\vectp{\VnormLp{\txdfPr-\xdfPr}^2-\penSv/7}\leq
\cst{}\ssY^{-1}\big(\miSv[\Di_{\ydf}]\Di_{\ydf}+ \miSv[\ssY_{o}]\big)$;
\item\label{ak:re:nd:rest2}
  $\sum_{\Di\in\nset{\ssY}}\penSv\nVg\big(\VnormLp{\txdfPr-\xdfPr}^2\geq\penSv/7\big)\leq\cst{}\ssY^{-1}\big(\miSv[\Di_{\ydf}]^2\Di_{\ydf}^3+\miSv[\ssY_{o}]^2\big)$;
\item\label{ak:re:nd:rest3}
  $\nVg\big(\VnormLp{\txdfPr[\mdDi]-\xdfPr[\mdDi]}^2\geq\penSv[\mdDi]/7\big)\leq    \cst{} \big(\exp\big(\tfrac{-\cmiSv[\mdDi]\mdDi}{200\Vnormlp[1]{\fydf}}\big)+\ssY^{-1}\big)$.
\end{resListeN}
\end{lem}

\begin{pro}[Proof of \cref{ak:re:nd:rest}.]
  We show below that for  $\DipenSv=\cmSv \Di \miSv$
  with  $\cmiSv\geq1$ as in \eqref{ak:de:LiSy} there is a numerical constant
  $\cst{}$ such that for all $\ssY\in\Nz$ and
  $\Di\in\nset{\ssY}$ hold
  \begin{resListeN}[\renewcommand{\theListeN}{(\alph{ListeN})}]
  \item\label{ak:re:rest:i}
    $\sum_{\Di \in \nset{\ssY}}\nEx  \vectp{\VnormLp{\txdfPr-\xdfPr}^2-12\DipenSv/\ssY}
    \leq \cst{}\ssY^{-1}\big(\miSv[\Di_{\ydf}]\Di_{\ydf}+ \miSv[\ssY_{o}]\big);$
  \item\label{ak:re:rest:ii}
    $\sum_{\Di \in \nset{\ssY}}\DipenSv\nVg\big(\VnormLp{\txdfPr-\xdfPr}^2
    \geq12\DipenSv/\ssY\big)\leq\cst{}\big(\miSv[\Di_{\ydf}]^2\Di_{\ydf}^3+\miSv[\ssY_{o}]^2\big);$
  \item\label{ak:re:rest:iii}
  $\nVg\big(\VnormLp{\txdfPr-\xdfPr}^2 \geq 12\DipenSv/\ssY\big)\leq
    \cst{} \big(\exp\big(\tfrac{-\cmiSv\Di}{200\Vnormlp[1]{\fydf}}\big)+\ssY^{-1}\big).$
  \end{resListeN}
  Since
  $\penSv/7\geq12\DipenSv\ssY^{-1}$ for all
  $\Di\in\nset{\ssY}$ the bounds \ref{ak:re:rest:i},
  \ref{ak:re:rest:ii} and \ref{ak:re:rest:iii}, respectively, imply
  immediately \cref{ak:re:nd:rest} \ref{ak:re:nd:rest1},
  \ref{ak:re:nd:rest2} and \ref{ak:re:nd:rest3}.  In the sequel we use without further reference that $\Di\miSv\leq \exp
 (\sqrt{\cmSv}\log(\Di+2))$ and $\cmiSv\geq1$  for each $\Di\in\Nz$.
  Considering
  \ref{ak:re:rest:i}  we show that
 \begin{equation}\label{ak:re:rest:pr1}
   \sum_{\Di \in \nset{\ssY}}\miSv\exp\big(\tfrac{-\cmSv\Di}{3\Vnormlp[1]{\fydf}}\big)\leq
   9\miSv[\Di_{\ydf}]{\Vnormlp[1]{\fydf}}
   \quad\text{ and }\quad \sum_{\Di \in \nset{\ssY}}\tfrac{\Di\miSv}{\ssY}\exp\big(\tfrac{-\sqrt{n\cmSv}}{200}\big)\leq \miSv[\ssY_{o}]\ssY_{o}.
\end{equation}
hold  for all $\ssY\in\Nz$, where a combination of the last bounds  and
\cref{re:conc} \ref{re:conc:i} implies directly \ref{ak:re:rest:i}.
We decompose the first sum in \eqref{ak:re:rest:pr1}  into two parts which we bound separately.
Exploiting that $\sum_{\Di\in\Nz}\exp(-\mu\Di)\leq \mu^{-1}$
for any $\mu>0$ and setting
  $\tDi_{\ydf}:=\floor{3({6\Vnormlp[1]{\fydf}})^2}$ holds
  \begin{equation}\label{ak:re:rest:pr2}
    \sum_{\Di \in \nset{\tDi_{\ydf}}}\miSv\exp\big(\tfrac{-\cmSv\Di}{3\Vnormlp[1]{\fydf}}\big)
    \leq \miSv[\tDi_{\ydf}]\sum_{\Di \in \nset{\tDi_{\ydf}}}
    \exp\big(\tfrac{-\Di}{3\Vnormlp[1]{\fydf}}\big)
   \leq \miSv[\tDi_{\ydf}]{3\Vnormlp[1]{\fydf}}.
 \end{equation}
  On the other hand for any
  $\Di>\tDi_{\ydf}$
  holds
  $\tfrac{\sqrt{\cmSv}\Di}{6\Vnormlp[1]{\fydf}}\geq\log(\Di+2)$
  implying
  $\miSv\exp\big(\tfrac{-\cmSv\Di}{3\Vnormlp[1]{\fydf}}\big)
    \leq\exp\big(-\tfrac{\Di}{6\Vnormlp[1]{\fydf}}\big)$
    and hence
    $\sum_{\Di \in \nsetlo{\tDi_{\ydf}, \ssY}}\miSv\exp\big(\tfrac{-\cmSv\Di}{3\Vnormlp[1]{\fydf}}\big)\leq
    \sum_{\Di\in \nsetlo{\tDi_{\ydf}, \ssY}}\exp\big(-\tfrac{\Di}{6\Vnormlp[1]{\fydf}}\big)
    \leq {6\Vnormlp[1]{\fydf}}$.
  The last bound, \eqref{ak:re:rest:pr2} and
  $\tDi_{\ydf}\leq\Di_{\ydf}$ imply together
  the first bound in \eqref{ak:re:rest:pr1}.
 Considering the second bound
for $\ssY\in\Nz$  we distinguish the following two
 cases, \begin{inparaenum}[i]\renewcommand{\theenumi}{\dgrau\rm(\alph{enumi})}\item\label{ak:re:rest:p:c1}
$ \ssY>\tn_{o}:=15({200})^4$ and \item\label{ak:re:rest:p:c2}
$\ssY\in\nset{\tn_{o}}$. \end{inparaenum} Firstly, consider
\ref{ak:re:rest:p:c1}, where
 $\sqrt{\ssY}\geq{200}\log(\ssY+2)$
 and hence
 \begin{equation}\label{ak:re:rest:pr3}
   \sum_{\Di \in \nset{\ssY}}\tfrac{\Di\miSv}{\ssY}\exp\big(\tfrac{-\sqrt{n\cmSv}}{200}\big)
   \leq
   \sum_{\Di \in \nset{\ssY}}\tfrac{1}{\ssY}\exp\big(-\sqrt{\cmSv}[\tfrac{\sqrt{\ssY}}{200}-\log(\Di+2)]\big)\leq \sum_{\Di \in \nset{\ssY}}\tfrac{1}{\ssY}=1.
 \end{equation}
Secondly, considering \ref{ak:re:rest:p:c2} $\ssY\in\nset{\tn_{o}}$  holds
   $\sum_{\Di \in \nset{\ssY}}\tfrac{\Di\miSv}{\ssY}\exp\big(\tfrac{-\sqrt{\ssY\cmSv}}{200}\big)
   \leq\tn_{o}\miSv[\tn_{o}]\leq\ssY_{o}\miSv[\ssY_{o}]$, since
 $\miSv[\ssY]\leq\miSv[\tn_{o}]\leq\miSv[\ssY_{o}]$.
  Combining
    \eqref{ak:re:rest:pr3} and the last bound   for the two cases \ref{ak:re:rest:p:c1} $\ssY>\tn_{o}$
 and \ref{ak:re:rest:p:c2} $\ssY\in\nset{\tn_{o}}$ we obtain the
 second bound in \eqref{ak:re:rest:pr1}.
Consider  \ref{ak:re:rest:ii}. We show  that
 \begin{multline}\label{ak:re:rest:pr4}
  \sum_{\Di \in \nset{\ssY}}\Di\cmSv\miSv\exp\big(\tfrac{-\cmSv\Di}{200\Vnormlp[1]{\fydf}}\big)\leq \miSv[\Di_{\ydf}]^2\Di_{\ydf}^3\quad\text{ and }\quad\\  \sum_{\Di \in \nset{\ssY}}\Di\cmSv\miSv\exp\big(\tfrac{-\sqrt{\ssY\cmSv}}{200}\big)\leq \miSv[\ssY_{o}]^2\ssY_{o}^2
\end{multline}
hold for all $\ssY\in\Nz$.  Combining the last bounds  and
\cref{re:conc} \ref{re:conc:ii}  we obtain \ref{ak:re:rest:ii}.
We decompose the first sum in \eqref{ak:re:rest:pr4} into two parts
which we bound separately. Note that $\log(\Di\miSv)\leq
\tfrac{1}{e}\Di\miSv$, and hence $\cmSv\leq\Di\miSv$. Setting $\Di_{\ydf}=\floor{3({400\Vnormlp[1]{\fydf}})^2}$ holds
\begin{multline}\label{ak:re:rest:pr5}
  \sum_{\Di\in \nset{\Di_{\ydf}}}\Di\cmSv\miSv\exp\big(\tfrac{-\cmSv\Di}{200\Vnormlp[1]{\fydf}}\big)\leq
  \cmSv[\Di_{\ydf}]\miSv[\Di_{\ydf}]\Di_{\ydf}\sum_{\Di\in \nset{\Di_{\ydf}}}\exp\big(\tfrac{-\Di}{200\Vnormlp[1]{\fydf}}\big)\\\leq
  \Di_{\ydf}^2\miSv[\Di_{\ydf}]^2({200\Vnormlp[1]{\fydf}})
\end{multline}
On the other hand for any  $\Di\geq 3({400\Vnormlp[1]{\fydf}})^2$ holds
$\Di\geq  ({400\Vnormlp[1]{\fydf}})\log(\Di+2)$, and
hence
$\Di-{200\Vnormlp[1]{\fydf}}\log(\Di+2)\geq{200\Vnormlp[1]{\fydf}}\log(\Di+2)$
or equivalently,
$\tfrac{\Di}{200\Vnormlp[1]{\fydf}}-\log(\Di+2)\geq\log(\Di+2)\geq1$,
which implies
$\Di\cmSv\miSv\exp\big(\tfrac{-\cmSv\Di}{200\Vnormlp[1]{\fydf}}\big)
\leq
(\Di+2)\exp\big(-\tfrac{\Di}{200\Vnormlp[1]{\fydf}}\big)$.
Consequently,  exploiting that for any $\mu>0$ holds
$\sum_{\Di\in\Nz}(\Di+2)\exp(-\mu\Di)\leq \exp(\mu)\mu^{-2}+2\mu^{-1}$ we
obtain
$\sum_{\Di\in \nsetlo{\Di_{\ydf}, \ssY}}\Di\cmSv\miSv\exp\big(\tfrac{-\cmSv\Di}{200\Vnormlp[1]{\fydf}}\big)\leq\exp(\tfrac{1}{200\Vnormlp[1]{\fydf}})(200\Vnormlp[1]{\fydf})^2+2(200\Vnormlp[1]{\fydf})$.
  The last bound and \eqref{ak:re:rest:pr5}  imply together
  the first bound in \eqref{ak:re:rest:pr4}.
 Considering the second bound,
for $\ssY\in\Nz$  we distinguish the following two
 cases, \begin{inparaenum}[i]\renewcommand{\theenumi}{\dgrau\rm(\alph{enumi})}\item\label{ak:re:rest:p:c3}
$ \ssY>n_{o}=15({600})^4$ and \item\label{ak:re:rest:p:c4}
$\ssY\in\nset{n_{o}}$. \end{inparaenum} Firstly, consider
\ref{ak:re:rest:p:c3} ,
where $\sqrt{\ssY}\geq{600}\log(\ssY+2)$, and hence together  with
$\cmSv\leq\Di\miSv$  it follows
\begin{multline}\label{ak:re:rest:pr6}
\sum_{\Di \in \nset{\ssY}}\Di\cmSv\miSv\exp\big(\tfrac{-\sqrt{\ssY\cmSv}}{200}\big)\leq
\sum_{\Di \in \nset{\ssY}}\Di^2\miSv^2\exp\big(\tfrac{-\sqrt{\ssY\cmSv}}{200}\big)\\
\leq\sum_{\Di \in \nset{\ssY}}\tfrac{1}{\ssY}\exp\big(-3\sqrt{\cmSv}[\tfrac{\sqrt{\ssY}}{600}-\log(\ssY+2)]\big)
\leq \sum_{\Di \in \nset{\ssY}}\tfrac{1}{\ssY}=1.
  \end{multline}
Secondly, consider \ref{ak:re:rest:p:c4}. Since  $\ssY^b\exp(-a\ssY^{1/c})\leq
(\tfrac{cb}{ea})^{cb}$ for all $c>0$ and $a,b\geq0$ it follows
\begin{multline*}
\sum_{\Di \in \nset{\ssY}}\Di\cmSv\miSv\exp\big(\tfrac{-\sqrt{\ssY\cmSv}}{200}\big)\leq\ssY^2\cmSv[\ssY]\miSv[\ssY]\exp\big(\tfrac{-\sqrt{\ssY}}{200}\big)\leq
\miSv[\ssY]^2\ssY^3\exp\big(\tfrac{-\sqrt{\ssY}}{200}\big)\\\leq \miSv[n_{o}]^2\big({600}\big)^6\leq\miSv[n_{o}]^2n_{o}^2.
\end{multline*}
  Combining
    \eqref{ak:re:rest:pr6} and the last bound   for the two cases \ref{ak:re:rest:p:c3} $\ssY>n_{o}$
 and \ref{ak:re:rest:p:c4} $\ssY\in\nset{n_{o}}$ we obtain the
 second bound in \eqref{ak:re:rest:pr4}.
Consider \ref{ak:re:rest:iii}. Since
$\tfrac{\sqrt{\ssY\cmiSv}}{200}\geq\tfrac{\sqrt{\ssY}}{200}$
and $\ssY\exp(-\tfrac{\sqrt{\ssY}}{200})\leq(200)^2$
from \cref{re:conc} \ref{re:conc:ii} follows immediately \ref{ak:re:rest:iii},
which  completes the proof.\proEnd\end{pro}

\begin{lem}\label{ak:ag:ub:p}Let the assumptions of
  \cref{ak:ag:ub:pnp} \ref{ak:ag:ub:pnp:p} be satisfied.
There is a finite numeric constant $\cst{}>0$ such that for all
$n\in\Nz$ with $\ssY_{o}:=15(600)^4$ holds
\begin{multline}\label{ak:ag:ub:p:e1}
  \nEx\VnormLp{\txdfAg[{\We[]}]-\xdf}^2\leq
  \cst{}\VnormLp{\ProjC[0]\xdf}^2\big[  \ssY^{-1}\vee\exp\big(\tfrac{-\cmiSv[\sDi{\ssY}]\sDi{\ssY}}{\Di_{\ydf}}\big)\big]\\
  +\cst{}\big([1\vee K\vee  c_{\xdf}K^2\miSv[K]^2
  ](\miSv[1]^2+\VnormLp{\ProjC[0]\xdf}^2)
  +\miSv[\Di_{\ydf}]^2\Di_{\ydf}^3+\miSv[\ssY_{o}]^2 \big)\ssY^{-1}.
\end{multline}
\end{lem}

\begin{pro}[Proof of \cref{ak:ag:ub:p}.]The proof is based on the
  upper bound
\eqref{ak:ag:ub:p3} which holds for any $\mdDi,\pdDi\in\nset{n}$ and associated
$\mDi,\pDi\in\nset{n}$ as defined in  \eqref{ak:de:*Di:ag}.
Consider first the case $K=0$, where $\bF[0]=0$
and hence $\VnormLp{\ProjC[0]\xdf}^2=0$. From \eqref{ak:ag:ub:p3}
follows
 \begin{equation}\label{ak:ag:ub:pnp:p2}
     \nEx\VnormLp{\txdfAg[{\We[]}]-\xdf}^2\leq \tfrac{2}{7}\penSv[\pDi]
    +\cst{}\big(\miSv[\Di_{\ydf}]^2\Di_{\ydf}^3+\miSv[\ssY_{o}]^2 \big)\ssY^{-1}
\end{equation}
Setting  $\pdDi:=1$ it follows from the definition
\eqref{ak:de:*Di:ag} of  $\pDi$ that
$\penSv[\pDi]\leq4\penSv[1]=4\cpen \LiSv[1]\ssY^{-1}$ and
$\LiSv[1]=\cmSv[1]\miSv[1]\leq\miSv[1]^2$. Thereby (keep in mind
$\cpen\geq84$) \eqref{ak:ag:ub:pnp:p2} implies
 \begin{equation}\label{ak:ag:ub:pnp:p3}
     \nEx\VnormLp{\txdfAg[{\We[]}]-\xdf}^2\leq\cst{}\big(\miSv[1]^2+\miSv[\Di_{\ydf}]^2\Di_{\ydf}^3+\miSv[\ssY_{o}]^2 \big)\ssY^{-1}
\end{equation}
Consider now  $K\in\Nz$, and hence
$\VnormLp{\ProjC[0]\xdf}^2\sbF[{[K-1]}]>0$.  Setting
$\ssY_{\xdf}:=[K\vee\gauss{c_{\xdf}\DipenSv[K]}]\in\Nz$ we distinguish for $\ssY\in\Nz$ the following two
 cases, \begin{inparaenum}[i]\renewcommand{\theenumi}{\dgrau\rm(\alph{enumi})}\item\label{ak:ag:ub:pnp:p:c1}
$\ssY\in\nset{\ssY_{\xdf}}$ and \item\label{ak:ag:ub:pnp:p:c2}
$\ssY> \ssY_{\xdf}$. \end{inparaenum} Firstly, consider
\ref{ak:ag:ub:pnp:p:c1} with $\ssY\in\nset{\ssY_{\xdf}}$, then setting $\mdDi:=1$, $\pdDi:=1$ we have
$\mDi=1$, $1\geq\bF[1]$ and from the definition
\eqref{ak:de:*Di:ag} of  $\pDi$ also
$\penSv[\pDi]\leq2(3\VnormLp{\ProjC[0]\xdf}^2\sbF[1]+2\penSv[1])\leq
6\VnormLp{\ProjC[0]\xdf}^2 +4\cpen\miSv[1]^2$. Thereby,  from \eqref{ak:ag:ub:p3}
follows
 \begin{multline*}
     \nEx\VnormLp{\txdfAg[{\We[]}]-\xdf}^2\leq \tfrac{8\cpen}{7}\miSv[1]^2   +\tfrac{26}{7}\VnormLp{\ProjC[0]\xdf}^2
    +\cst{}\big(\miSv[\Di_{\ydf}]^2\Di_{\ydf}^3+\miSv[\ssY_{o}]^2
    \big)\ssY^{-1}\\
    \leq \cst{}\big(\miSv[1]^2\ssY+\VnormLp{\ProjC[0]\xdf}^2\ssY+\miSv[\Di_{\ydf}]^2\Di_{\ydf}^3+\miSv[\ssY_{o}]^2\big)\ssY^{-1}.
\end{multline*}
Moreover, for all $\ssY\in\nset{\ssY_{\xdf}}$ with
$\ssY_{\xdf}=[K\vee\floor{c_{\xdf}\DipenSv[K]}]$ and
$\DipenSv[K]=K\cmSv[K] \miSv[K]\leq K^2\miSv[K]^2$ holds
$\ssY\leq[K\vee c_{\xdf}K^2\miSv[K]^2]$ and thereby,
\begin{equation}\label{ak:ag:ub:pnp:p4}
  \nEx\VnormLp{\txdfAg[{\We[]}]-\xdf}^2\leq
  \cst{}\big([K\vee c_{\xdf}K^2\miSv[K]^2](\miSv[1]^2+\VnormLp{\ProjC[0]\xdf}^2)+\miSv[\Di_{\ydf}]^2\Di_{\ydf}^3+\miSv[\ssY_{o}]^2\big)\ssY^{-1}.
\end{equation}
Secondly, consider \ref{ak:ag:ub:pnp:p:c2}, i.e., $\ssY>
\ssY_{\xdf}$. Setting
$\pdDi:=K\leq [K\vee\floor{c_{\xdf}\DipenSv[K]}]=\ssY_{\xdf}$, i.e.,
$\pdDi\in\nset{\ssY}$, it follows $\bF[\pdDi]=0$ and the
definition \eqref{ak:de:*Di:ag} of $\pDi$ implies
$\penSv[\pDi]\leq4\penSv[\pdDi]=4\cpen\DipenSv[K]\ssY^{-1}\leq4\cpen K^2\miSv[K]^2\ssY^{-1}$. From
\eqref{ak:ag:ub:p3} follows for all $\ssY> \ssY_{\xdf}$ thus
\begin{multline}\label{ak:ag:ub:pnp:p5}
  \nEx\VnormLp{\txdfAg[{\We[]}]-\xdf}^2\leq 2\VnormLp{\ProjC[0]\xdf}^2\sbF[\mDi]
    + \cst{}\VnormLp{\ProjC[0]\xdf}^2\Ind{\{\mDi>1\}}\big[\ssY^{-1}\vee
    \exp\big(\tfrac{-\cmiSv[\mdDi]\mdDi}{\Di_{\ydf}}\big)\big]\\
    +\cst{}\big(K^2\miSv[K]^2+\miSv[\Di_{\ydf}]^2\Di_{\ydf}^3+\miSv[\ssY_{o}]^2 \big)\ssY^{-1}.
  \end{multline}
Note that for all $\ssY>\ssY_{\xdf}$ holds
  $\sDi{\ssY}=\max\{\Di\in\nset{K,\ssY}:\ssY>c_{\xdf}\DipenSv\}$,
  since the defining set containing $K$ is not empty.
Consequently,  $\sDi{\ssY}\geq
K$ and, hence
$\bias[\sDi{\ssY}](\xdf)=0$, and
$\DipenSv[\sDi{\ssY}]\ssY^{-1}<c_{\xdf}^{-1}=\tfrac{\VnormLp{\ProjC[0]\xdf}^2\sbF[(K-1)]}{4\cpen}$,
it follows
$\VnormLp{\ProjC[0]\xdf}^2\sbF[(K-1)]>4\cpen\DipenSv[\sDi{\ssY}]\ssY^{-1}=4\penSv[\sDi{\ssY}]+\VnormLp{\ProjC[0]\xdf}^2\sbF[\sDi{\ssY}]$
and trivially
$\VnormLp{\ProjC[0]\xdf}^2\sbF[{K}]=0<4\penSv[\sDi{\ssY}]+\VnormLp{\ProjC[0]\xdf}^2\sbF[\sDi{\ssY}]$. Therefore,
setting $\mdDi:=\sDi{\ssY}$ the definition \eqref{ak:de:*Di:ag} of
$\mDi$ implies $\mDi=K$ and hence
$\sbF[\mDi]=\sbF[K]=0$. From \eqref{ak:ag:ub:pnp:p5}  follows
now for all $\ssY> \ssY_{\xdf}$ thus
\begin{multline}\label{ak:ag:ub:pnp:p6}
  \nEx\VnormLp{\txdfAg[{\We[]}]-\xdf}^2\leq  \cst{}\VnormLp{\ProjC[0]\xdf}^2\big[\ssY^{-1}\vee\exp\big(-\tfrac{\cmiSv[\sDi{\ssY}]\sDi{\ssY}}{\Di_{\ydf}}\big)\big]\\
  +\cst{}\big(K^2\miSv[K]^2 +\miSv[\Di_{\ydf}]^2\Di_{\ydf}^3+\miSv[\ssY_{o}]^2 \big)\ssY^{-1}.
\end{multline}
Combining  \eqref{ak:ag:ub:pnp:p4} and
    \eqref{ak:ag:ub:pnp:p6}  for $K\geq1$ with \ref{ak:ag:ub:pnp:p:c1}
$\ssY\in\nset{\ssY_{\xdf}}$ and \ref{ak:ag:ub:pnp:p:c2}
$\ssY\geq \ssY_{\xdf}$, respectively, and \eqref{ak:ag:ub:pnp:p3}  for
$K=0$ implies for all $K\in\Nz_0$ and for all $\ssY\in\Nz$ the claim
\eqref{ak:ag:ub:p:e1},
which  completes the
proof of \cref{ak:ag:ub:p}.\proEnd\end{pro}
\subsubsection{Proof of \cref{ak:ag:ub:pnp:mm} and \cref{ak:ag:ub2:pnp:mm}}\label{a:ak:mrb}
\begin{te}We present first the main arguments to prove
  \cref{ak:ag:ub:pnp:mm} which makes use of \cref{ak:re:nd:rest:mm}
  deferred to the end of this section.\
\end{te}

\begin{te}
  Considering an aggregation
  $\txdf[{\We[]}]=\sum_{\Di \in \nset{\ssY}} \We\txdfPr$  using either
  Bayesian weights $\We[]:=\rWe[]$
  as in \eqref{ak:de:rWe} or model selection weights $\We[]:=\msWe[]$
  as in \eqref{ak:de:msWe} we make use of the upper  bounds
  \eqref{co:agg:ag} and \eqref{co:agg:ms}, respectively. In
  \cref{ak:re:nd:rest:mm} we bound the
  last three terms in \eqref{co:agg:ag} and \eqref{co:agg:ms}  uniformly over $\rwCxdf$ and $\rwCedf$. Moreover, we
  note that  the definition \eqref{ak:de:*Di:ag} of $\pDi$
  and $\mDi$ implies
  $\penSv[\pDi]
  \leq(6\xdfCr+4\cpen\zeta_{\edfCr})\dmRaL{\pdDi}$ and
  $\VnormLp{\ProjC[0]\xdf}^2\sbFxdf[\mDi]
  \leq(\xdfCr+4\cpen\zeta_{\edfCr})\dmRaL{\mdDi}$. Combining \eqref{co:agg:ag} and \eqref{co:agg:ms}, the last
  bounds, $\VnormLp{\ProjC[0]\xdf}^2\leq\xdfCr$, $\rWc\geq1$, $\tfrac{3\rWc\cpen}{14}\DipenSv[\mdDi]\geq\tfrac{1}{\Di_{\xdfCw[]\edfCw[]}}\liCw[\mdDi]\mdDi$  and
  \cref{ak:re:nd:rest:mm} we obtain for all
  $\ssY\in\Nz$
\begin{multline}\label{co:agg:mm}
    \nmRi{\txdfAg[{\We[]}]}{\rwCxdf}{\rwCedf}
    \leq
    \tfrac{2}{7}(6\xdfCr+4\cpen\zeta_{\edfCr})\dmRaL{\pdDi}
    +2(\xdfCr+4\cpen)\dmRaL{\mdDi} \\\hfill+ \cst{} \xdfCr\exp\big(\tfrac{-\liCw[\mdDi]\mdDi}{\Di_{\xdfCw[]\edfCw[]}}\big)
     +\cst{}\ssY^{-1}\{\xdfCr+\edfCr^2\big(\edfCw[\Di_{\xdfCw[]\edfCw[]}]^2\Di_{\xdfCw[]\edfCw[]}^3+\edfCw[\ssY_{o}]^2\big)\}.
  \end{multline}
For $\oDi{\ssY}:=\nmDiL\in\nset{n}$ and $\dmRaL{\Di}$
as in \eqref{oo:de:doRao}  we set $\pdDi:=\oDi{\ssY}$, then for all
$\mdDi\in\nset{\ssY}$ holds
$\dmRaL{\mdDi}\geq\dmRaL{\oDi{\ssY}}=\nmRaL=\min\Nset[\Di\in\Nz]{\dmRaL{\Di}}\geq\ssY^{-1}$. Combining
the last bound  and \eqref{co:agg:mm}  implies the assertion
\eqref{ak:ag:ub:pnp:mm:e1}, that is for all $\mdDi\in\nset{\ssY}$ holds
\begin{multline}\label{ak:ag:ub:pnp:mm:p1}
    \nmRi{\txdfAg[{\We[]}]}{\rwCxdf}{\rwCedf}
    \leq\cst{}(\xdfCr+\zeta_{\edfCr})\big[\dmRaL{\mdDi}\vee\exp\big(\tfrac{-\liCw[\mdDi]\mdDi}{\Di_{\xdfCw[]\edfCw[]}}\big)\big]\\\hfill+
     +\cst{}\ssY^{-1}\{\xdfCr+\edfCr^2\big(\edfCw[\Di_{\xdfCw[]\edfCw[]}]^2\Di_{\xdfCw[]\edfCw[]}^3+\edfCw[\ssY_{o}]^2\big)\}
\end{multline}
with $\ssY_{o}=15(600)^4$,  which completes the proof of \cref{ak:ag:ub:pnp:mm}.
\end{te}

\begin{pro}[Proof of \cref{ak:ag:ub2:pnp:mm}.]
  Under  \ref{ak:ag:ub2:pnp:mm:c} for $\oDi{\ssY}:=\nmDiL$      as in \eqref{oo:de:doRao}  holds  $\exp\big(-\liCw[\oDi{\ssY}]\oDi{\ssY}/\Di_{\xdfCw[]\edfCw[]}\big)\leq
  \nmRaL$ while for $\ssY\in\nset{\ssY_{\xdfCw[],\edfCw[]}}$ we have
  $\exp\big(-\liCw[\oDi{\ssY}]\oDi{\ssY}/\Di_{\xdfCw[]\edfCw[]}\big)\leq1\leq
  \ssY\nmRaL\leq \ssY_{\xdfCw[],\edfCw[]}
  \nmRaL$. Thereby, from \eqref{ak:ag:ub:pnp:mm:e1} with
  $\nmRaL=\min_{\Di\in\nset{\ssY}}\dmRaL{\Di}$
  follows immediately the claim, which completes the proof of
  \cref{ak:ag:ub2:pnp:mm}.\proEnd\end{pro}

\begin{co}\label{ak:re:nd:rest:mm}Consider $(\penSv)_{\Di\in\nset{n}}$ as in \eqref{ak:de:penSv} with $\cpen\geq84$.
Let $ \ssY_{o}:=15({600})^4$ and
$\Di_{\xdfCw[]\edfCw[]}:=\floor{3(400)^2\xdfCr\zeta_{\edfCr}\,\Vnormlp[1]{\xdfCwS/\edfCwS}}$.
There exists a finite numerical constant  $\cst{}>0$ such that for
each $\xdf\in\rwCxdf$ and $\edf\in\rwCedf$
and for all $n\in\Nz$ and $\Di\in\nset{\ssY}$  hold
\begin{resListeN}
\item\label{ak:re:nd:rest:mm1}
$\sum_{\Di \in \nset{\ssY}}\nEx\vectp{\VnormLp{\txdfPr-\xdfPr}^2-\penSv/7}\leq
\cst{}\ssY^{-1}\edfCr\big(\edfCw[\Di_{\xdfCw[]\edfCw[]}]\Di_{\xdfCw[]\edfCw[]}+ \edfCw[\ssY_{o}]\big)$;
\item\label{ak:re:nd:rest:mm2}
  $\sum_{\Di \in \nset{\ssY}}\penSv\nVg\big(\VnormLp{\txdfPr-\xdfPr}^2\geq\penSv/7\big)\leq\cst{}\ssY^{-1}\edfCr^2\big(\edfCw[\Di_{\xdfCw[]\edfCw[]}]^2\Di_{\xdfCw[]\edfCw[]}^3+\edfCw[\ssY_{o}]^2\big)$;
\item\label{ak:re:nd:rest:mm3}
  $\nVg\big(\VnormLp{\txdfPr[\Di]-\xdfPr[\Di]}^2\geq\penSv/7\big)\leq   \cst{} \big(\exp\big(\tfrac{-\liCw[\Di]\Di}{\Di_{\xdfCw[]\edfCw[]}}\big)+\ssY^{-1}\big)$.
\end{resListeN}
\end{co}

\begin{pro}[Proof of \cref{ak:re:nd:rest:mm}.]The result follows
  immediately from \ref{ak:re:rest:i}-\ref{ak:re:rest:iii} in the proof
  of \cref{ak:re:nd:rest} by using that for all
  $\xdf\in\rwCxdf$, $\edf\in\rwCedf$ and $\Di\in\Nz$ hold $\edfCr^{-1}\cmiSv[\Di]\geq \zeta_{\edfCr}^{-1}\liCw[\Di]$, $\miSv[\Di]\leq\edfCr\edfCw[\Di]$  and
  $\Vnormlp[1]{\fydf}^2\leq
 \edfCr\Vnormlp[1]{\xdfCwS/\edfCwS}\Vnorm[1/{\xdfCw[]}]{\xdf}^2\leq
  \xdfCr\edfCr\,\Vnormlp[1]{\edfCwS\xdfCwS}$, and
   we
  omit the details.\proEnd\end{pro}

\subsection{Proofs of \cref{au}}\label{a:au}

\begin{pro}[Proof of \cref{co:agg:au}.]
We start the proof with the observation that $\fhxdf{0}-\fxdf[0]=0$
and for each $j\in\Zz$ holds
$\ofhxdfPr[{\We[]}]{j}-\ofxdf[j]=\fhxdfPr[{\We[]}]{-j}-\fxdf[-j]$, where
$\fhxdfPr[{\We[]}]{j}-\fxdf[j]=-\fxdf[j]$ for all $|j|>\ssY$, and
\begin{multline*}
  \fhxdfPr[{\We[]}]{j}-\fxdf[j]=
\hfedfmpI[j](\hfydf[j]-\fydf[j])\FuVg{\We[]}(\nset{|j|,n})
+\hfedfmpI[j](\fedf[j]-\hfedf[j])\fxdf[j]\FuVg{\We[]}(\nset{|j|,n})\\
-\Ind{\xEv}\fxdf[j]\FuVg{\We[]}(\nsetro{|j|})-\Ind{\xEv^c}\fxdf[j]\text{ for all }|j|\in\nset{\ssY}
\end{multline*}
with $\xEv:= \{|\hfedf[j]|^2\geq1/\ssE\}$ and
$\xEv^c:= \{|\hfedf[j]|^2<1/\ssE\}$. Consequently,  we  have
  \begin{multline}\label{pro:au:key:e1}
    \VnormLp{\hxdfPr[{\We[]}]-\xdf}^2
\leq3\sum_{|j|\in\nset{\ssY}}|\hfedfmpI[j]|^2|\hfydf[j]-\fydf[j]|^2\FuVg{\We[]}(\nset{|j|,n})
\\\hfill+3\sum_{|j|\in\nset[1,]{\ssY}}\Ind{\xEv}|\fxdf[j]|^2\FuVg{\We[]}(\nsetro[1,]{|j|})+\sum_{|j|>n}|\fxdf[j]|^2\\\hfill
+3\sum_{|j|\in\nset[1,]{\ssY}}|\hfedfmpI[j]|^2|\fedf[j]-\hfedf[j]|^2|\fxdf[j]|^2
+\sum_{|j|\in\nset[1,]{\ssY}}\Ind{\xEv^c}|\fxdf[j]|^2.
\end{multline}
where we consider the first  and the second and third term on the
right hand side separately. Considering the first term
from $\VnormLp{\hxdfPr-\pxdfPr}^2=\sum_{|j|\in\nset[1,]{\Di}}|\hfedfmpI[j]|^2|\hfydf[j]-\fydf[j]|^2$ follows
\begin{multline}\label{pro:au:key:e2}
\sum_{|j|\in\nset[1,]{\ssY}}|\hfedfmpI[j]|^2(\hfydf[j]-\fydf[j])^2
\FuVg{\We[]}(\nset{|j|,\ssY})\\
\hfill\leq\VnormLp{\hxdfPr[\pDi]-\pxdfPr[\pDi]}^2
+\sum_{l\in\nsetlo{\pDi,\ssY}}\We[l]\vectp{\VnormLp{\hxdfPr[l]-\pxdfPr[l]}^2-\pen[l]/7}\\
+\tfrac{1}{7}\sum_{l\in\nsetlo{\pDi,\ssY}}\We[l]\pen[l]\Ind{\{\VnormLp{\hxdfPr[l]-\pxdfPr[l]}^2\geq\pen[l]/7\}}
+\tfrac{1}{7}\sum_{l\in\nsetlo{\pDi,\ssY}}\We[l]\pen[l]\Ind{\{\VnormLp{\hxdfPr[l]-\pxdfPr[l]}^2<\pen[l]/7\}}
\end{multline}
Considering the second and third term  we split the first sum into two parts and obtain
\begin{multline}\label{pro:au:key:e3}
\sum_{|j|\in\nset[1,]{\ssY}}\Ind{\xEv}|\fxdf[j]|^2\FuVg{\We[]}(\nsetro[1,]{|j|})+\sum_{|j|>n}|\fxdf[j]|^2\\
\hspace*{5ex}\leq  \sum_{|j|\in\nset[1,]{\mDi}}|\fxdf[j]|^2\Ind{\xEv}\FuVg{\We[]}(\nsetro[1,]{|j|})+ \sum_{|j|\in\nsetlo{\mDi,n}}|\fxdf[j]|^2+
 2\sum_{|j|>n}|\fxdf[j]|^2\\\hfill
\leq \VnormLp{\ProjC[0]\xdf}^2\{\FuVg{\We[]}(\nsetro[1,]{\mDi})+\sbFxdf[\mDi]\}
\end{multline}
Combining  \eqref{pro:au:key:e1} and \eqref{pro:au:key:e2}, \eqref{pro:au:key:e3} we obtain   the assertion, which completes the proof.\proEnd
\end{pro}
\subsubsection{Proof of \cref{au:ag:ub:pnp} and \cref{au:ag:ub2:pnp}}\label{a:au:rb}

\begin{te}
  We present first the main arguments of the proof of
  \cref{au:ag:ub:pnp}. More technical details are gathered in
  \cref{au:re:SrWe:ag,au:re:SrWe:ms,au:re:nd:rest,au:ag:ub:p} in the end of this section. Keeping in mind the
  definitions \eqref{ak:de:penSv} and \eqref{au:de:peneSv} let us for
  $l\in\nset[1,]{\ssY}$ introduce the event
  $\aixEv[l]:=\setB{1/4\leq\iSv[j]^{-1}\eiSv[j]\leq9/4,\;\forall\;j\in\nset[1,]{l}}$
  and its complement $\aixEv[l]^c$, where due to \cref{re:aixEv}
  holds
  $\tfrac{1}{50}\penSv\Ind{\aixEv[l]}\leq\peneSv\Ind{\aixEv[l]}\leq7\penSv$
  for all $\Di\in\nset[1,]{l}$.\
\end{te}

\begin{te}
  For any $\pdDi,\mdDi\in\nset[1,]{\ssY}$
  (to be choosen suitable below) let us define
  \begin{multline}\label{au:de:*Di:ag}
    \mDi:=\min\set{\Di\in\nset[1,]{\mdDi}:
      \VnormLp{\ProjC[0]\xdf}^2\sbFxdf\leq
      \VnormLp{\ProjC[0]\xdf}^2\sbFxdf[\mdDi]+104\penSv[\mdDi]}
    \quad\text{and}\\\pDi:=\max\set{\Di\in\nset{\pdDi,\ssY}: \peneSv
      \leq
      6\VnormLp{\ProjC[\pdDi]\pxdfPr[\ssY]}^2+4\peneSv[\pdDi]\big)}
  \end{multline}
  where   $\VnormLp{\ProjC[\Di]\pxdfPr[\ssY]}^2
  =\sum_{|j|\in\nsetlo{\Di,\ssY}}\eiSv[j]\iSv[j]^{-1}|\fxdf[j]|^2$ and
  the defining set obviously contains $\mdDi$ and $\pdDi$,
  respectively, and hence, they are not empty. Note that by
  construction the random dimension $\pDi$ is independent of the
  sample $(\rY_i)_{i\in\nset[1,]{\ssY}}$.  We intend to combine the
  upper bound in \eqref{co:agg:au:e1} and the bounds for \pcw
  $\We[]=\erWe[]$ as in \eqref{au:de:erWe} and \msw $\We[]=\msWe[]$ as
  in \eqref{au:de:msWe} given in \cref{au:re:SrWe:ag} and
  \cref{au:re:SrWe:ms}, respectively.  Conditionally on
  $(\rE_i)_{i\in\nset[1,]{\ssE}}$ the r.v.'s
  $(\rY_i)_{i\in\nset[1,]{\ssY}}$ are \iid and we denote by
  $\FuVg[\ssY]{\rY|\rE}$ and $\FuEx[\ssY]{\rY|\rE}$ their joint
  conditional distribution and expectation, respectively.\
\end{te}

\begin{te}
  Exploiting  \cref{au:re:SrWe:ag}  \ref{au:re:SrWe:ag:i} and
  \ref{au:re:SrWe:ag:ii}, where \ref{au:re:SrWe:ag:i} implies
  \begin{multline*}
    \FuEx[\ssY]{\rY|\rE}\FuVg{\erWe[]}(\nsetro[1,]{\mDi})\leq\Ind{\{\mDi>1\}}
    \big(\tfrac{50}{\rWc\cpen}\exp\big(-\tfrac{\rWc\cpen}{2}\DipenSv[\mdDi]\big)\\+
    \FuVg[\ssY]{\rY|\rE}\big(\VnormLp{\hxdfPr[\mdDi]-\pxdfPr[\mdDi]}^2
    \geq\peneSv[\mdDi]/7\big) \Ind{\aixEv[\mdDi]} + \Ind{\aixEv[\mdDi]^c}\big),
  \end{multline*}
   from \eqref{co:agg:au:e1} for  \pcw $\We[]=\erWe[]$ as in \eqref{au:de:erWe} follows immediately
  \begin{multline}\label{co:agg:au:ag}
   \FuEx[\ssY]{\rY|\rE}\VnormLp{\hxdf[{\We[]}]-\xdf}^2\leq
    3\FuEx[\ssY]{\rY|\rE}\VnormLp{\hxdfPr[\pDi]-\pxdfPr[\pDi]}^2
    +3 \VnormLp{\ProjC[0]\xdf}^2\bias[\mDi]^2(\xdf)\\\hfill
    +\tfrac{150}{\rWc\cpen}\VnormLp{\ProjC[0]\xdf}^2\Ind{\{\mDi>1\}}
    \exp\big(-\tfrac{\rWc\cpen}{2}\DipenSv[\mdDi]\big)
    +\ssY^{-1}\tfrac{48}{7\rWc}\\\hfill
    +3 \VnormLp{\ProjC[0]\xdf}^2\Ind{\{\mDi>1\}}
    \big(\FuVg[\ssY]{\rY|\rE}\big(\VnormLp{\hxdfPr[\mdDi]-\pxdfPr[\mdDi]}^2
    \geq\peneSv[\mdDi]/7\big) \Ind{\aixEv[\mdDi]} + \Ind{\aixEv[\mdDi]^c}\big)
    \\\hfill
    +3\sum_{l\in\nsetlo{\pDi,\ssY}}\FuEx[\ssY]{\rY|\rE}
    \vectp{\VnormLp{\hxdfPr[l]-\pxdfPr[l]}^2-\peneSv[l]/7}
    +\tfrac{3}{7}\sum_{l\in\nsetlo{\pDi,\ssY}}\peneSv[l]\FuVg[\ssY]{\rY|\rE}\big(
    \VnormLp{\hxdfPr[l]-\pxdfPr[l]}^2\geq\peneSv[l]/7\big)\\
    +6\sum_{j\in\nset[1,]{\ssY}}|\hfedfmpI[j]|^2|\fedf[j]-\hfedf[j]|^2|\fxdf[j]|^2
    +2\sum_{j\in\nset[1,]{\ssY}}\Ind{\{|\hfedf[j]|^2<1/\ssE\}}|\fxdf[j]|^2.
  \end{multline}
\end{te}

\begin{te}
  On the other hand \eqref{co:agg:au:ag} holds also true  for \msw
  $\We[]=\msWe[]$ by a combination of  the
  upper bound in \eqref{co:agg:au:e1} and the bounds given in
  \cref{au:re:SrWe:ms}.\
\end{te}

\begin{te}
  The deviations of the last three terms in \eqref{co:agg:au:ag}
  we bound in \cref{au:re:nd:rest}, which implies
\end{te}

\begin{te}
    \begin{multline}\label{au:ag:ub:p1}
    \FuEx[\ssY]{\rY|\rE}\VnormLp{\hxdf[{\We[]}]-\xdf}^2\leq
    3\FuEx[\ssY]{\rY|\rE}\VnormLp{\hxdfPr[\pDi]-\pxdfPr[\pDi]}^2
    +3 \VnormLp{\ProjC[0]\xdf}^2\bias[\mDi]^2(\xdf)\\\hfill
    +\cst{} \VnormLp{\ProjC[0]\xdf}^2\Ind{\{\mDi>1\}} \big(
    \tfrac{1}{\rWc}\exp\big(-\tfrac{\rWc\cpen}{2}\DipenSv[\mdDi]\big)+
    \exp\big(\tfrac{-\cmeiSv[\mdDi]\mdDi}{200\Vnormlp[1]{\fydf}}\big)
    \Ind{\aixEv[\mdDi]} + \Ind{\aixEv[\mdDi]^c}\big)
    \\\hfill
    +6\sum_{j\in\nset[1,]{\ssY}}|\hfedfmpI[j]|^2|\fedf[j]-\hfedf[j]|^2|\fxdf[j]|^2
    +2\sum_{j\in\nset[1,]{\ssY}}\Ind{\{|\hfedf[j]|^2<1/\ssE\}}|\fxdf[j]|^2\\
    +\cst{}\ssY^{-1}\big(\tfrac{1}{\rWc}+[1\vee\meiSv[\Di_{\ydf}]^2]\Di_{\ydf}^3
    +[1\vee\meiSv[\ssY_{o}]^2]+\VnormLp{\ProjC[0]\xdf}^2\Ind{\{\mDi>1\}} \big).
  \end{multline}
  Keeping \eqref{ak:de:LiSy} and
  $\peneSv/7\geq12\DipeneSv\ssY^{-1}$ in mind on the one hand holds
  $\FuEx[\ssY]{\rY|\rE}\VnormLp{\hxdfPr[\pDi]-\pxdfPr[\pDi]}^2
  =2\sum_{j=1}^{\pDi}\eiSv[j]/\ssY\leq2\DipeneSv[\pDi]\ssY^{-1}
  \leq\tfrac{1}{42}\peneSv[\pDi]$ and due to \cref{re:aixEv}
  $\FuEx[\ssY]{\rY|\rE}\VnormLp{\hxdfPr[\pDi]-\pxdfPr[\pDi]}^2\leq2\ssE$
  for $\pDi\in\nset{\ssY}$
  due to \cref{re:aixEv}
  \ref{re:aixEv:ii}.   Consequently,
  $\FuEx[\ssY]{\rY|\rE}\VnormLp{\hxdfPr[\pDi]-\pxdfPr[\pDi]}^2\leq
  2\ssE\Ind{\aixEv[\pdDi]^c}+\tfrac{1}{42}\peneSv[\pDi]\Ind{\aixEv[\pdDi]}$
  and hence
  $\FuEx[\ssY]{\rY|\rE}\VnormLp{\hxdfPr[\pDi]-\pxdfPr[\pDi]}^2\leq
  2\ssE\Ind{\aixEv[\pdDi]^c}
  +\tfrac{1}{42}(6\VnormLp{\ProjC[\pdDi]\pxdfPr[\ssY]}^2
  +4\peneSv[\pdDi])\Ind{\aixEv[\pdDi]}$ exploiting the definition
  \eqref{au:de:*Di:ag} of $\pDi$. Thereby, with  $\meiSv[j]\leq\ssE$, $j\in\Nz$, $\rWc\geq1$ and
  $\cpen\geq1$ from \eqref{au:ag:ub:p1} follows
  \begin{multline*}
    \FuEx[\ssY]{\rY|\rE}\VnormLp{\hxdf[{\We[]}]-\xdf}^2 \leq
    \tfrac{2}{7}\peneSv[\pdDi]\Ind{\aixEv[\pdDi]}
    +\tfrac{3}{7}\VnormLp{\ProjC[\pdDi]\pxdfPr[\ssY]}^2\Ind{\aixEv[\pdDi]}
    +6\ssE\Ind{\aixEv[\pdDi]^c} +3
    \VnormLp{\ProjC[0]\xdf}^2\bias[\mDi]^2(\xdf)\\\hfill +\cst{}
    \VnormLp{\ProjC[0]\xdf}^2\Ind{\{\mDi>1\}} \big(
    \tfrac{1}{\rWc}\exp\big(-\tfrac{\rWc\cpen}{2}\DipenSv[\mdDi]\big)+
    \exp\big(\tfrac{-\cmeiSv[\mdDi]\mdDi}{200\Vnormlp[1]{\fydf}}\big)
    \Ind{\aixEv[\mdDi]} + \Ind{\aixEv[\mdDi]^c}\big)\\\hfill
    +\cst{}\ssY^{-1}\big([1\vee\meiSv[\Di_{\ydf}]^2]\Di_{\ydf}^3
    \Ind{\aixEv[\Di_{\ydf}]} +
    \Di_{\ydf}^3\ssE^2\Ind{\aixEv[\Di_{\ydf}]^c}+[1\vee\meiSv[\ssY_{o}]^2]
    \Ind{\aixEv[\ssY_{o}]}+\ssE^2\Ind{\aixEv[\ssY_{o}]^c}
    +\VnormLp{\ProjC[0]\xdf}^2\Ind{\{\mDi>1\}}\big)\\\hfill
    +6\sum_{j\in\nset{\ssY}}|\hfedfmpI[j]|^2|\fedf[j]-\hfedf[j]|^2|\fxdf[j]|^2
    +2\sum_{j\in\nset{\ssY}}\Ind{\{|\hfedf[j]|^2<1/\ssE\}}|\fxdf[j]|^2.
  \end{multline*}
  Exploiting \cref{re:aixEv} \ref{re:aixEv:ii},
  $\LiSv[\mdDi]\geq\liSv[\mdDi]$ and
  $\tfrac{\rWc\cpen}{2}>\tfrac{9}{20000\Vnormlp[1]{\fydf}}>\tfrac{1}{\Di_{\ydf}}$
  it  follows
  \begin{multline*}
    \FuEx[\ssY]{\rY|\rE}\VnormLp{\hxdf[{\We[]}]-\xdf}^2
    \leq  2\penSv[\pdDi]
    +\tfrac{3}{7}\VnormLp{\ProjC[\pdDi]\pxdfPr[\ssY]}^2
    +3 \VnormLp{\ProjC[0]\xdf}^2\bias[\mDi]^2(\xdf)\\\hfill
    +\cst{} \VnormLp{\ProjC[0]\xdf}^2\Ind{\{\mDi>1\}} \big(
    \exp\big(\tfrac{-\cmiSv[\mdDi]\mdDi}{\Di_{\ydf}}\big)
    + \Ind{\aixEv[\mdDi]^c}\big)\\\hfill
    +\cst{}\big(\ssE\Ind{\aixEv[\pdDi]^c}+\ssY^{-1}\big(
    \Di_{\ydf}^3\ssE^2\Ind{\aixEv[\Di_{\ydf}]^c}
    +\ssE^2\Ind{\aixEv[\ssY_{o}]^c}\big)\big)\\\hfill
    +6\sum_{j\in\nset{\ssY}}|\hfedfmpI[j]|^2|\fedf[j]-\hfedf[j]|^2|\fxdf[j]|^2
    +2\sum_{j\in\nset{\ssY}}\Ind{\{|\hfedf[j]|^2<1/\ssE\}}|\fxdf[j]|^2\\
    +\cst{}\ssY^{-1}\big(\miSv[\Di_{\ydf}]^2\Di_{\ydf}^3+
    \miSv[\ssY_{o}]^2+\VnormLp{\ProjC[0]\xdf}^2\Ind{\{\mDi>1\}}\big).
  \end{multline*}
  Bounding the second term and the two sums on the right hand side due
  to  \cref{oSv:re}
  implies
  \begin{multline*}
    \nmEx\VnormLp{\hxdf[{\We[]}]-\xdf}^2\leq
    2\penSv[\pdDi]
    +\tfrac{12}{7}\VnormLp{\ProjC[0]\xdf}^2\bias[\pdDi]^2(\xdf)
    +3 \VnormLp{\ProjC[0]\xdf}^2\bias[\mDi]^2(\xdf)\\\hfill
    +\cst{} \VnormLp{\ProjC[0]\xdf}^2\Ind{\{\mDi>1\}} \big(
    \exp\big(\tfrac{-\cmiSv[\mdDi]\mdDi}{\Di_{\ydf}}\big)
    + \mVg(\aixEv[\mdDi]^c)\big)\\\hfill
    + \cst{}\big(\ssE\mVg(\aixEv[\pdDi]^c)
    + \ssY^{-1}\{\Di_{\ydf}^3\ssE^2\mVg(\aixEv[\Di_{\ydf}]^c)
    +\ssE^2\mVg(\aixEv[\ssY_{o}]^c)\}\big)
    \\\hfill
    +\cst{}\Vnorm[1\wedge\iSv/\ssE]{\ProjC[0]\xdf}^2
    +\cst{}\ssY^{-1}\{\miSv[\Di_{\ydf}]^2\Di_{\ydf}^3+\miSv[\ssY_{o}]^2
    +\VnormLp{\ProjC[0]\xdf}^2\Ind{\{\mDi>1\}}\}.
  \end{multline*}
  Due to \cref{re:evrest} \ref{re:evrest:ii} there is a numerical
  constant $\cst{}$ such that for all $\ssE,\Di\in\Nz$ holds
  $\mVg(\aixEv^c)\leq\cst{}\Di\miSv^2\ssE^{-2}$ and hence,
  $\ssE^2\mVg(\aixEv[\Di_{\ydf}]^c)\leq
  \cst{}\Di_{\ydf}\miSv[\Di_{\ydf}]^2$ and
  $\ssE^2\mVg(\aixEv[\ssY_{o}]^c)\leq
  \cst{}\ssY_{o}\miSv[\ssY_{o}]^2$.
  Consequently,  there is a numerical constant $\cst{}>0$ such
  that for all $\ssY,\ssE\in\Nz$ holds
  \begin{multline}\label{au:ag:ub:e1}
    \nmEx\VnormLp{\hxdf[{\We[]}]-\xdf}^2\leq
    2\penSv[\pdDi]
    +\tfrac{12}{7}\VnormLp{\ProjC[0]\xdf}^2\bias[\pdDi]^2(\xdf)+3
    \VnormLp{\ProjC[0]\xdf}^2\bias[\mDi]^2(\xdf)\\\hfill
    +\cst{} \VnormLp{\ProjC[0]\xdf}^2\Ind{\{\mDi>1\}} \big(
    \exp\big(\tfrac{-\cmiSv[\mdDi]\mdDi}{\Di_{\ydf}}\big)
    + \mVg(\aixEv[\mdDi]^c)\big)+ \cst{}\ssE\mVg(\aixEv[\pdDi]^c)
    \\\hfill
    +\cst{}\moRa
    +\cst{}\ssY^{-1}\big(\miSv[\Di_{\ydf}]^2\Di_{\ydf}^4+\miSv[\ssY_{o}]^2
    +\VnormLp{\ProjC[0]\xdf}^2\Ind{\{\mDi>1\}}\big)
  \end{multline}
  (keep in mind that
  $\ssY_{o}$ is a numerical constant).\
\end{te}

\begin{te}
  From the upper bound \eqref{au:ag:ub:e1} for a suitable choice of the
  dimension parameters $\mdDi,\pdDi\in\nset{\ssY}$ we derive separately
  the risk bound in the two cases \ref{au:ag:ub:pnp:p} and
  \ref{au:ag:ub:pnp:np} considered in \cref{au:ag:ub:pnp}.  The
  tedious case-by-case analysis for \ref{au:ag:ub:pnp:p} is deferred
  to \cref{au:ag:ub:p} in the end of this section.\
\end{te}

\begin{te}
  In case \ref{au:ag:ub:pnp:np} we destinguish for $\ssE\in\Nz$ with
  $\ssE_{\iSv}:=\floor{289(\log3)\cmiSv[1]\miSv[1]}$ the following
  two cases,
  \begin{inparaenum}[i]\renewcommand{\theenumi}{\dgrau\rm(\alph{enumi})}
  \item\label{au:ag:ub:pnp:np:m1}
    $\ssE\in\nset{\ssE_{\iSv}}$ and
  \item\label{au:ag:ub:pnp:np:m2}
    $\ssE>\ssE_{\iSv}$.
  \end{inparaenum} Consider firstly the case \ref{au:ag:ub:pnp:np:m1}
  $\ssE\in\nset{\ssE_{\iSv}}$. We set $\pdDi=\mdDi=1$, and hence
  $\mDi=1$, $\sbF[1]\leq1$, $\penSv[1]\leq\cpen\miSv[1]^2\ssY^{-1}$,
  $\miSv[1]^2\leq \miSv[\ssY_{o}]^2$,
  $\ssE_{\iSv}\leq \cst{}\miSv[1]^2$ and due to \cref{re:evrest}
  \ref{re:evrest:ii}
  $\FuVg[\ssE]{\rE}(\aixEv[1]^c)\leq
  \cst{}\miSv[1]^2\ssE^{-2}$. Thereby, from \eqref{au:ag:ub:e1}
  for all $\ssY\in\Nz$ and $\ssE\in\nset[1,]{\ssE_{\iSv}}$ follows
  \begin{multline}\label{au:ag:ub:pnp:p10}
    \nmEx\VnormLp{\hxdf[{\We[]}]-\xdf}^2\leq
    \cst{}\moRa
    + \cst{}[1\vee\VnormLp{\ProjC[0]\xdf}^2]\miSv[1]^2\ssE^{-1}
    +\cst{}\big(\miSv[\Di_{\ydf}]^2\Di_{\ydf}^4+\miSv[\ssY_{o}]^2\big)\ssY^{-1}.
  \end{multline}
  Consider secondly \ref{au:ag:ub:pnp:np:m2} $\ssE>\ssE_{\iSv}$ with
  $\sDi{\ssE}:=\max\{\Di\in\nset{\ssE}:289\log(\Di+2)\cmiSv[\Di]\miSv[\Di]\leq\ssE\}$.
  For each
  $\Di\in\nset{\sDi{\ssE}}$ holds
  $\ssE\geq289\log(\Di+2)\cmiSv\miSv$, and thus from \cref{re:evrest}
  \ref{re:evrest:iii} follows
  $\FuVg[\ssE]{\rE}(\aixEv[\Di]^c)\leq 11226\ssE^{-2}$. For
  $\oDi{\ssY}:=\noDiL\in\nset{\ssY}$ as in \eqref{oo:de:doRao} setting
  $\pdDi:=\oDi{\ssY}\wedge\sDi{\ssE}$, where
  $\ssE\FuVg[\ssE]{\rE}(\aixEv[\pdDi]^c)\leq \cst{}\ssE^{-1}$,
  $\penSv[\pdDi]\leq
  \cpen\noRaL$
  and $\sbFxdf[\pdDi]\leq \noRaL+\sbFxdf[\sDi{\ssE}]$,
   \eqref{au:ag:ub:e1} implies
  \begin{multline}\label{au:ag:ub:pnp:p11}
    \nmEx\VnormLp{\hxdf[{\We[]}]-\xdf}^2\leq
    \cst{}[1\vee\VnormLp{\ProjC[0]\xdf}^2]\noRaL
    +\tfrac{12}{7}\VnormLp{\ProjC[0]\xdf}^2\sbFxdf[\sDi{\ssE}]\\\hfill
    +3 \VnormLp{\ProjC[0]\xdf}^2\sbFxdf[\mDi]
    +\cst{}\VnormLp{\ProjC[0]\xdf}^2\Ind{\{\mDi>1\}}\big(
    \exp\big(\tfrac{-\cmiSv[\mdDi]\mdDi}{\Di_{\ydf}}\big)
    +  \FuVg[\ssE]{\rE}(\aixEv[\mdDi]^c) \big)
    \\\hfill
    +\cst{}\moRa+\cst{}\ssE^{-1}
    +\cst{}\ssY^{-1}\big(\miSv[\Di_{\ydf}]^2\Di_{\ydf}^4
    +\miSv[\ssY_{o}]^2
    +\VnormLp{\ProjC[0]\xdf}^2\Ind{\{\mDi>1\}}\big).
  \end{multline}
  Let
  $\sDi{\ssY}:=\argmin\{\doRaL{\Di}\vee
  \exp\big(\tfrac{-\cmiSv[\Di]\Di}{\Di_{\ydf}}\big):\Di\in\nset{\ssY}\}$. Setting
  $\mdDi:=\sDi{\ssY}\wedge\sDi{\ssE}$ from \cref{re:evrest}
  \ref{re:evrest:iii} follows
  $\FuVg[\ssE]{\rE}(\aixEv[\mdDi]^c)\leq 53\ssE^{-1}$, while $\mDi$ as
  in definition \eqref{au:de:*Di:ag} satisfies
  \begin{multline*}
    \VnormLp{\ProjC[0]\xdf}^2\sbFxdf[\mDi] \leq
    \VnormLp{\ProjC[0]\xdf}^2\sbFxdf[\mdDi]+104\penSv[\mdDi]
    \\\leq
    \VnormLp{\ProjC[0]\xdf}^2\sbFxdf[\sDi{\ssE}]
    +(\VnormLp{\ProjC[0]\xdf}^2+104\cpen)\doRaL{\sDi{\ssY}},
  \end{multline*}
  where $\ssY^{-1}\leq\noRaL\leq\doRaL{\sDi{\ssY}}$ by
  \eqref{oo:de:doRao} and
  $\moRa\geq\tfrac{1}{2}\VnormLp{\ProjC[0]\xdf}^2 \ssE^{-1}$ (see
  \cref{oo:rem:nm}).  Thereby, we obtain from \eqref{au:ag:ub:pnp:p11}
  for all $\ssY\in\Nz$ and $\ssE>\ssE_{\iSv}$
  \begin{multline}\label{au:ag:ub:pnp:p12}
    \nmEx\VnormLp{\hxdf[{\We[]}]-\xdf}^2
    \leq\cst{}[1\vee\VnormLp{\ProjC[0]\xdf}^2]
    \min_{\Di\in\nset{\ssY}}\{\big[\doRaL{\Di}\vee
    \exp\big(\tfrac{-\cmiSv[\Di]\Di}{\Di_{\ydf}}\big)\big]\}\\\hfill
    +\cst{}\VnormLp{\ProjC[0]\xdf}^2\big[\sbFxdf[\sDi{\ssE}]
    \vee\exp\big(\tfrac{-\cmiSv[\sDi{\ssE}]\sDi{\ssE}}{\Di_{\ydf}}\big)\big]
    +\cst{}\moRa
    \\\hfill
    + \cst{}\ssE^{-1}
    +\cst{}\ssY^{-1}\big(\miSv[\Di_{\ydf}]^2\Di_{\ydf}^4+\miSv[\ssY_{o}]^2\big).
  \end{multline}
  Combining \eqref{au:ag:ub:pnp:p10} and \eqref{au:ag:ub:pnp:p12} for
  the cases \ref{au:ag:ub:pnp:np:m1}
  and
  \ref{au:ag:ub:pnp:np:m2}
   for all $\ssY,\ssE\in\Nz$ holds
  \begin{multline}\label{au:ag:ub:pnp:p13}
    \nmEx\VnormLp{\hxdf[{\We[]}]-\xdf}^2\leq
    \cst{}[1\vee\VnormLp{\ProjC[0]\xdf}^2]
    \min_{\Di\in\nset{\ssY}}\{[\doRaL{\Di}\vee
    \exp\big(\tfrac{-\cmiSv[\Di]\Di}{\Di_{\ydf}}\big)]\}
    \Ind{\{\ssE>\ssE_{\iSv}\}}\\\hfill
    +\cst{}\VnormLp{\ProjC[0]\xdf}^2\big[\bias[\sDi{\ssE}]^2(\xdf)
    \vee\exp\big(\tfrac{-\cmiSv[\sDi{\ssE}]\sDi{\ssE}}{\Di_{\ydf}}\big)\big]
    \Ind{\{\ssE>\ssE_{\iSv}\}}
    +\cst{}\moRa  \\\hfill
    + \cst{}[1\vee\VnormLp{\ProjC[0]\xdf}^2]\miSv[1]^2\ssE^{-1}
    +\cst{}\{\miSv[\Di_{\ydf}]^2\Di_{\ydf}^4
    +\miSv[\ssY_{o}]^2\}\ssY^{-1},
  \end{multline}
  which shows \eqref{au:ag:ub:pnp:e2} and completes the proof of \cref{au:ag:ub:pnp}.\proEnd
\end{te}

\begin{pro}[Proof of \cref{au:ag:ub2:pnp}.]
  Consider the case \ref{au:ag:ub2:pnp:p}. In the proof of
  \cref{au:ag:ub2:pnp} we have shown, that under the  additional
  assumption \ref{ak:ag:ub2:pnp:pc} holds
  $\noRi{\txdfAg[{\We[]}]}{\xdf}{\iSv} \leq \cst{\xdf,\iSv}\ssY^{-1}$
  for all $\ssY\in\Nz$. If in addition \ref{au:ag:ub2:pnp:pc:b} is
  satisfied for $\sDi{\ssE}$ as in
  \cref{au:ag:ub:pnp}, then we have for all $\ssE>\ssE_{\xdf,\edf}$ trivially
  $\exp\big(\tfrac{-\cmiSv[\sDi{\ssE}]\sDi{\ssE}}{\Di_{\ydf}}\big)\leq\ssE^{-1}
  $ while for $\ssY\in\nset[1,]{\ssE_{\xdf\edf}}$ we
  have
  $\exp\big(\tfrac{-\cmiSv[\sDi{\ssE}]\sDi{\ssE}}{\Di_{\ydf}}\big)\leq
  1\leq \ssE_{\xdf,\edf} \ssE^{-1}$. Combining both bounds we obtain
  the assertion \ref{au:ag:ub2:pnp:p}.   On the other hand side, in case
  \ref{au:ag:ub2:pnp:np} under the  additional assumption
  \ref{ak:ag:ub2:pnp:npc} holds
  $\min_{\Di\in\nset{\ssY}}\setB{\big[\dmRaL{\Di}\vee\exp\big(\tfrac{-\cmiSv\Di}{\Di_{\edf}}\big)\big]}\leq\ssY_{\xdf,\iSv}
  \noRaL $
  (cf. \cref{ak:ag:ub2:pnp} \ref{ak:ag:ub2:pnp:np}). A combination of the
  last bound and
  $\exp\big(\tfrac{-\cmiSv[\sDi{\ssE}]\sDi{\ssE}}{\Di_{\ydf}}\big)\leq
  \ssE_{\xdf,\edf} \ssE^{-1}$ due to \ref{au:ag:ub2:pnp:pc:b} implies
  the assertion \ref{au:ag:ub2:pnp:np}, which completes the proof of
  \cref{au:ag:ub2:pnp}.\proEnd
  \end{pro}

\begin{te}
 Below  we state  and prove the technical
 \cref{au:re:SrWe:ag,au:re:SrWe:ms,au:re:nd:rest} used in the proof of \cref{au:ag:ub:pnp}. The
  proof of \cref{au:re:SrWe:ag} is based on \cref{re:erWe} given first.
\end{te}

\begin{lem}\label{re:erWe}
  Consider \pcw $\erWe[]$ as in
  \eqref{au:de:erWe} and let  $l\in\nset{\ssY}$.
  \begin{resListeN}[]
    \item\label{re:erWe:i} For
    $\aixEv[l]:=\setB{\tfrac{1}{4}\leq\iSv[j]^{-1}\eiSv[j]\leq\tfrac{9}{4},\;\forall\;j\in\nset{l}}$ and $k\in\nsetro{l}$ holds\\
    $\erWe[k]\Ind{\setB{\VnormLp{7\hxdfPr[l]\dxdfPr[l]}^2<\peneSv[l]}}\Ind{\aixEv[l]}$\\
    \null\hfill
    $\leq\exp\big(\rWn\big\{\tfrac{25}{2}\penSv[l]+\tfrac{1}{8}\VnormLp{\ProjC[0]\xdf}^2\bias[l]^2(\xdf)-\tfrac{1}{8}\VnormLp{\ProjC[0]\xdf}^2\bias^2(\xdf)-\tfrac{1}{50}\penSv\big\}\big)$;
  \item\label{re:erWe:ii} For
    $\VnormLp{\ProjC[l]\pxdfPr[\ssY]}^2=\sum_{|j|\in\nsetlo{l,\ssY}}\iSv[j]^{-1}\eiSv[j]|\fxdf[j]|^2$ and $\Di\in\nsetlo{l,\ssY}$ holds\\
    $\erWe\Ind{\setB{7\VnormLp{\hxdfPr-\dxdfPr}^2<\peneSv}} \leq
   \exp\big(\rWn\big\{-\tfrac{1}{2}\peneSv+\tfrac{3}{2}\VnormLp{\ProjC[l]\pxdfPr[\ssY]}^2+\peneSv[l]\big\}\big)$.
  \end{resListeN}
\end{lem}

\begin{pro}[Proof of \cref{re:erWe}.]
Given $\Di,l\in\nset{\ssY}$ and an event $\dmEv{\Di}{l}$ (to be specified below) it follows
\begin{equation}\label{p:re:erWe:e1}
 \erWe\Ind{\dmEv{\Di}{l}}
\leq
\exp\big(\rWn\big\{\VnormLp{\hxdfPr}^2-\VnormLp{\hxdfPr[l]}^2+(\peneSv[l]-\peneSv)\big\}\big)\Ind{\dmEv{\Di}{l}}.
\end{equation}
We distinguish the two cases \ref{re:erWe:i} $\Di\in\nsetro{1,l}$ and
\ref{re:erWe:ii} $\Di\in\nsetlo{l,\ssY}$.
Consider first  \ref{re:erWe:i} $\Di\in\nsetro{1,l}$. From \ref{re:contr:e1} in \cref{re:contr}
(with
$\dxdfPr[]:=\hxdfPr[\ssY]$
and  $\pxdfPr[]:=\pxdfPr[\ssY]$) follows
\begin{equation}\label{p:re:erWe:e2}
 \erWe\Ind{\dmEv{\Di}{l}}
\leq
\exp\big(\rWn\big\{\tfrac{11}{2}\VnormLp{\hxdfPr[l]-\pxdfPr[l]}^2-\tfrac{1}{2}\VnormLp{\dProj{\Di}{l}\pxdfPr[\ssY]}^2
+(\peneSv[l]-\peneSv)\big\}\big)\Ind{\dmEv{k}{l}}
\end{equation}
Setting $\dmEv{k}{l}:=\{\VnormLp{\hxdfPr[l]-\pxdfPr[l]}^2<\peneSv[l]/7\}\cap\aixEv[l]$
the last bound together with \cref{re:aixEv} \ref{re:aixEv:i}
and \ref{re:aixEv:iv} implies the  assertion \ref{re:erWe:i}.
Consider secondly \ref{re:erWe:ii} $\Di\in\nsetlo{l,\ssY}$. From
\ref{re:contr:e2} in \cref{re:contr} (with $\dxdfPr[]:=\hxdfPr[\ssY]$
and  $\pxdfPr[]:=\pxdfPr[\ssY]$) and \eqref{p:re:erWe:e1} follows
 \begin{equation*}
  \erWe\Ind{\dmEv{l}{\Di}}
\leq \exp\big(\rWn\big\{\tfrac{7}{2}\VnormLp{\hxdfPr[k]-\pxdfPr[k]}^2+\tfrac{3}{2}\VnormLp{\dProj{l}{k}\pxdfPr[\ssY]}^2+(\peneSv[l]-\peneSv)\big\}\big)\Ind{\dmEv{l}{k}}
\end{equation*}
Setting $\dmEv{l}{\Di}:=\{\VnormLp{\hxdfPr-\pxdfPr}^2<\peneSv/7\}$ the
last bound together with \cref{re:aixEv} \ref{re:aixEv:i} implies
\ref{re:erWe:ii},
 which completes the proof.\proEnd
\end{pro}

\begin{lem}\label{au:re:SrWe:ag}Consider \pcw $\erWe[]$
  as in \eqref{au:de:erWe} and penalties $(\peneSv)_{\Di\in\nset{\ssY}}$ as in \eqref{au:de:peneSv}.
 For any $\mdDi,\pdDi\in\nset{\ssY}$ and
  associated $\pDi,\mDi\in\nset{\ssY}$ as in \eqref{au:de:*Di:ag}
  hold
  \begin{resListeN}
  \item\label{au:re:SrWe:ag:i}
    $\FuVg{\rWe[]}(\nsetro{\mDi})\leq
    \tfrac{50}{\rWc\cpen}\Ind{\{\mDi>1\}}
    \exp\big(-\tfrac{\rWc\cpen}{2} \DipenSv[\mdDi]\big)
    +\Ind{\{\VnormLp{\hxdfPr[\mdDi]-\pxdfPr[\mdDi]}^2
      \geq\peneSv[\mdDi]/7\}\cup\aixEv[\mdDi]^c}$;
  \item\label{au:re:SrWe:ag:ii}
    $\sum_{\Di\in\nsetlo{\pDi,n}}\peneSv\erWe
    \Ind{\{\VnormLp{\hxdfPr-\pxdfPr}^2<\peneSv/7\}}
    \leq \tfrac{16}{\rWc}\ssY^{-1}$.
  \end{resListeN}
\end{lem}

\begin{pro}[Proof of \cref{au:re:SrWe:ag}.]
  Consider \ref{au:re:SrWe:ag:i}. Let $\mDi\in\nset{\mdDi}$ as in
  \eqref{au:de:*Di:ag}. For the non trivial case $\mDi>1$ from
  \cref{re:erWe} \ref{re:erWe:i} with $l=\mdDi$ follows for all
  $\Di<\mDi\leq \mdDi$
  \begin{multline*}
    \erWe\Ind{\setB{\VnormLp{\hxdfPr[\mdDi]-\pxdfPr[\mdDi]}^2
        <\peneSv[\mdDi]/7}\cap\aixEv[l]}
    \\\leq
    \exp\big(\rWn\big\{-\tfrac{1}{8}\VnormLp{\ProjC[0]\xdf}^2\bias^2(\xdf)
    +(\tfrac{25}{2}\penSv[\mdDi]
    +\tfrac{1}{8}\VnormLp{\ProjC[0]\xdf}^2\bias[\mdDi]^2(\xdf))
    -\tfrac{1}{50}\penSv\big\}\big),
  \end{multline*}
  and hence by exploiting the definition \eqref{au:de:*Di:ag} of
  $\mDi$, that is
  $\VnormLp{\ProjC[0]\xdf}^2\sbF\geq
  \VnormLp{\ProjC[0]\xdf}^2\sbF[(\mDi-1)]>
  \VnormLp{\ProjC[0]\xdf}^2\sbF[\mdDi]+104\penSv[\mdDi]$, we
  obtain for each $\Di\in\nsetro{\mDi}$
  \begin{equation*}
    \erWe\Ind{\setB{\VnormLp{\hxdfPr[\mdDi]-\pxdfPr[\mdDi]}^2
        <\peneSv[\mdDi]/7}\cap\aixEv[l]}
    \leq\exp\big(-\tfrac{1}{2}\rWn\penSv[\mdDi]-\tfrac{1}{50}\rWn\penSv\big).
  \end{equation*}
  The last upper bound together with
  $\penSv=\cpen \DipenSv \ssY^{-1}\geq \cpen\Di\ssY^{-1}$,
  $\Di\in\nset{\ssY}$, as in \eqref{ak:de:LiSy} gives
  \begin{equation*}
    \FuVg{\rWe[]}(\nsetro{\mDi})
    \leq\exp\big(-\tfrac{\rWc}{2}\ssY\penSv[\mdDi]\big)
    \sum_{\Di\in\nsetro{\mDi}}\exp(-\tfrac{\rWc\cpen}{50}\Di)
    +\Ind{\setB{\VnormLp{\hxdfPr[\mdDi]-\pxdfPr[\mdDi]}^2
        \geq\peneSv[\mdDi]/7}\cup\aixEv[\mdDi]^c}
  \end{equation*}
  which combined with $\sum_{\Di\in\Nz}\exp(-\mu\Di)\leq \mu^{-1}$ for any $\mu>0$
  implies \ref{au:re:SrWe:ag:i}.
  Consider \ref{au:re:SrWe:ag:ii}. Let $\pDi\in\nset{\pdDi,\ssY}$ as
  in \eqref{au:de:*Di:ag}. For the non trivial case $\pDi<\ssY$ from
  \cref{re:erWe} \ref{re:erWe:ii} with $l=\pdDi$ follows for all
  $\Di>\pDi\geq \pdDi$
  \begin{equation*}
    \erWe\Ind{\setB{\VnormLp{\hxdfPr-\pxdfPr}^2<\peneSv/7}}
    \leq\exp\big(\rWn\big\{-\tfrac{1}{2}\peneSv
    +\tfrac{3}{2}\VnormLp{\ProjC[\pdDi]\pxdfPr[\ssY]}^2
    +\peneSv[\pdDi]\big\}\big),
  \end{equation*}
  and hence by employing the definition \eqref{au:de:*Di:ag} of $\pDi$, that is,
  $\tfrac{1}{4}\peneSv\geq \tfrac{1}{4}\peneSv[(\pDi+1)] >
  \peneSv[\pdDi]$ $+\tfrac{3}{2}\VnormLp{\ProjC[\pdDi]\pxdfPr[\ssY]}^2$, we
  obtain for each $\Di\in\nsetlo{\pDi,\ssY}$
  \begin{equation*}
    \erWe\Ind{\setB{\VnormLp{\txdfPr-\pxdfPr}^2<\peneSv/7}}
      \leq \exp\big(\rWn\big\{-\tfrac{1}{4} \peneSv\big\}\big).
  \end{equation*}
  The last bound together with \cref{re:aixEv} \ref{re:aixEv:iv},
  i.e., $\peneSv=\peneSv\Ind{\{\meiSv\geq1\}}$, implies
  \begin{multline}\label{p:au:re:SrWe:ag:e1}
    \sum_{\Di\in\nsetlo{\pDi,\ssY}}\peneSv\erWe
    \Ind{\{\VnormLp{\hxdfPr[k]-\pxdfPr[k]}^2<\peneSv/7\}}
    \leq \sum_{\Di\in\nsetlo{\pDi,\ssY}}\peneSv
    \exp\big(-\tfrac{\rWc}{4}\ssY\peneSv\big)
    \\\hfill
    =\sum_{\Di\in\nsetlo{\pDi,\ssY}}\peneSv
    \exp\big(-\tfrac{\rWc}{4}\ssY\peneSv\big)
    \Ind{\{\meiSv\geq1\}}\\
    =\cpen\ssY^{-1}\sum_{\Di\in\nsetlo{\pDi,\ssY}}
    \Di\cmeiSv\meiSv\exp\big(-\tfrac{\rWc\cpen}{4}\Di\cmeiSv\meiSv\big)
    \Ind{\{\meiSv\geq1\}}
  \end{multline}
  Comparing the last bound with \eqref{ak:re:SrWe:ag:pe1} the
  remainder of the proof of \ref{au:re:SrWe:ag:ii} follows line by
  line the arguments used to prove of \cref{ak:re:SrWe:ag}
  \ref{ak:re:SrWe:ag:ii} starting by \eqref{ak:re:SrWe:ag:pe1}, and we
  omit the details, which completes the proof.\proEnd
\end{pro}

\begin{lem}\label{au:re:SrWe:ms}Consider \msw $\msWe[]$
  as in \eqref{au:de:msWe} and penalties  $(\peneSv)_{\Di\in\nset{\ssY}}$ as in \eqref{au:de:peneSv}.  For any $\mdDi,\pdDi\in\nset{\ssY}$ and associated
  $\pDi,\mDi\in\nset{\ssY}$ as in \eqref{au:de:*Di:ag} hold
  \begin{resListeN}[]
  \item\label{au:re:SrWe:ms:i}
    $\FuVg{\msWe[]}(\nsetro{\mDi})\Ind{\{\VnormLp{\hxdfPr[\mdDi]-\pxdfPr[\mdDi]}^2<\peneSv[\mdDi]/7\}\cap\aixEv[\mdDi]}=0$;
  \item\label{au:re:SrWe:ms:ii}
    $\sum_{\Di\in\nsetlo{\pDi,\ssY}}\peneSv\msWe\Ind{\{\VnormLp{\hxdfPr-\pxdfPr}^2<\peneSv/7\}}=0$.
  \end{resListeN}
\end{lem}

\begin{pro}[Proof of \cref{au:re:SrWe:ms}.]  The assertions can be directly deduced from \cref{au:re:SrWe:ag} by letting  $\rWc\to\infty$ or  following line by line the proof of
  \cref{ak:re:SrWe:ms}, and we omit the details.\proEnd
\end{pro}

\begin{lem}\label{au:re:nd:rest}Consider $(\peneSv)_{\Di\in\nset{1,\ssY}}$ as in \eqref{au:de:peneSv} with $\cpen\geq84$. Let $\Di_{\ydf}:=\floor{3(400\Vnormlp[1]{\fydf})^2}$ and
  $ \ssY_{o}:=15(600)^4$.  There exists a finite numerical constant
  $\cst{}>0$ such that for all $\ssY\in\Nz$ and all
  $\mdDi\in\nset{\ssY}$ hold
  \begin{resListeN}
  \item\label{au:re:nd:rest1}
    $\sum_{\Di\in\nset{\ssY}}\FuEx[\ssY]{\rY|\rE}\vectp{\VnormLp{\hxdfPr-\dxdfPr}^2-\peneSv/7}\leq
    \cst{}\ssY^{-1}\big([1\vee\meiSv[\Di_{\ydf}]]\Di_{\ydf}+[1\vee\meiSv[\ssY_{o}]]\big)$;
  \item\label{au:re:nd:rest2}
    $\sum_{\Di\in\nset{\ssY}}\peneSv\FuVg[\ssY]{\rY|\rE}\big(\VnormLp{\hxdfPr-\dxdfPr}^2\geq\peneSv/7\big)\leq\cst{}\ssY^{-1}\big([1\vee\meiSv[\Di_{\ydf}]^2]\Di_{\ydf}^3+[1\vee\meiSv[\ssY_{o}]^2]\big)$;
  \item\label{au:re:nd:rest3}
    $\FuVg[\ssY]{\rY|\rE}\big(\VnormLp{\hxdfPr[\mdDi]-\dxdfPr[\mdDi]}^2\geq\peneSv[\mdDi]/7\big)\leq
    \cst{}
    \big(\exp\big(\tfrac{-\cmeiSv[\mdDi]\mdDi}{200\Vnormlp[1]{\fydf}}\big)
    +\ssY^{-1}\big)$.
  \end{resListeN}
\end{lem}

\begin{pro}[Proof of \cref{au:re:nd:rest}.]
  By using \cref{re:cconc} rather than \cref{re:conc} together with
  $\setB{\meiSv[l]<1}=\setB{\meiSv[l]=0}$ for all $l\in\Nz$ due to
  \cref{re:aixEv} \ref{re:aixEv:ii} the proof follows line by line the proof of \cref{ak:re:nd:rest}, and we omit the details.\proEnd
\end{pro}

\begin{lem}\label{au:ag:ub:p}Let the assumptions of
  \cref{au:ag:ub:pnp} \ref{au:ag:ub:pnp:p} be satisfied.
    There is a numerical constant $\cst{}$
    such that for all $\ssY,\ssE\in\Nz$ with $\ssY_{o}:=15(600)^4$ holds
    \begin{multline}\label{au:ag:ub:p:e1}
    \nmEx\VnormLp{\hxdf[{\We[]}]-\xdf}^2\leq
    \cst{}\VnormLp{\ProjC[0]\xdf}^2\big[\ssY^{-1}\vee \ssE^{-1}
    \vee\exp\big(\tfrac{-\cmiSv[{[\sDi{\ssY}\wedge\sDi{\ssE}]}]
      [\sDi{\ssY}\wedge\sDi{\ssE}]}{\Di_{\ydf}}\big)\big]\\
   \hfill + \cst{}\big(\big[1\vee
    K\vee c_{\xdf}K^2\miSv[K]^2\big]
    \big(\miSv[1]^2+\VnormLp{\ProjC[0]\xdf}^2\big)
    +\miSv[\Di_{\ydf}]^2\Di_{\ydf}^4
    +\miSv[\ssY_{o}]^2\big)\ssY^{-1}\\
    + \cst{}\big(\miSv[1]^2
    +K\miSv[K]^2
    +\VnormLp{\ProjC[0]\xdf}^2\miSv[K] \big)\ssE^{-1}.
  \end{multline}\end{lem}

\begin{pro}[Proof of \cref{au:ag:ub:p}.]The proof follows a long the
  lines of the proof of  \cref{ak:ag:ub:p} by using the upper bound
  \eqref{au:ag:ub:e1} instead of \eqref{ak:ag:ub:p3}  which
  hold for any $\mdDi,\pdDi\in\nset{\ssY}$ and associated
  $\mDi,\pDi\in\nset{\ssY}$ as defined in \eqref{au:de:*Di:ag} contrarily to \eqref{ak:de:*Di:ag}. We
  present exemplary the case \ref{ak:ag:ub:pnp:p:c2} $\ssY>
\ssY_{\xdf}:=[K\vee\floor{c_{\xdf}\DipenSv[K]}]$ with $K\in\Nz$ and $c_{\xdf}:=\tfrac{104\cpen}{\VnormLp{\ProjC[0]\xdf}^2\sbF[(K-1)]}$,
  and omit the details for the others.
Setting
$\pdDi:=K\leq
\ssY_{\xdf}$, i.e.,
$\pdDi\in\nset{\ssY}$, it follows $\bF[\pdDi]=0$ and
$\penSv[\pdDi]=\cpen\DipenSv[K]\ssY^{-1}\leq
\cpen K^2\miSv[K]^2\ssY^{-1}$. From
\eqref{au:ag:ub:e1} follows for all $\ssY> \ssY_{\xdf}$ thus
\begin{multline*}
   \nmEx\VnormLp{\hxdf[{\erWe[]}]-\xdf}^2\leq
   3 \VnormLp{\ProjC[0]\xdf}^2\bias[\mDi]^2(\xdf)\\\hfill
     +\cst{} \VnormLp{\ProjC[0]\xdf}^2\Ind{\{\mDi>1\}} \big(
    \exp\big(\tfrac{-\cmiSv[\mdDi]\mdDi}{\Di_{\ydf}}\big)
    + \mVg(\aixEv[\mdDi]^c)\big)+ \cst{}\ssE\mVg(\aixEv[K]^c)\\\hfill
    +\cst{}\moRa
    +\cst{}\ssY^{-1}\big(K^2\miSv[K]^2+\miSv[\Di_{\ydf}]^2\Di_{\ydf}^3
    +\miSv[\ssY_{o}]^2+\VnormLp{\ProjC[0]\xdf}^2\Ind{\{\mDi>1\}}\big).
  \end{multline*}
  Exploiting \cref{re:evrest} \ref{re:evrest:ii} there is a numerical
  constant $\cst{}$ such that for all $\ssE\in\Nz$ holds
  $\mVg(\aixEv[K]^c)\leq\cst{}K\miSv[K]^2\ssE^{-2}$, which together
  with $\moRa\leq \VnormLp{\ProjC[0]\xdf}^2\miSv[K]\ssE^{-1}$ implies
  \begin{multline}\label{au:ag:ub:pnp:p5}
    \nmEx\VnormLp{\hxdf[{\erWe[]}]-\xdf}^2\leq
    +\cst{}\ssY^{-1}\big(K^2\miSv[K]^2+\miSv[\Di_{\ydf}]^2\Di_{\ydf}^3
    +\miSv[\ssY_{o}]^2+\VnormLp{\ProjC[0]\xdf}^2\Ind{\{\mDi>1\}}\big).\\\hfill
    + 3 \VnormLp{\ProjC[0]\xdf}^2\sbF[\mDi]
    +\cst{} \VnormLp{\ProjC[0]\xdf}^2\Ind{\{\mDi>1\}} \big(
    \exp\big(\tfrac{-\cmiSv[\mdDi]\mdDi}{\Di_{\ydf}}\big)
    + \mVg(\aixEv[\mdDi]^c)\big)\\\hfill
    +\cst{}\ssE^{-1}\big(K\miSv[K]^2+ \VnormLp{\ProjC[0]\xdf}^2\miSv[K]\big)
  \end{multline}
  In order to control the terms involving $\mdDi$ and $\mDi$ we
  destinguish for $\ssE\in\Nz$ with
  $\ssE_{\xdf,\iSv}:=\floor{289\log(K+2)\cmiSv[K]\miSv[K]}$ the
  following two cases,
  \begin{inparaenum}[i]\renewcommand{\theenumi}{\dgrau\rm(b-\roman{enumi})}
  \item\label{au:ag:ub:pnp:p:m1}
    $\ssE\in\nset{\ssE_{\xdf,\iSv}}$ and
  \item\label{au:ag:ub:pnp:p:m2}
    $\ssE>\ssE_{\xdf,\iSv}$.
  \end{inparaenum}
  Consider first \ref{au:ag:ub:pnp:p:m1}
  $\ssE\in\nset{\ssE_{\xdf,\iSv}}$.  We set $\mdDi=1$ and hence
  $\mDi=1$. Thereby, with $\bias[1]^2(\xdf)\leq1$,
  $\log(K+2)\leq \tfrac{K+2}{e}\leq 2K$,
  $\cmiSv[\Di]\miSv[\Di]\leq K\miSv[K]^2$, and hence
  $\ssE_{\xdf,\iSv}\leq\cst{}K^2\miSv[K]^2$, from
  \eqref{au:ag:ub:pnp:p5} follows for all
  $\ssE\in\nset{\ssE_{\xdf,\iSv}}$
  \begin{multline}\label{au:ag:ub:pnp:p6}
    \nmEx\VnormLp{\hxdf[{\erWe[]}]-\xdf}^2\leq
    \cst{}\ssY^{-1}\big(K^2\miSv[K]^2+\miSv[\Di_{\ydf}]^2\Di_{\ydf}^3
    +\miSv[\ssY_{o}]^2\big)\\\hfill
    +\cst{}\ssE^{-1}\big(K\miSv[K]^2
    + \VnormLp{\ProjC[0]\xdf}^2(K^2\miSv[K]^2+\miSv[K])\big)
  \end{multline}
 Consider \ref{au:ag:ub:pnp:p:m2} $\ssE>\ssE_{\xdf,\iSv}$ ensuring the defining set of
  $\sDi{\ssE}=\max\{\Di\in\nset[1,]{\ssE}:289\log(\Di+2)\cmiSv[\Di]\miSv[\Di]\leq\ssE\}$
  is not empty and $\sDi{\ssE}\geq K$. For each
  $\mdDi\in\nset{K,\sDi{\ssE}}$ it follows
  $\FuVg[\ssE]{\rE}(\aixEv[\mdDi]^c)\leq 53\ssE^{-1}$ due to
  \cref{re:evrest} \ref{re:evrest:iii}.
  Since
  $\ssY> \ssY_{\xdf}=[K\vee\floor{c_{\xdf}\DipenSv[K]}]$ with
  $c_{\xdf}:=\tfrac{104\cpen}{\VnormLp{\ProjC[0]\xdf}^2\sbF[(K-1)]}$
  the defining set of
  $\sDi{\ssY}=\max\{\Di\in\nset[1,]{\ssY}:\ssY>c_{\xdf}\DipenSv\}$ is
  not empty and $\sDi{\ssY}\geq K$. For each $\mdDi\in\nset{K,\sDi{\ssY}}$ we have $\bias[\mdDi](\xdf)=0$, and
  $\penSv[\mdDi]=\DipenSv[\mdDi]\ssY^{-1}<c_{\xdf}^{-1}
  =\tfrac{\VnormLp{\ProjC[0]\xdf}^2\sbF[(K-1)]}{104\cpen}$. It follows
    $\VnormLp{\ProjC[0]\xdf}^2\sbF[(K-1)]>
  \VnormLp{\ProjC[0]\xdf}^2\sbF[{\mdDi}]+104\penSv[\mdDi]$
  and trivially
  $\VnormLp{\ProjC[0]\xdf}^2\sbF[{K}]=0
  < \VnormLp{\ProjC[0]\xdf}^2\sbF[{\mdDi}]+104\penSv[\mdDi]$. Therefore,
  $\mDi$ as in \eqref{au:de:*Di:ag} satisfies $\mDi=K$ and hence
  $\bF[\mDi]=0$. Finally, setting
  $\mdDi:=\sDi{\ssY}\wedge\sDi{\ssE}$ it follows
  $\mVg(\aixEv[\mdDi]^c)\leq 53\ssE^{-1}$, $\mDi=K$ and
  $\bF[\mDi]=0$. From \eqref{au:ag:ub:pnp:p5}
  follows for all $\ssE>\ssE_{\xdf,\iSv}$ and $\ssY>\ssY_{\xdf,\iSv}$ thus
  \begin{multline}\label{au:ag:ub:pnp:p7}
    \nmEx\VnormLp{\hxdf[{\erWe[]}]-\xdf}^2
    \leq \cst{}\ssY^{-1}\big(K^2\miSv[K]^2\ssY^{-1}
    +\miSv[\Di_{\ydf}]^2\Di_{\ydf}^3+\miSv[\ssY_{o}]^2
    +\VnormLp{\ProjC[0]\xdf}^2\big)\\\hfill
    +\cst{}\VnormLp{\ProjC[0]\xdf}^2
    \exp\big(\tfrac{-\cmiSv[\mdDi]\mdDi}{\Di_{\ydf}}\big)
    +\cst{}\ssE^{-1}\big(K\miSv[K]^2
    + \VnormLp{\ProjC[0]\xdf}^2\miSv[K]\big).
  \end{multline}
  By combining \eqref{au:ag:ub:pnp:p6} and \eqref{au:ag:ub:pnp:p7} for
  the cases \ref{au:ag:ub:pnp:p:m1} $\ssE\in\nset{\ssE_{\xdf,\iSv}}$
  and \ref{au:ag:ub:pnp:p:m2} $\ssE>\ssE_{\xdf,\iSv}$ the upper bound
  \eqref{au:ag:ub:p:e1} holds in case \ref{ak:ag:ub:pnp:p:c2}, i.e., for all $\ssE\in\Nz$ and for all
  $\ssY>\ssY_{\xdf,\iSv}$,
  which completes the proof of \cref{au:ag:ub:p}.\proEnd
\end{pro}
\subsubsection{Proof of \cref{au:mrb:ag:ub:pnp} and
  \cref{au:mrb:ag:ub2:pnp}}\label{a:au:mrb}

\begin{pro}[Proof of \cref{au:mrb:ag:ub:pnp}.]
  \begin{te}
    Keeping
  \eqref{a:ak:mrb:LiSy} in mind  for
all  $\xdf\in\rwCxdf$, $\edf\in\rwCedf$ and $\Di,\ssY,\ssE\in\Nz$ we have
$\Vnormlp[1]{\fydf}^2\leq  \xdfCr\edfCr\,\Vnormlp[1]{\edfCwS\xdfCwS}$, hence
$\Di_{\ydf}=\floor{3(400)^2\Vnormlp[1]{\fydf}^2}\leq
\floor{3(400)^2\xdfCr\zeta_{\edfCr}\,\Vnormlp[1]{\edfCwS\xdfCwS}}=\Di_{\xdfCw[]\edfCw[]}$
and $\tfrac{1}{\Di_{\ydf}}\liSv\geq\tfrac{1}{\Di_{\xdfCw[]\edfCw[]}}\liCw$,
  $\VnormLp{\ProjC[0]\xdf}^2\leq\xdfCr$,
$\VnormLp{\ProjC[0]\xdf}^2\sbFxdf\leq
\xdfCr\xdfCw$, $\penSv\leq \cpen\zeta_{\edfCr}\Di\LiCw/\ssY$,
$\moRa\leq\xdfCr\edfCr\mmRa$,  and
 $\mDi$ as in \eqref{au:de:*Di:ag} satisfies
   $\VnormLp{\ProjC[0]\xdf}^2\bias[\mDi]^2(\xdf) \leq
   \xdfCr\xdfCw[\mdDi]+104\cpen\zeta_{\edfCr}\mdDi\LiCw[\mdDi]/\ssY$. Combining
   the
   last bounds together with the upper bound \eqref{au:ag:ub:e1} there is a numerical constant $\cst{}>0$ such
  that uniformely for
all  $\xdf\in\rwCxdf$, $\edf\in\rwCedf$, $\ssY,\ssE\in\Nz$ and $\mdDi,\pdDi\in\nset{\ssY}$ holds
  \begin{multline}\label{au:mrb:ag:ub:e1}
   \nmEx\VnormLp{\hxdf[{\We[]}]-\xdf}^2\leq
    2\cpen\zeta_{\edfCr}\LiCw[\pdDi]\pdDi/\ssY
    +\tfrac{12}{7}\xdfCr\xdfCw[\pdDi]+3\xdfCr\xdfCw[\mdDi]+312\cpen\zeta_{\edfCr}\mdDi\LiCw[\mdDi]/\ssY\\\hfill
    +\cst{}\xdfCr\big(
    \exp\big(\tfrac{-\liCw[\mdDi]\mdDi}{\Di_{\xdfCw[]\edfCw[]}}\big)
    + \mVg(\aixEv[\mdDi]^c)\big)+ \cst{}\ssE\mVg(\aixEv[\pdDi]^c)
    \\\hfill
    +\cst{}\xdfCr\edfCr\mmRa
    +\cst{}\ssY^{-1}\big(\edfCr^2\edfCw[\Di_{\xdfCw[]\edfCw[]}]^2\Di_{\xdfCw[]\edfCw[]}^4+\edfCr^2\edfCw[\ssY_{o}]^2
    +\xdfCr\big)
  \end{multline}
  We destinguish for $\ssE\in\Nz$ with
  $\ssE_{\edfCw[]}:=\floor{289(\log3)\zeta_{\edfCr}\liCw[1]\edfCw[1]}$ the
  two cases,
  \begin{inparaenum}[i]\renewcommand{\theenumi}{\dgrau\rm(\alph{enumi})}
  \item\label{au:mrb:ag:ub:pnp:m1}
    $\ssE\in\nset{\ssE_{\edfCw[]}}$ and
  \item\label{au:mrb:ag:ub:pnp:m2}
    $\ssE>\ssE_{\edfCw[]}$.
  \end{inparaenum}
\end{te}

\begin{te}
  Consider \ref{au:mrb:ag:ub:pnp:m1}. We set $\pdDi=\mdDi=1$. Since
  $\FuVg[\ssE]{\rE}(\aixEv[1]^c)\leq
  \cst{}\miSv[1]^2\ssE^{-2}\leq \cst{}\edfCr^2\edfCw[1]^2\ssE^{-2}$   due to \cref{re:evrest}
  \ref{re:evrest:ii}, \eqref{au:mrb:ag:ub:e1} implies
  for all $\ssY\in\Nz$ and $\ssE\in\nset{\ssE_{\edfCw[]}}$
 \begin{multline}\label{au:mrb:ag:ub:pnp:p1}
   \mnmRi{\hxdfAg[{\We[]}]}{\rwCxdf}{\rwCedf}
   \leq\cst{}\xdfCr\edfCr\mmRa    \\\hfill
    +   \cst{}
   \ssE^{-1}(\xdfCr\zeta_{\edfCr}\xdfCw[1]+\edfCr^2)\edfCw[1]^2
    +\cst{}\ssY^{-1}\big(\edfCr^2\edfCw[\Di_{\xdfCw[]\edfCw[]}]^2\Di_{\xdfCw[]\edfCw[]}^4+\edfCr^2\edfCw[\ssY_{o}]^2
    +\xdfCr\big).
  \end{multline}
\end{te}

\begin{te}
  Consider secondly \ref{au:ag:ub:pnp:np:m2}. Since $\ssE>\ssE_{\edfCw[]}$
  the defining set of
  $\sDi{\ssE}:=\max\{\Di\in\nset{\ssE}:289\log(\Di+2)\zeta_{\edfCr}\liCw[\Di]\edfCw[\Di]\leq\ssE\}$
  is not empty. Keeping in mind, that due to \eqref{a:ak:mrb:LiSy}
   for all
  $\edf\in\rwCedf$ and for each
  $\Di\in\nset{\sDi{\ssE}}$ holds
  $\zeta_{\edfCr}\liCw\edfCw\geq \liSv\miSv$, and hence
  $\ssE\geq289\log(\Di+2)\cmiSv\miSv$ and
$\FuVg[\ssE]{\rE}(\aixEv[\Di]^c)\leq 11226\ssE^{-2}$ applying \cref{re:evrest}
  \ref{re:evrest:iii}.
   For
  $\oDi{\ssY}:=\nmDiL\in\nset{\ssY}$ as in \eqref{oo:de:doRao} let
  $\pdDi:=\oDi{\ssY}\wedge\sDi{\ssE}$ and hence
  $\ssE\FuVg[\ssE]{\rE}(\aixEv[\pdDi]^c)\leq \cst{}\ssE^{-1}$.  Since
  $\LiCw[\pdDi]\pdDi/\ssY\leq\dmRaL{\oDi{\ssY}}=\nmRaL$ and $\xdfCw[\pdDi]\leq \nmRaL+\xdfCw[\sDi{\ssE}]$
  from \eqref{au:mrb:ag:ub:e1} follows
    \begin{multline}\label{au:mrb:ag:ub:pnp:p2}
   \mnmRi{\hxdfAg[{\We[]}]}{\rwCxdf}{\rwCedf}\leq
    (2\cpen\zeta_{\edfCr}+2\xdfCr) \nmRaL +\tfrac{12}{7}\xdfCr\xdfCw[\sDi{\ssE}]\\\hfill+3\xdfCr\xdfCw[\mdDi]+312\cpen\zeta_{\edfCr}\mdDi\LiCw[\mdDi]/\ssY
    +\cst{}\xdfCr\big(
    \exp\big(\tfrac{-\liCw[\mdDi]\mdDi}{\Di_{\xdfCw[]\edfCw[]}}\big)
    + \mVg(\aixEv[\mdDi]^c)\big)  \\\hfill + \cst{}\ssE^{-1}
    +\cst{}\xdfCr\edfCr\mmRa
    +\cst{}\ssY^{-1}\big(\edfCr^2\edfCw[\Di_{\xdfCw[]\edfCw[]}]^2\Di_{\xdfCw[]\edfCw[]}^4+\edfCr^2\edfCw[\ssY_{o}]^2
    +\xdfCr\big)
  \end{multline}
  For
  $\sDi{\ssY}:=\argmin\{\dmRaL{\Di}\vee\exp\big(\tfrac{-\liCw[\Di]\Di}{\Di_{\xdfCw[]\edfCw[]}}\big):\Di\in\nset[1,]{\ssY}\}$
  with $\dmRaL{\sDi{\ssY}}\leq\nmRaA$  let
$\mdDi:=\sDi{\ssY}\wedge\sDi{\ssE}$ and hence
$\FuVg[\ssE]{\rE}(\aixEv[\mdDi]^c)\leq 53\ssE^{-1}$. Since $
\xdfCr\xdfCw[\mdDi]+\zeta_{\edfCr}\LiCw[\mdDi]\mdDi\ssY^{-1}
\leq \xdfCr\xdfCw[\sDi{\ssE}]+
(\xdfCr+\zeta_{\edfCr})\dmRaL{\sDi{\ssY}}$ and
 $\ssY^{-1}\leq\nmRaL\leq\dmRaL{\sDi{\ssY}}$ from
 \eqref{au:mrb:ag:ub:pnp:p2} follows
 for all $\ssY\in\Nz$, $\ssE>\ssE_{\edfCw[]}$
    \begin{multline}\label{au:mrb:ag:ub:pnp:p3}
   \mnmRi{\hxdfAg[{\We[]}]}{\rwCxdf}{\rwCedf}\leq
   \cst{}(\xdfCr+\zeta_{\edfCr})\nmRaA +\cst{}\xdfCr\edfCr\mmRa\\\hfill  +5\xdfCr\big[\xdfCw[\sDi{\ssE}]^2\vee
    \exp\big(\tfrac{-\liCw[\sDi{\ssE}]\sDi{\ssE}}{\Di_{\xdfCw[]\edfCw[]}}\big)\big] + \cst{}\xdfCr\ssE^{-1}
    +\cst{}\ssY^{-1}\big(\edfCr^2\edfCw[\Di_{\xdfCw[]\edfCw[]}]^2\Di_{\xdfCw[]\edfCw[]}^4+\edfCr^2\edfCw[\ssY_{o}]^2 \big)
  \end{multline}
\end{te}

\begin{te}
Combining \eqref{au:mrb:ag:ub:pnp:p1} and \eqref{au:mrb:ag:ub:pnp:p3}
for   the cases  \ref{au:mrb:ag:ub:pnp:m1}
 and \ref{au:mrb:ag:ub:pnp:m2} for all $\ssY,\ssE\in\Nz$ holds
    \begin{multline}\label{au:mrb:ag:ub:pnp:p4}
   \mnmRi{\hxdfAg[{\We[]}]}{\rwCxdf}{\rwCedf}\leq
\cst{}(\xdfCr+\zeta_{\edfCr})\nmRaA\Ind{\{\ssE>\ssE_{\edf}\}}\\\hfill
+ \cst{}\xdfCr\big[\xdfCw[\sDi{\ssE}]^2\vee
    \exp\big(\tfrac{-\liCw[\sDi{\ssE}]\sDi{\ssE}}{\Di_{\xdfCw[]\edfCw[]}}\big)\big]\Ind{\{\ssE>\ssE_{\edf}\}} + \cst{}\xdfCr\edfCr\mmRa \\\hfill
    +   \cst{}\ssE^{-1}(\xdfCr\zeta_{\edfCr}\xdfCw[1]+\edfCr^2)\edfCw[1]^2
    +\cst{}\ssY^{-1}\big(\edfCr^2\edfCw[\Di_{\xdfCw[]\edfCw[]}]^2\Di_{\xdfCw[]\edfCw[]}^4+\edfCr^2\edfCw[\ssY_{o}]^2\big)
  \end{multline}
which  shows \eqref{au:mrb:ag:ub:pnp:e1} and completes the
proof of \cref{au:mrb:ag:ub:pnp}.
\end{te}
  \proEnd
\end{pro}

\begin{pro}[Proof of \cref{au:mrb:ag:ub2:pnp}.]
  The proof is similar to the proof of  \cref{ak:ag:ub2:pnp:mm} and
  \cref{au:ag:ub2:pnp}, and we omit the details.
  \proEnd
\end{pro}

\bibliography{DACD}
\end{document}